\documentstyle[twoside,psfig]{article}

\oddsidemargin=\evensidemargin 
\catcode`\@=11 \long\def\@makefntext#1{ \protect\noindent \hbox to
3.2pt {\hskip-.9pt
$^{{\eightrm\@thefnmark}}$\hfil}#1\hfill}       %CAN BE USED

\def\ps@myheadings{\let\@mkboth\@gobbletwo      %SIZE OF R/H NOS.
\def\@oddhead{\hbox{}
\rightmark\hfil\eightrm\thepage}
\def\@oddfoot{}\def\@evenhead{\eightrm\thepage\hfil
\leftmark\hbox{}}\def\@evenfoot{}
\def\sectionmark##1{}\def\subsectionmark##1{}}

\catcode`@=11                        %SIZE OF OPENING PAGE NO.
\def\ps@plain{\let\@mkboth\@gobbletwo
     \def\@oddhead{}\def\@oddfoot{\eightrm\hfil\thepage
     \hfil}\def\@evenhead{}\let\@evenfoot\@oddfoot}

\newcounter{sectionc}\newcounter{subsectionc}\newcounter{subsubsectionc}
\renewcommand{\section}[1] {\vspace{12pt}\addtocounter{sectionc}{1}
\setcounter{subsectionc}{0}\setcounter{subsubsectionc}{0}\noindent
    {\tenbf\thesectionc. #1}\par\vspace{5pt}}
\renewcommand{\subsection}[1] {\vspace{12pt}\addtocounter{subsectionc}{1}
    \setcounter{subsubsectionc}{0}\noindent
    {\bf\thesectionc.\thesubsectionc.
    {\kern1pt \bfit #1}}\par\vspace{5pt}}
\renewcommand{\subsubsection}[1] {\vspace{12pt}
    \addtocounter{subsubsectionc}{1}
    \noindent
    {\tenrm\thesectionc.\thesubsectionc.\thesubsubsectionc. {\kern1pt
    \it #1}}\par\vspace{5pt}}

\newcounter{appendixc}
\newcounter{subappendixc}[appendixc]
\newcounter{subsubappendixc}[subappendixc]

\renewcommand{\appendix}[1] {\vspace{12pt}  %Appendix A. Heading +
    \refstepcounter{appendixc}      %Appendix A (11/5/94)
    \setcounter{figure}{0}
    \setcounter{table}{0}
    \setcounter{lemma}{0}
    \setcounter{theorem}{0}
    \setcounter{corollary}{0}
    \setcounter{definition}{0}
    \setcounter{equation}{0}
    \renewcommand{\thefigure}{\Alph{appendixc}.\arabic{figure}}
    \renewcommand{\thetable}{\Alph{appendixc}.\arabic{table}}
    \renewcommand{\theappendixc}{\Alph{appendixc}}
    \renewcommand{\thelemma}{\Alph{appendixc}.\arabic{lemma}}
    \renewcommand{\thetheorem}{\Alph{appendixc}.\arabic{theorem}}
    \renewcommand{\thedefinition}{\Alph{appendixc}.\arabic{definition}}
    \renewcommand{\thecorollary}{\Alph{appendixc}.\arabic{corollary}}
    \renewcommand{\theequation}{\Alph{appendixc}.\arabic{equation}}
    \noindent{\tenbf Appendix \theappendixc #1}\par\vspace{5pt}}

\newcommand{\smalllineskip}{\baselineskip=10pt}

\newcommand{\copyrightheading}[1]
    {\vspace*{-2.5cm}\smalllineskip{\flushleft
    {\footnotesize }\\
    {\footnotesize \copyright\kern2pt }\\
         }}

\newcounter{itemlistc}
\newcounter{romanlistc}
\newcounter{alphlistc}
\newcounter{arabiclistc}

\newcommand{\fcaption}[1]{
        \refstepcounter{figure}
        \setbox\@tempboxa = \hbox{\footnotesize Fig.~\thefigure. #1}
        \ifdim \wd\@tempboxa > 5in
           {\begin{center}
        \parbox{5in}{\footnotesize\smalllineskip Fig.~\thefigure. #1}
            \end{center}}
        \else
             {\begin{center}
             {\footnotesize Fig.~\thefigure. #1}
              \end{center}}
        \fi}

\newcommand{\tcaption}[1]{
        \refstepcounter{table}
        \setbox\@tempboxa = \hbox{\footnotesize Table~\thetable. #1}
        \ifdim \wd\@tempboxa > 5in
           {\begin{center}
        \parbox{5in}{\footnotesize\smalllineskip Table~\thetable. #1}
            \end{center}}
        \else
             {\begin{center}
             {\footnotesize Table~\thetable. #1}
              \end{center}}
        \fi}

\def\pmb#1{\setbox0=\hbox{#1}
    \kern-.025em\copy0\kern-\wd0
    \kern.05em\copy0\kern-\wd0
    \kern-.025em\raise.0433em\box0}

\def\fnt#1#2{\footnotetext{\kern-.3em
    {$^{\mbox{\scriptsize #1}}$}{#2}}}

\def\fpage#1{\begingroup
\voffset=.3in
\thispagestyle{empty}\begin{table}[b]\centerline{\footnotesize #1}
    \end{table}\endgroup}

\def\runninghead#1#2{\pagestyle{myheadings}
\markboth{{\protect\footnotesize\it{\quad #1}}\hfill}
{\hfill{\protect\footnotesize\it{#2\quad}}}}

\font\tenrm=cmr10  \font\tenbf=cmbx10
\font\bfit=cmbxti10 at 10pt \font\ninerm=cmr9 
 \font\eightrm=cmr8

\newtheorem{definition}{Definition}

\def\@begintheorem#1#2{\trivlist    %6/9/94
    \item[\hskip\labelsep{\bf #1\ #2.}]}
\def\@opargbegintheorem#1#2#3{\trivlist
    \item[\hskip\labelsep{\bf #1\ #2\ (#3).}]}

        {\setcounter{itemlistc}{0}      %\NOINDENT\BULLET
     \begin{list}{$\bullet$}        %AUTOMATICALLY
    {\usecounter{itemlistc}         %UNLESS OTHERWISE SPECIFIED
     \leftmargin10pt           %IE.  0PT = MOVE OUTSIDE LEFT MARGIN
     \setlength{\parsep}{0pt}
     \setlength{\itemsep}{0pt}     %ADD EXTRA VSPACE BETWEEN ITEMS
    }}{\end{list}}
    {\setcounter{romanlistc}{0}     %\NOINDENT (i)
     \begin{list}{$($\roman{romanlistc}$)$} %AUTOMATICALLY
    {\usecounter{romanlistc}        %UNLESS OTHERWISE SPECIFIED
     \leftmargin18pt %(ii)      %IE. 21PT = MOVE MORE TO THE RIGHT-->
     \setlength{\parsep}{0pt}
     \setlength{\itemsep}{0pt}  %ADD EXTRA VSPACE BETWEEN ITEMS
     \settowidth{\labelwidth}{#1}
    }}{\end{list}}
    {\setcounter{enumii}{0}         %\NOINDENT (a)
     \begin{list}{$($\alph{enumii}$)$}  %AUTOMATICALLY
    {\usecounter{enumii}            %UNLESS OTHERWISE SPECIFIED
     \leftmargin18pt        %IE.  0PT = MOVE OUTSIDE LEFT MARGIN
     \setlength{\parsep}{0pt}
     \setlength{\itemsep}{0pt}  %ADD EXTRA VSPACE BETWEEN ITEMS
     \settowidth{\labelwidth}{#1}
    }}{\end{list}}
\def\qed{\hbox{${\vcenter{\vbox{            %HOLLOW SQUARE
   \hrule height 0.4pt\hbox{\vrule width 0.4pt height 6pt
   \kern5pt\vrule width 0.4pt}\hrule height 0.4pt}}}$}}

\def\theequation{\thesectionc.\arabic{equation}}  %FOR SETTING EQ.~(1.1)
\begin{document}

\runninghead{Tutte and Jones Polynomials of Link Families} {Tutte
and Jones Polynomials of Link Families}

%\normalsize\textlineskip
%\thispagestyle{empty}
\setcounter{page}{1}

\markboth{Slavik Jablan,  Ljiljana Radovi\' c, Radmila Sazdanovi\' c
}{}

%%%%%%%%%%%%%%%%%%%%% Publisher's Area please ignore %%%%%%%%%%%%%%
%\catchline{}{}{}{}{}
%%%%%%%%%%%%%%%%%%%%%%%%%%%%%%%%%%%%%%%%%%%%%%%%%%%%%%%%%%%%%%%%%%%

\fpage{1} \centerline{\bf TUTTE AND JONES POLYNOMIALS OF LINK
FAMILIES}
\bigskip

\centerline{\footnotesize SLAVIK JABLAN,  LJILJANA RADOVI\' C$^ \dag
$, RADMILA SAZDANOVI\' C$^{\dag \dag }$,}

\medskip

\centerline{\footnotesize\it The Mathematical Institute, Knez
Mihailova 36,}\centerline{\footnotesize\it P.O.Box 367, 11001
Belgrade,}\centerline{\footnotesize\it Serbia}
\centerline{\footnotesize\it sjablan@gmail.com}

\medskip

\centerline{\footnotesize\it Faculty of Mechanical
Engineering$^\dag$}\centerline{\footnotesize\it A.~Medvedeva 14 }
\centerline{\footnotesize\it 18 000 Ni\v s }
\centerline{\footnotesize\it Serbia} \centerline{\footnotesize\it
ljradovic@gmail.com}

\medskip

\centerline{\footnotesize\it The George Washington University$^{\dag
\dag }$}\centerline{\footnotesize\it Department of
Mathematics,}\centerline{\footnotesize\it Monroe Hall of
Government,} \centerline{\footnotesize\it 2115 G Street, NW,}
\centerline{\footnotesize\it Washington, D.C. 20052, }
\centerline{\footnotesize\it USA} \centerline{\footnotesize\it
radmilas@gmail.edu}

\bigskip

\begin{abstract}
\noindent This article contains general formulas for Tutte and Jones
polynomials for families of knots and links given in Conway notation
and "portraits of families"-- plots of zeroes of their corresponding
Jones polynomials.
\end{abstract}

\section{Introduction}

Knots and links (or shortly $KL$s) will be given in Conway notation
[Con, Ro, Cau, JaSa].

\begin{definition}
For a link or knot $L$ given in an unreduced\footnote{The Conway
notation is called unreduced if in symbols of polyhedral links
elementary tangles 1 in single vertices are not omitted.} Conway
notation $C(L)$ denote by $S$ a set of numbers in the Conway symbol
excluding numbers denoting basic polyhedron and zeros (determining
the position of tangles in the vertices of polyhedron) and let
$\tilde S=\{a_1,a_2, \ldots, a_k\}$ be a non-empty subset of $S$.
Family $F_{\tilde S}(L)$ of knots or links derived from $L$ consists
of all knots or links $L'$ whose Conway symbol is obtained by
substituting all $ a_i\neq \pm 1$, by $sgn(a_i) |a_i+k_{a_i}|$,
$|a_i+k_{a_i}| >1$, $k_{a_i} \in Z$. [JaSa].
\end{definition}

An infinite subset of a family is called {\it subfamily}. If all
$k_{a_i}$ are even integers, the number of components is preserved
within the corresponding subfamilies, i.e., adding full-twists
preserves the number of components inside the subfamilies.

\begin{definition}
A link given by Conway symbol containing only tangles $\pm 1$ and
$\pm 2$ is called a {\it source link}.
\end{definition}

A {\it graph} is defined as a pair $(V,E)$, where $V$ is the {\it
vertex set} and $E\subseteq V\times V$ the {\it edge set}. We
consider only {\it undirected} graphs, meaning $(x,y)$ is the same
as $(y,x)$. A {\it loop} is an edge $(x,x)$ between the same vertex,
and a bridge is an edge whose removal disconnects two or more
vertices (i.e. there is no longer a path between them) [ThoPe].

Two operations are essential to understanding the Tutte polynomial
definition. These are: {\it edge deletion} denoted by $G-e$, and
{\it edge contraction} $G/e$. The latter involves first deleting
$e$, and then merging its endpoints as follows:

\begin{figure}[th]
\centerline{\psfig{file=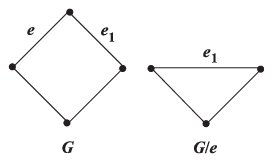,width=1.60in}}
\end{figure}

\begin{definition}
The {\it Tutte polynomial} of a graph $G(V,E)$ is a two-variable
polynomial defined as follows:

$$T(G)=\cases {1 & $E(\emptyset )$ \quad  \quad   \quad   \quad   \quad \quad  \quad   \quad   \quad   \quad  \quad  \quad   \quad   \quad   \quad \quad \,(1)  \cr
                xT(G/e) & $e\in E$ and $e$ is a bridge  \quad  \quad   \quad   \quad   \quad  \quad   \quad   \quad  \,(2)  \cr
                yT(G-e) &  $e\in E$ and $e$ is a loop \quad  \quad   \quad   \quad   \quad  \quad   \quad   \quad  \quad (3)  \cr
                T(G-e)+T(G/e)  &   $e\in E$ and $e$ is neither a loop or a bridge
                \,\,\,(4) \cr
                }$$
\end{definition}

The definition of a Tutte polynomial outlines a simple recursive
procedure for computing it, but the order in which rules are applied
is not fixed.

According to Thistlethwaite's Theorem, Jones polynomial of an
alternating link, up to a factor, can be obtained from Tutte
polynomial by replacements: $x\rightarrow -x$ and $y\rightarrow
-{1\over x}$ [Thi, Kau, Bo, ChaShro]. Moreover, from general
formulas for Tutte polynomials with negative values of parameters we
obtain Tutte polynomials expressed as Laurent polynomials. By the
same replacements we obtain, up to a factor, Jones polynomials of
non-alternating links.

A {\it cut-vertex} (or articulation vertex) of a connected graph is
a vertex whose removal disconnects the graph [Char]. In general, a
cut-vertex is a vertex of a graph whose removal increases the number
of components [Har]. A {\it block} is a maximal biconnected subgraph
of a given graph.

{\it One-point union} or {\it block sum} of two (disjoint) graphs
$G_1$ and $G_2$, neither of which is a vertex graph, and which we
shall denote as $G_1\ast G_2$ is of particular interest. This
one-point union is such that the intersection of $G_1$ and $G_2$ can
only consist of a vertex[KuMu].

Decomposition of a graph $G$ into a finite number of blocks
$G_1$,...,$G_n$, denoted by

$$G=G_1\ast G_2\ast \ldots G_n$$ is called the {\it block sum} of
$G_1$,...,$G_n$. The following formula holds for the Tutte
polynomial of the block sum:

$$T(G_1\ast G_2\ast \ldots G_n)=T(G_1)T(G_2)\ldots T(G_n).$$

A dual graph $\overline G$ of a given planar graph $G$ is a graph
which has a vertex for each plane region of $G$, and an edge for
each edge in $G$ joining two neighboring regions, for a certain
embedding of $G$. The Tutte polynomial of $\overline G$ can be
obtained from $T(G)$ by replacements $x\rightarrow y$, $y\rightarrow
x$, i.e. $T(\overline G)(x,y)=T(G)(y,x)$.

There is a nice bijective correspondence between $KL$s and graphs:
to obtain a graph from a projection of $KL$, first color every other
region of the $KL$ diagram black or white, so that the infinite
outermost region is black. In the {\it checker-board coloring} (or
{\it Tait coloring}) of the plane obtained, put a vertex at the
center of each white region. Two vertices of a graph are connected
by an edge if there was a crossing between  corresponding regions in
a $KL$ diagram. In addition, to each edge of a graph we can assign
the sign of its corresponding vertex of the $KL$ diagram. Family of
graphs corresponding to a family of link diagrams $L$ will be
denoted by $G(L)$.

Tutte polynomials were known for the following special families of
graphs corresponding to the knots and links: circuit graphs $C_p$
which correspond to the link family $p$, graphs corresponding to the
link family $p\,2$, wheel graphs $Wh(n+1)$ corresponding to the
family of antiprismatic basic polyhedra $(2n)^*$ ($6^*$, $8^*$,
$10^*$, $\ldots $), so-called "hammock" graphs corresponding to the
pretzel links $2,2,\ldots ,2$ where tangle $2$ occurs $n$ times
($n\ge 3$), and graphs corresponding to the links of the form
$(2n)^*:2\,0:2\,0\ldots :2\,0$, where tangle $2$ occurs $n$ times
($6^*2\,0:2\,0:2\,0$, $8^*2\,0:2\,0:2\,0:2\,0$, $\ldots $)
[ChaShro]. Recursive formulas for the computation of Tutte
polynomials corresponding to the link families $3\,1\,3$,
$3\,1\,2\,1\,3$, $3\,1\,2\,1\,2\,1\,3$, $3\,1\,2\,1\,2\,1\,2\,1\,3$,
$\ldots $, $2\,1\,2$, $2\,1\,2\,1\,2$, $2\,1\,2\,1\,2\,1\,2$,
$2\,1\,2\,1\,2\,1\,2\,1\,2$, $\ldots $, and the family of polyhedral
links $.2:2$, $9^*2:2$, $\ldots $ are derived by F.~Emmert-Streib
[Emm]. Moreover, general formulas for the Jones polynomial are known
for some other particular classes of knots and links, e.g., torus
knots, or repeating chain knot with $3n$ crossings, proposed by
L.~Kauffman, given by Conway symbols of the form $3,2\,1\,1,2$;
$(3,2)\,1,2\,1\,1,2$; $((3,2)\,1,2)\,1,2\,1\,1,2$;
$(((3,2)\,1,2)\,1,2)\,1,2\,1\,1,2$;
$((((3,2)\,1,2)\,1,2)\,1,2)\,1,2\,1\,1,2$; $\ldots $ [WuWa]. Zeros
of the Jones polynomial are analyzed by X-S.~Lin [Li], S.~Chang and
R.~Shrock [ChaShro], F.~Wu and J.~Wang [WuWa], and X.~Jin and
F.~Zhang [JiZha]. Experimental observations of many authors who have
studied the distribution of roots of Jones polynomials for various
families of knots and links are explained by A.~Champanarekar and
I.~Kofman [ChaKo].

\section{Tutte polynomials of $KL$ families}

\subsection{Family \it{\textbf{p}}}

\begin{figure}[th]
\centerline{\psfig{file=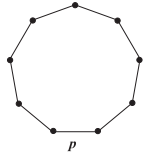,width=0.80in}} \vspace*{8pt}
\caption{Cycle graph $G(p)$.\label{f1.1}}
\end{figure}

The first family we consider is the family \emph{\textbf{p}} ($p\ge
1$), which consists of the knots and links $1_1$, $2_1^2$, $3_1$,
$4_1^2$, $5_1$, $\ldots$ [ChaShro]. Graphs corresponding to links of
this family are cycles of length $p$, which we can denote by $G(p)$.
By deleting one edge $G(p)$ gives the chain of edges of the length
$p-1$ with the Tutte polynomial $x^{p-1}$, and by contraction it
gives $G(p-1)$. Hence, $T(G(p))-T(G(p-1))=x^{p-1}$, and $T(G(1))=y$,
so the general formula for the Tutte polynomial of the graph $G(p)$
is

$$T(G(p))={{x^p-1} \over {x-1}} +y-1.$$

\subsection{Family \it{\textbf{p\,q}}}

\begin{figure}[th]
\centerline{\psfig{file=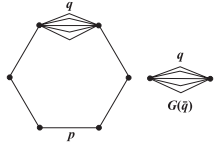,width=1.20in}} \vspace*{8pt}
\caption{Graph $G(p\,q)$.\label{f1.2}}
\end{figure}

The link family \emph{\textbf{p\,q}} gives the family of graphs,
illustrated in Fig. 2, satisfying the relations

$$T(G(p\,q))-T(G((p-1)\,q))=x^{p-1}T(G(\overline q)),$$ where
$G(\overline q)$ is the dual of the graph $G(q)$. Since the Tutte
polynomial of the graph $G(0\,q)$ is $T(G(0\,q))=y^q$, the general
formula for the Tutte polynomial of the graphs $G(p\,q)$ is

$$T(G(p\,q))={{(x^p-1)(y^q-1)}\over{(x-1)(y-1)}}+x^p+y^q-1.$$

\subsection{Family \it{\textbf{p\,1\,q}}}

\begin{figure}[th]
\centerline{\psfig{file=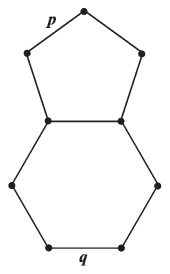,width=0.80in}} \vspace*{8pt}
\caption{Graph $G(p\,1\,q)$.\label{f1.3}}
\end{figure}

The family of graphs (Fig. 3) corresponds  to the link family
\emph{\textbf{p\,1\,q}}. Since

$$T(G(p\,1\,q))-T(G(p\,1\,(q-1)))=x^{q-1}G(p+1),$$

the general formula for the Tutte polynomial of the graphs
$G(p\,1\,q)$ is

$$T(G(p\,1\,q))={x(x^p-1)(x^q-1)\over (x-1)^2}+{{(x^p+x^q+xy-x-y-1)y}\over (x-1)}.$$

\subsection{Family \it{\textbf{p\,q\,r}}}

\begin{figure}[th]
\centerline{\psfig{file=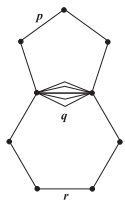,width=0.80in}} \vspace*{8pt}
\caption{Graph $G(p\,q\,r)$.\label{f1.4}}
\end{figure}

The family of graphs (Fig. 4) corresponds to the family of rational
links \emph{\textbf{p\,q\,r}} and their Tutte polynomials satisfy
the following relations

$$T(G(p\,q\,r))-T(G((p-1)\,q\,r))=x^{p-1}T(G(r\,q)).$$

\noindent The general formula for the Tutte polynomial of the graphs
$G(p\,q\,r)$ is

$$T(G(p\,q\,r))={(x+y)(x^p-1)(x^r-1)\over (x-1)^2}+{{y^q(x^{r+1}+x^p-x-1)}\over (x-1)}$$
$$+{(x^p-1)(x^r-1)(y^q-y^2)\over (x-1)^2(y-1)}-(x^r-y)y^q.$$

\subsection{Family \it{\textbf{p\,1\,1\,q}}}

\begin{figure}[th]
\centerline{\psfig{file=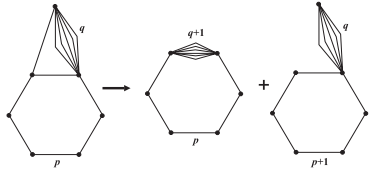,width=2.00in}} \vspace*{8pt}
\caption{Graph $G(p\,1\,1\,q)$.\label{f1.5}}
\end{figure}

The graph of the link family \emph{\textbf{p\,1\,1\,q}}, illustrated
in Fig. 5, resolves into the graph $G(p\,(q+1))$ and the block sum
of the graphs $G(p+1)$ and $G({\overline q})$. The general formula
for the Tutte polynomial of the graphs $G(p\,1\,1\,q)$ is

$$T(G(p\,1\,1\,q))=T(G(p\,(q+1))+({{x^{p+1}-1}\over {x-1}}+y-1)({{y^q-1}\over {y-1}}+x-1).$$

\subsection{Family \it{\textbf{p,q,r}}}

\begin{figure}[th]
\centerline{\psfig{file=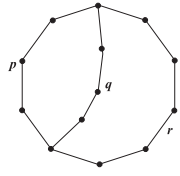,width=1.00in}} \vspace*{8pt}
\caption{Graph $G(p,q,r)$.\label{f1.6}}
\end{figure}

The graph family from Fig. 6 corresponds to the family of pretzel
links \emph{\textbf{p,q,r}}, whose Tutte polynomials satisfy the
relations

$$T(G(p,q,r))-T(G(p-1,q,r))=x^{p-1}G(q+r)$$

\noindent The general formula for the Tutte polynomial of the graphs
$G(p,q,r)$ is\footnote{In this paper we derived the general formula
for Tutte polynomial of pretzel links with five parameters, which
can be easily generalized to general formulas for pretzel links with
an arbitrary number of parameters, by using the relation
$T(G(p_1,\ldots ,p_n))={{x^{p_n}-1}\over {x-1}}T(G(p_1),\ldots
,G(p_{n-1}))+T(G(p_1))\ldots T(G(p_{n-1}))$ ($p\ge 3$).}

$$T(G(p,q,r))={{x^{p+q+r}+(x^{p+1}+x^{q+1}+x^{r+1})(y-1)-(x^p+x^q+x^r)y}\over (x-1)^2}+$$

$${(xy-x-y)(xy-x-y-1)\over (x-1)^2}.$$

\subsection{Antiprismatic basic polyhedra \it{\textbf{(2n)$^*$}}}

\begin{figure}[th]
\centerline{\psfig{file=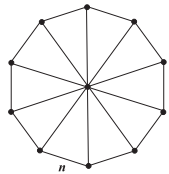,width=1.00in}} \vspace*{8pt}
\caption{Wheel graph $Wh(n+1)$.\label{f1.7}}
\end{figure}

Basic polyhedron $6^*$ is the first member of the class of
antiprismatic basic polyhedra \emph{\textbf{(2n)$^*$}} ($n\ge 3$):
$6^*$, $8^*$, $10^*$, $\ldots $. The corresponding graphs are wheel
graphs (Fig. 7) denoted by $Wh(n+1)$. Their Tutte polynomials are
given by the general formula [ChaShro]:

$$T(G((2n)^*)=T(Wh(n+1))=[{1\over 2}[(1+x+y)+[(1+x+y)^2-4xy]^{1/2}]]^n+$$

$$[{1\over 2}[(1+x+y)-[(1+x+y)^2-4xy]^{1/2}]]^n+{xy-x-y-1}.$$

\subsection{Family \it{\textbf{.p}}}

\begin{figure}[th]
\centerline{\psfig{file=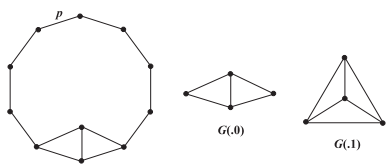,width=2.20in}} \vspace*{8pt}
\caption{Graphs $G(.p)$.\label{f1.8}}
\end{figure}

The graphs from Fig. 8. correspond to the link family
\emph{\textbf{.p}}. The Tutte polynomials of this graph family
satisfy the relations

$$T(G(.p))-T(G(.(p-1)))=x^{p-1}T(G(.0)),$$

\noindent $T(G(.0))=x+2x^2+x^3+y+2xy+y^2$  and
$T(G(.1))=2x+3x^2+x^3+2y+4xy+3y^2+y^3$. The general formula for the
Tutte polynomial of the graphs $G(.p)$ is

$$T(G(.p))={{x^p(2x^3+2x^2+y^2+2xy+y)-2x^2-2x-y^2-y-2xy}\over {x-1}}$$
$$-x^{p+1}(x+1)+y+2xy+2y^2+y^3.$$

\subsection{Family \it{\textbf{p\,q\,1\,r}}}

\begin{figure}[th]
\centerline{\psfig{file=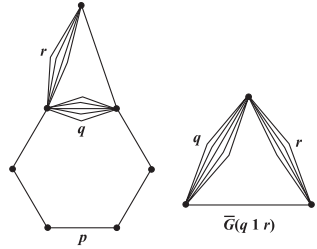,width=1.80in}} \vspace*{8pt}
\caption{Graphs $G(p\,q\,1\,r)$ and ${\overline
G}(q\,1\,r)$.\label{f1.9}}
\end{figure}

The graphs illustrated in Fig. 9  correspond to the link family
\emph{\textbf{p\,q\,1\,r}}. In order to obtain formula for the Tutte
polynomial we use the relations

$$T(G(p\,q\,1\,r))-T(G((p-1)\,q\,1\,r))=x^{p-1}T({\overline G}(q\,1\,r)),$$

\noindent where ${\overline G}(q\,1\,r)$ denotes the graph from Fig.
9. Since its Tutte polynomial is

$$T({\overline G}(q\,1\,r))=({{y^{r+1}-1}\over {y-1}}+x-1){{y^q-1}\over {y-1}}+x({{y^r-1}\over {y-1}}+x-1),$$

\noindent for the Tutte polynomial of the graphs $G(p\,q\,1\,r)$ we
obtain the general formula

$$T(G(p\,q\,1\,r))=({{y^{r+1}-1}\over {y-1}}+x-1){{y^q-1}\over
{y-1}}+x({{y^r-1}\over {y-1}}+x-1){{x^p-1}\over {x-1}}+
y^q({{y^{r+1}-1}\over {y-1}}+x-1).$$

\subsection{Family \it{\textbf{p\,1,q,r}}}

\begin{figure}[th]
\centerline{\psfig{file=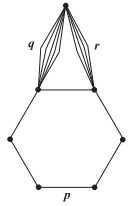,width=0.80in}} \vspace*{8pt}
\caption{Graph $G(p\,1,q,r)$.\label{f1.10}}
\end{figure}

The graphs illustrated in Fig. 10 correspond  to the link family
\emph{\textbf{p\,1,q,r}}. In order to obtain general formula for the
Tutte polynomial we use the relations

$$T(G(p\,1,q,r))-T(G((p-1)\,1,q,r))=x^{p-1}T({\overline G}(q\,1\,r)),$$

\noindent where ${\overline G}(q\,1\,r)$ denotes the graph from Fig.
9, which now occupies a different position with regards to the chain
of edges $p$. The general formula for the Tutte polynomial of the
graphs $G(p\,1,q,r)$ is

$$T(G(p\,1,q\,r))=({{y^{r+1}-1}\over {y-1}}+x-1){{y^q-1}\over
{y-1}}+x({{y^r-1}\over {y-1}}+x-1){{x^p-1}\over {x-1}}+
y({{y^{q+r}-1}\over {y-1}}+x-1).$$

\subsection{Family \it{\textbf{p,q,r+}}}

\begin{figure}[th]
\centerline{\psfig{file=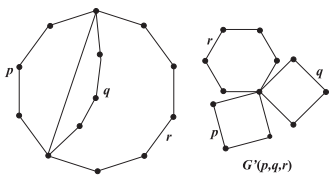,width=1.80in}} \vspace*{8pt}
\caption{Graphs $G(p,q,r+)$ and $G'(p,q,r)$.\label{f1.11}}
\end{figure}

The Tutte polynomial of the graphs corresponding to the link family
\emph{\textbf{p,q,r+}} (Fig. 11) we obtain from the Tutte polynomial
of the graphs $G(p,q,r)$, where the additional term is the Tutte
polynomial of the graphs $G'(p,q,r)$. The general formula for the
Tutte polynomial of the graphs $G(p,q,r+)$ is

$$T(G(p,q,r+))={T(G(p,q,r))+({{x^p-1}\over {x-1}}+y-1)({{x^q-1}\over {x-1}}+y-1)({{x^r-1}\over {x-1}}+y-1)}.$$

\subsection{Family \it{\textbf{p\,1\,1\,1\,q}}}

\begin{figure}[th]
\centerline{\psfig{file=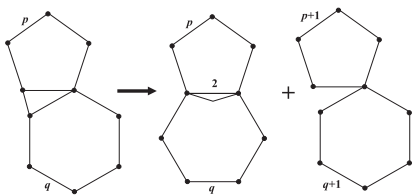,width=2.40in}} \vspace*{8pt}
\caption{Resolving the graph $G(p\,1\,1\,1\,q)$ into graphs
$G(p\,2\,q)$ and $G'(p+1,q+1)$.\label{f1.12}}
\end{figure}

By resolving the graph $G(p\,1\,1\,1\,q)$ into the graphs
$G(p\,2\,q)$ and $G'(p+1,q+1)$ (Fig. 12) we obtain the general
formula for the Tutte polynomial of the graphs $G(p\,1\,1\,1\,q)$

$$T(G(p\,1\,1\,1\,q))={T(G(p\,2\,q))+({{x^{p+1}-1}\over {x-1}}+y-1)({{x^{q+1}-1}\over {x-1}}+y-1)}.$$

\subsection{Family \it{\textbf{.p\,1}}}

\begin{figure}[th]
\centerline{\psfig{file=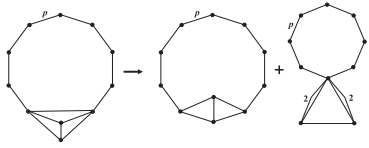,width=2.00in}} \vspace*{8pt}
\caption{Resolving the graph $G(.p\,1)$ into the graph $G(.p)$ and
the block sum of graphs ${\overline G}(2\,1\,2)$ and
$G(p)$.\label{f1.13}}
\end{figure}

The graph of the link family \emph{\textbf{.p\,1}}  (Fig. 13)
resolves into the graph $G(.p)$ and the block sum of graphs
${\overline G}(2\,1\,2)$ and $G(p)$. The general formula for the
Tutte polynomial of the graphs $G(.p\,1)$ is

$$T(G(.p\,1))=T(G(.p))+(x+x^2+y+2xy+2y^2+y^3)({{{x^p-1}\over
{x-1}}+y-1}).$$

\subsection{Family \it{\textbf{.p:q}}}

\begin{figure}[th]
\centerline{\psfig{file=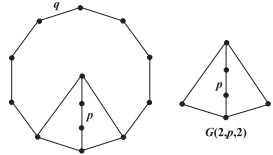,width=1.50in}} \vspace*{8pt}
\caption{Graphs $G(.p:q)$ and $G(.p:0)=G(2,p,2)$.\label{f1.14}}
\end{figure}

The graph family illustrated in Fig. 14   corresponds  to the link
family \emph{\textbf{.p:q}}. In order to obtain general formula for
the Tutte polynomial we can use the relations

$$T(G(.p:q))-T(G(.p:(q-1)))=x^{q-1}T(G(2,p,2)).$$

\noindent where $G(.p:0)$ is the graph $G(2,p,2)$ with the Tutte
polynomial

$$T(G(.p:0))=T(G(2,p,2))={{x^p(x^3+x^2+x+y)+(2x+y+1)(xy-x-y)}\over {x-1}}.$$

\noindent The general formula for the Tutte polynomial of the graph
$G(.p:q)$ is

$$T(G(.p:q))={{{x^q-1}\over {(x-1)^2}}(x^p(x^3+x^2+x+y)+(2x+y+1)(xy-x-y))}+$$

$${{x^p-1}\over {x-1}}(x+y)^2+x+y+y^2+y^3.$$

\subsection{Family \it{\textbf{.p.q}}}

\begin{figure}[th]
\centerline{\psfig{file=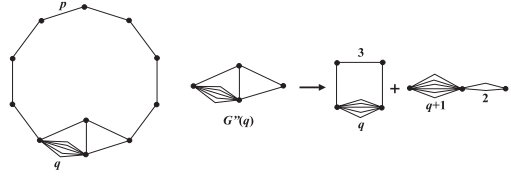,width=2.80in}} \vspace*{8pt}
\caption{Graph $G(.p.q)$ and resolving of the graph
$G''(q)$.\label{f1.15}}
\end{figure}

The graph family illustrated in Fig. 15  corresponds  to the link
family \emph{\textbf{.p.q}}. In order to obtain formula for the
Tutte polynomial we can use the relations

$$T(G(.p.q))-T(G(.(p-1).q))=x^{p-1}T(G''(q)).$$

\noindent where $G''(q)$ is the graph which resolves into the graph
$G(3\,q)$ and the block sum of graphs $G({\overline {q+1}})$ and
$G({\overline 2})$. The general formula for the Tutte polynomial of
the graph $G''(q)$ is

$$T(G''(q))={{(1+2x+x^2+y)(xy-x-y)+y^q(x+x^2+y+xy+y^2)}\over
{y-1}},$$

\noindent and the Tutte polynomial of the graph $G(.p.q)$ is

$$T(G(.p.q))={{{x^p-1}\over {x-1}}T(G''(q))}+2y{{y^{q+2}-1}\over {y-1}}-y-y^{q+2}+$$

$$xy{{y^q-1}\over {y-1}}+xy+x+x^2.$$

\subsection{Family \it{\textbf{.p:q\,0}}}

\begin{figure}[th]
\centerline{\psfig{file=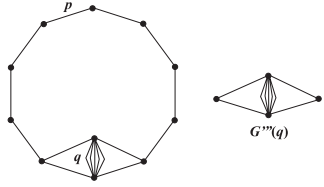,width=1.80in}} \vspace*{8pt}
\caption{Graphs $G(.p:q\,0)$ and $G'''(q)$.\label{f1.16}}
\end{figure}

The graph family $G(.p:q\,0)$ illustrated in Fig. 16   corresponds
to the link family \emph{\textbf{.p:q\,0}}. In order to obtain
general formula for the Tutte polynomial we can use the relations

$$T(G(.p.q))-T(G(.(p-1).q))=x^{p-1}T(G'''(q)).$$

\noindent The Tutte polynomial of the graph $G'''(q)$ is

$$T(G'''(q))={{{y^q-1}\over {y-1}}(x+y)^2+x+x^2+x^3+y},$$

\noindent and the general formula for the Tutte polynomial of the
graph $G(.p:q\,0)$ is

$$T(G(.p:q\,0))={{{x^p-1}\over {x-1}}T(G'''(q))}+2y{{y^{q+2}-1}\over {y-1}}-2y-y^{q+2}-xy^q+$$

$$xy{{y^q-1}\over {y-1}}+2xy+x+x^2+{{y^{q+1}-1}\over {y-1}}-1.$$

\subsection{Family \it{\textbf{.p.q\,0}}}

\begin{figure}[th]
\centerline{\psfig{file=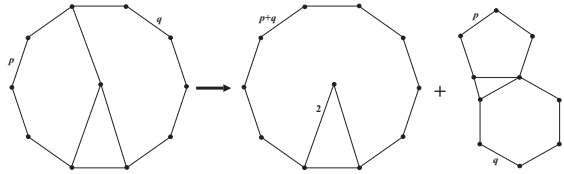,width=3.00in}} \vspace*{8pt}
\caption{Resolving the graph $G(.p.q\,0)$ into the graphs
$G((p+q)\,1\,2)$ and $G(p\,1\,1\,1\,q)$.\label{f1.17}}
\end{figure}

The graph family $G(.p.q\,0)$ illustrated in Fig. 16  corresponds to
the link family \emph{\textbf{.p.q\,0}}. The graph $G(.p.q\,0)$
resolves into the graphs $G((p+q)\,1\,2)$ and $G(p\,1\,1\,1\,q)$.
The general formula for the Tutte polynomial of the graphs
$G(.p.q\,0)$ is

$$T(G(.p.q\,0))=T(G((p+q)\,1\,2))+T(G(p\,1\,1\,1\,q)).$$

\subsection{Family \it{\textbf{p\,1,q\,1,r}}}

\begin{figure}[th]
\centerline{\psfig{file=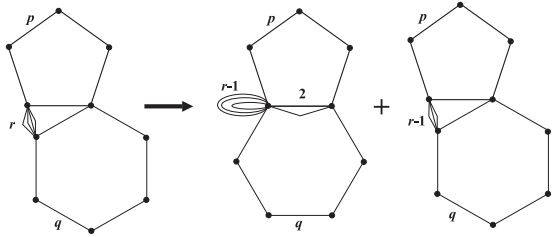,width=3.00in}} \vspace*{8pt}
\caption{Resolving the graph $G(p\,1,q\,1,r)$.\label{f1.18}}
\end{figure}

The graph family $G(p\,1,q\,1,r)$ illustrated in Fig. 18 corresponds
to the link family \emph{\textbf{p\,1,q\,1,r}}. In order to obtain
general formula for the Tutte polynomial we can use the relations

$$T(G(p\,1,q\,1,r))-T(G(p\,1,q\,1,(r-1)))=y^{r-1}T(G(p\,2\,q)).$$

\noindent The general formula for the Tutte polynomial of the graphs
$G(p\,1,q\,1,r)$ is

$$T(G(p\,1,q\,1,r))=T(G(p\,2\,q)){{y^r-1}\over {y-1}}+$$

$$({{x^{p+1}-1}\over {x-1}}+y-1)({{x^{q+1}-1}\over {x-1}}+y-1).$$

\subsection{Family \it{\textbf{p\,1,q,r+}}}

\begin{figure}[th]
\centerline{\psfig{file=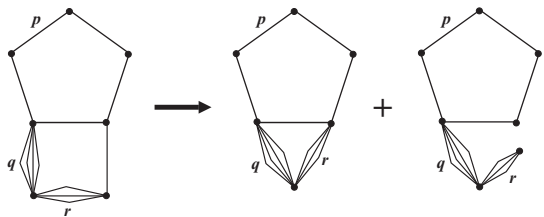,width=3.00in}} \vspace*{8pt}
\caption{Resolving the graph $G(p\,1,q,r+)$.\label{f1.19}}
\end{figure}

The graph family $G(p\,1,q,r+)$ illustrated in Fig. 19 corresponds
to the link family \emph{\textbf{p\,1,q,r+}}. The graph
$G(p\,1,q,r+)$ resolves into the graph $G(p\,1,q,r)$ and the block
sum of the graphs $G(p+1)$, $G(\overline q)$, and $G(\overline r)$.
The general formula for the Tutte polynomial of the graphs
$G(p\,1,q,r+)$ is

$$T(G(p\,1,q,r+))=T(G(p\,1,q,r))+({{x^{p+1}-1}\over {x-1}}+y-1)({{y^q-1}\over {y-1}}+x-1)({{y^r-1}\over {y-1}}+x-1).$$

\subsection{Family \it{\textbf{p\,1\,1,q,r}}}

\begin{figure}[th]
\centerline{\psfig{file=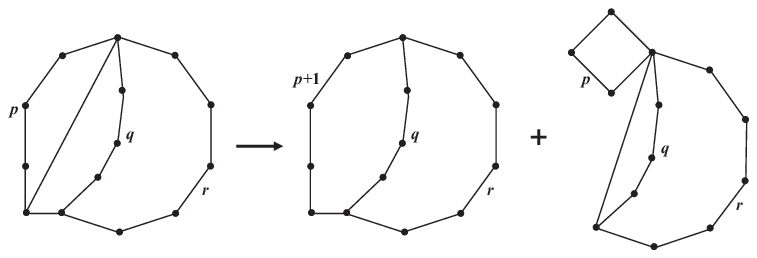,width=3.00in}} \vspace*{8pt}
\caption{Resolving the graph $G(p\,1\,1,q,r)$.\label{f1.20}}
\end{figure}

The graph family $G(p\,1\,1,q,r)$ illustrated in Fig. 20 corresponds
to the link family \emph{\textbf{p\,1\,1,q,r}}. The graph
$G(p\,1\,1,q,r)$ resolves into the graphs $G((p+1),q,r)$ and the
block sum of the graphs $G(1,q,r)$ and $G(p)$. The general formula
for the Tutte polynomial of the graphs $G(p\,1\,1,q,r)$ is

$$T(G(p\,1\,1,q,r))=T(G((p+1),q,r))+T(G(1,q,r))T(G(p)).$$

\subsection{Family \it{\textbf{p\,1\,q\,1\,r}}}

\begin{figure}[th]
\centerline{\psfig{file=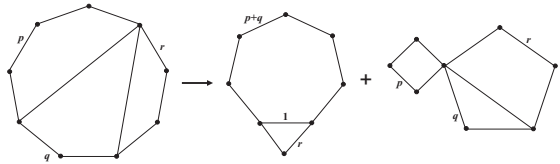,width=3.00in}} \vspace*{8pt}
\caption{Resolving the graph $G(p\,1\,q\,1\,r)$.\label{f1.21}}
\end{figure}

The graph family $G(p\,1\,q\,1\,r)$ illustrated in Fig. 21
corresponds  to the link family \emph{\textbf{p\,1\,q\,1\,r}}. The
graph $G(p\,1\,q\,1\,r)$  resolves into the graph $G((p+q)\,1\,r)$
and the block sum of the graphs $G(q\,1\,r)$ and $G(p)$. The general
formula for the Tutte polynomial of the graph $G(p\,1\,q\,1\,r)$ is

$$T(G(p\,1\,q\,1\,r))=T(G((p+q)\,1\,r))+T(G(q\,1\,r))T(G(p)).$$

\subsection{Family \it{\textbf{p,q,r,s}}}

\begin{figure}[th]
\centerline{\psfig{file=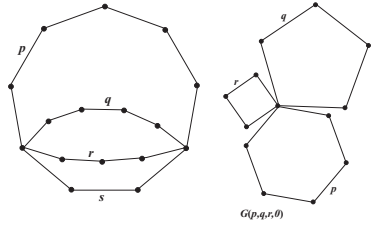,width=2.00in}} \vspace*{8pt}
\caption{Graphs $G(p,q,r,s)$ and $G(p,q,r,0)$.\label{f1.22}}
\end{figure}

The graph family $G(p,q,r,s)$ illustrated in Fig. 22  corresponds to
the link family \emph{\textbf{p,q,r,s}}. In order to obtain general
formula for the Tutte polynomial we can use the relations

$$T(G(p,q,r,s))-T(G(p,q,r,(s-1)))=x^{s-1}T(G(p,q,r)).$$

\noindent Since the Tutte polynomial of the graph $G(p,q,r,0)$ is
$$T(G(p,q,r,0))=({{x^p-1}\over {x-1}}+y-1)({{x^q-1}\over
{x-1}}+y-1)({{x^r-1}\over {x-1}}+y-1),$$

\noindent the general formula for the Tutte polynomial of the graphs
$G(p,q,r,s)$ is

$$T(G(p,q,r,s))={{x^s-1}\over {x-1}}T(G(p,q,r))+$$

$$({{x^p-1}\over {x-1}}+y-1)({{x^q-1}\over {x-1}}+y-1)({{x^r-1}\over {x-1}}+y-1).$$

\subsection{Family \it{\textbf{p\,1\,1\,1\,1\,q}}}

\begin{figure}[th]
\centerline{\psfig{file=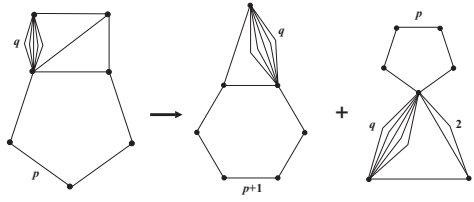,width=2.50in}} \vspace*{8pt}
\caption{Resolving the graph $G(p\,1\,1\,1\,1\,q)$.\label{f1.23}}
\end{figure}

The graph family $G(p\,1\,1\,1\,1\,q)$ illustrated in Fig. 23
corresponds  to the link family \emph{\textbf{p\,1\,1\,1\,1\,q}}.
The graph  $G(p\,1\,1\,1\,1\,q)$ resolves into the graph
$G((p+1)\,1\,1\,q)$ and the block sum of the graphs $\overline
G(q\,1\,2)$ and $G(p)$. The general formula for the Tutte polynomial
of the graphs $G(p\,1\,1\,1\,1\,q)$ is

$$T(G(p\,1\,1\,1\,1\,q))=T(G((p+1)\,1\,1\,q))+T(\overline G(q\,1\,2))T(G(p)).$$

\subsection{Family \it{\textbf{p\,q\,1\,1\,r}}}

\begin{figure}[th]
\centerline{\psfig{file=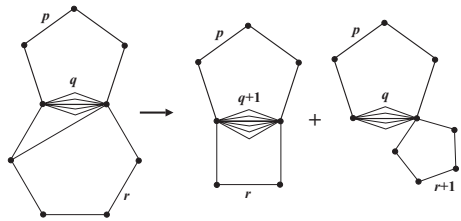,width=2.50in}} \vspace*{8pt}
\caption{Resolving the graph $G(p\,q\,1\,1\,r)$.\label{f1.24}}
\end{figure}

The graph family $G(p\,q\,1\,1\,r)$ illustrated in Fig. 24
corresponds  to the link family \emph{\textbf{p\,q\,1\,1\,r}}. The
graph $G(p\,q\,1\,1\,r)$ resolves into the graph $G(p\,(q+1)\,r)$
and the block sum of the graphs $G(p\,q)$ and $G(r+1)$. The general
formula for the Tutte polynomial of the graphs $G(p\,q\,1\,1\,r)$ is

$$T(G(p\,q\,1\,1\,r))=T(G(p\,(q+1)\,r))+T(G(p\,q))T(G(r+1)).$$

\subsection{Family \it{\textbf{p\,q\,r\,s}}}

\begin{figure}[th]
\centerline{\psfig{file=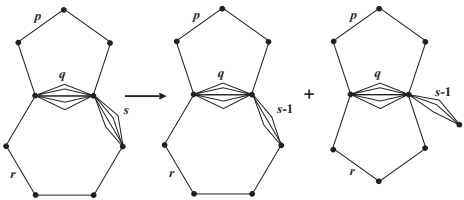,width=2.50in}} \vspace*{8pt}
\caption{The relations for the graph $G(p\,q\,r\,s)$.\label{f1.25}}
\end{figure}

The graph family $G(p\,q\,r\,s)$ illustrated in Fig. 25  corresponds
to the link family \emph{\textbf{p\,q\,r\,s}}. In order to obtain
general formula for the Tutte polynomial we can use the relations

$$T(G(p\,q\,r\,s))-T(G(p\,q\,r\,(s-1)))=y^{s-1}T(G(p\,q\,r)).$$

\noindent The general formula for the Tutte polynomial of the graphs
$G(p\,q\,r\,s)$ is

$$T(G(p\,q\,r\,s))={{y^s-1}\over {y-1}}T(G(p\,q\,r))+x^r({{x^{p+1}-1}\over {x-1}}+y-1)+$$

$$x^ry{{y^{q-1}-1}\over {y-1}}({{x^p-1}\over {x-1}}+y-1).$$

\subsection{Family \it{\textbf{p\,q,r,s}}}

\begin{figure}[th]
\centerline{\psfig{file=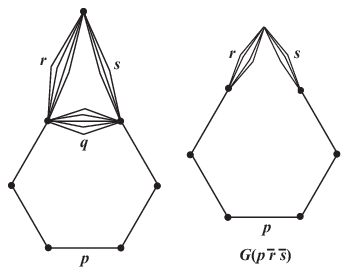,width=1.80in}} \vspace*{8pt}
\caption{Graphs $G(p\,q,r,s)$ and $G(p\,\overline r\,\overline
s)$.\label{f1.26}}
\end{figure}

The graph family $G(p\,q,r,s)$ illustrated in Fig. 26  corresponds
to the link family \emph{\textbf{p\,q,r,s}}. In order to obtain
general formula for the Tutte polynomial we can use the relations

$$T(G(p\,q,r,s))-T(G(p\,(q-1),r,s))=y^{q-1}T(G(p))T(\overline
G(r+s)).$$

\noindent After computing the Tutte polynomial of the graphs
$G(p\,\overline r\,\overline s)$

$$T(G(p\,\overline r\,\overline s))={{x^p-1}\over {x-1}}({{y^r-1}\over
{y-1}}+x-1)({{y^s-1}\over {y-1}}+x-1)+({{y^{r+s}-1}\over
{y-1}}+x-1),$$

\noindent we conclude that the general formula for the Tutte
polynomial of the graphs $G(p\,q,r,s)$ is

$$T(G(p\,q,r,s))=T(G(p\,\overline r\,\overline s))+{{y^q-1}\over {y-1}}({{x^p-1}\over {x-1}}+y-1)({{y^{r+s}-1}\over
{y-1}}+x-1).$$

\subsection{Family \it{\textbf{p,q,r+s}}}

\begin{figure}[th]
\centerline{\psfig{file=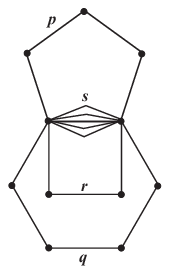,width=1.00in}} \vspace*{8pt}
\caption{Graph $G(p,q,r+s)$.\label{f1.27}}
\end{figure}

The graph family $G(p,q,r+s)$ illustrated in Fig. 27  corresponds to
the link family \emph{\textbf{p,q,r+s}}. In order to obtain general
formula for the Tutte polynomial we can use the relations

$$T(G(p,q,r+s))-T(G(p,q,r+(s-1)))=y^{s-1}T(G(p))T(G(q))T(G(r))$$

\noindent and $T(G(p,q,r+0))=T(G(p,q,r)).$ The general formula for
the Tutte polynomial of the graphs $G(p,q,r+s)$ is

$$T(p,q,r+s)={{y^s-1}\over {y-1}}({{x^p-1}\over {x-1}}+y-1)({{x^q-1}\over {x-1}}+y-1)({{x^r-1}\over
{x-1}}+y-1)+T(G(p,q,r)).$$

\subsection{Family \it{\textbf{(p,q)\,(r,s)}}}

\begin{figure}[th]
\centerline{\psfig{file=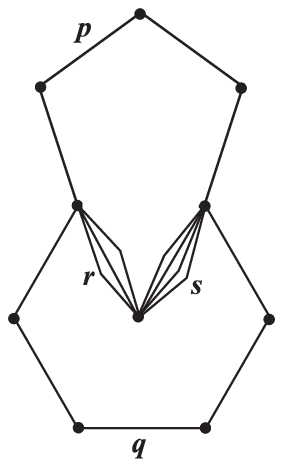,width=1.00in}} \vspace*{8pt}
\caption{Graph $G(p,q)\,(r,s)$.\label{f1.28}}
\end{figure}

The graph family $G(p,q)\,(r,s)$ illustrated in Fig. 28  corresponds
to the link family \emph{\textbf{(p,q)\,(r,s)}}. In order to obtain
general formula for the Tutte polynomial we can use the relations

$$T(G((p,q)\,(r,s)))-T(G((p,(q-1))\,(r,s)))=x^{q-1}T(G(p\,\overline r\,\overline
s)).$$

\noindent The general formula for the Tutte polynomial of the graphs
$G(p,q)\,(r,s)$ is

$$T(G((p,q)\,(r,s)))={{x^q-1}\over {x-1}}T(G(p\,\overline r\,\overline
s))+({{x^p-1}\over {x-1}}+y-1)({{y^{r+s}-1}\over {y-1}}+x-1).$$

\subsection{Family \it{\textbf{p\,q\,1\,r\,s}}}

\begin{figure}[th]
\centerline{\psfig{file=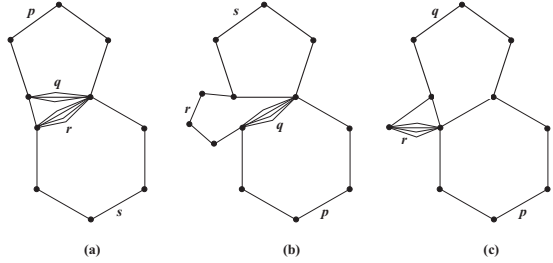,width=3.00in}} \vspace*{8pt}
\caption{Graphs (a) $G(p\,q\,1\,r\,s)$; (b) $G(p\,q\,r\,1\,s)$; (c)
$G(p\,1\,q\,1\,1\,r)$.\label{f1.29}}
\end{figure}

Graph $G(p\,q\,1\,r\,s)$ (Fig. 29a) resolves into the graph
$G(p\,(q+r)\,s)$ and the block sum of graphs $G(p\,q)$ and
$G(r\,s)$. The general formula for the Tutte polynomial of the
graphs $G(p\,q\,1\,r\,s)$ is

$$T(G(p\,q\,1\,r\,s))=T(G(p\,(q+r)\,s))+T(G(p\,q))T(G(s\,r)).$$

\subsection{Family \it{\textbf{p\,q\,r\,1\,s}}}

Graph $G(p\,q\,r\,1\,s)$ (Fig. 29b) resolves into the graph
$G(p\,q\,(r+s))$ and the block sum of graphs $G(p\,q\,r)$ and
$G(s)$. The general formula for the Tutte polynomial of the graphs
$G(p\,q\,r\,1\,s)$ is

$$T(G(p\,q\,r\,1\,s))=T(G(p\,q\,(r+s)))+T(G(p\,q\,r))T(G(s)).$$

\subsection{Family \it{\textbf{p\,1\,q\,1\,1\,r}}}

Graph $G(p\,1\,q\,1\,1\,r)$ (Fig. 29c) resolves into the graph
$G((p+q)\,1\,1\,r)$ and the block sum of graphs $G(q\,1\,1\,r)$ and
$G(p)$. The general formula for the Tutte polynomial of the graphs
$G(p\,1\,q\,1\,1\,r)$ is

$$T(G(p\,1\,q\,1\,1\,r))=T(G((p+q)\,1\,1\,r))+T(G(q\,1\,1\,r))T(G(p)).$$

\subsection{Family \it{\textbf{p\,q\,1\,1\,1\,r}}}

\begin{figure}[th]
\centerline{\psfig{file=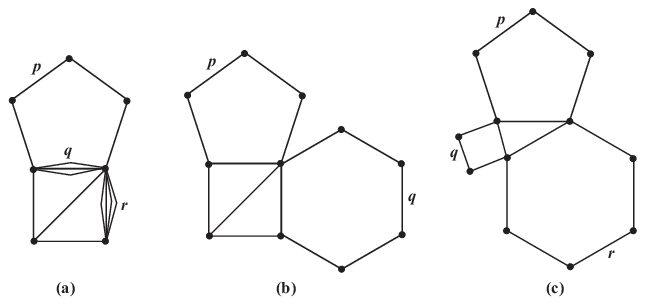,width=3.00in}} \vspace*{8pt}
\caption{Graphs (a) $G(p\,q\,1\,1\,1\,r)$; (b)
$G(p\,1\,1\,1\,1\,1\,q)$; (c) $G(p\,1,q\,1,r\,1)$.\label{f1.30}}
\end{figure}

Graph $G(p\,q\,1\,1\,1\,r)$ (Fig. 30a) resolves into the graph
$G(p\,(q+1)\,1\,r)$ and the block sum of graphs $G(p\,q)$ and
$G(2\,r)$. The general formula for the Tutte polynomial of the
graphs $G(p\,q\,1\,1\,1\,r)$ is

$$T(G(p\,q\,1\,1\,1\,r))=T(G(p\,(q+1)\,1\,r))+T(G(p\,q))T(G(2\,r)).$$

\subsection{Family \it{\textbf{p\,1\,1\,1\,1\,1\,q}}}

Graph $G(p\,1\,1\,1\,1\,1\,q)$ (Fig. 30b) resolves into the graph
$G(p\,1\,1\,1\,(q+1))$ and the block sum of graphs $G(p\,1\,1\,2)$
and $G(q)$. The general formula for the Tutte polynomial of the
graphs $G(p\,1\,1\,1\,1\,1\,q)$ is

$$T(G(p\,1\,1\,1\,1\,1\,q))=T(G(p\,1\,1\,1\,(q+1)))+T(G(p\,1\,1\,2))T(G(q)).$$

\subsection{Family \it{\textbf{p\,1,q\,1,r\,1}}}

Graph $G(p\,1,q\,1,r\,1)$ (Fig. 30c) resolves into the graph
$G(p\,1\,q\,1\,r)$ and the block sum of graphs $G(p\,2\,r)$ and
$G(q)$. The general formula for the Tutte polynomial of the graphs
$G(p\,1,q\,1,r\,1)$ is

$$T(G(p\,1,q\,1,r\,1))=T(G(p\,1\,q\,1\,r))+T(G(p\,2\,r))T(G(q)).$$

\subsection{Family \it{\textbf{p\,1\,q,r,s}}}

\begin{figure}[th]
\centerline{\psfig{file=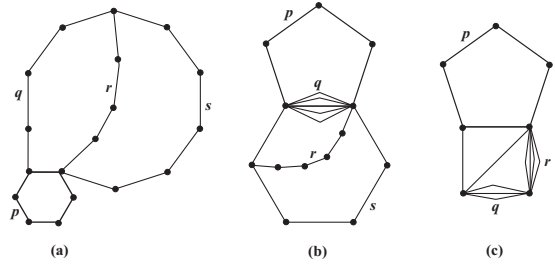,width=3.00in}} \vspace*{8pt}
\caption{Graphs (a) $G(p\,1\,q,r,s)$; (b) $G(p\,q\,1,r,s)$; (c)
$G(p\,1\,1\,1,q,r)$.\label{f1.31}}
\end{figure}

Graph $G(p\,1\,q,r,s)$ (Fig. 31a) resolves into the graph
$G((p+q),r,s)$ and the block sum of graphs $G(q,r,s)$ and $G(p)$.
The general formula for the Tutte polynomial of the graphs
$G(p\,1\,q,r,s)$ is

$$T(G(p\,1\,q,r,s))=T(G((p+q),r,s))+T(G(q,r,s))T(G(p)).$$

\subsection{Family \it{\textbf{p\,q\,1,r,s}}}

Graph $G(p\,q\,1,r,s)$ (Fig. 31b) resolves into the graph
$G(p,r,s+q)$ and the block sum of graphs $G(p\,q)$ and $G(r+s)$. The
general formula for the Tutte polynomial of the graphs
$G(p\,q\,1,r,s)$ is

$$T(G(p\,q\,1,r,s))=T(G(p,r,s+q))+T(G(p\,q))T(G(r+s)).$$

\subsection{Family \it{\textbf{p\,1\,1\,1,q,r}}}

Graph $G(p\,1\,1\,1,q,r)$ (Fig. 31c) resolves into the graph
$G(p\,2,q,r)$ and the block sum of graphs $\overline G(q\,1\,r)$ and
$G(p+1)$. The general formula for the Tutte polynomial of the graphs
$G(p\,1\,1\,1,q,r)$ is

$$T(G(p\,1\,1\,1,q,r))=T(G(p\,2,q,r))+T(\overline G(q\,1\,r))T(G(p+1)).$$

\subsection{Family \it{\textbf{p\,1,q,r,s}}}

\begin{figure}[th]
\centerline{\psfig{file=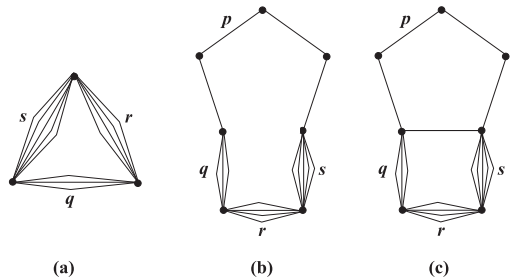,width=2.80in}} \vspace*{8pt}
\caption{Graphs (a) $\overline G(q\,r\,s)$; (b) $G(p\,\overline
q\,\overline r\,\overline s)$; (c) $G(p\,1,q,r,s)$.\label{f1.32}}
\end{figure}

The Tutte polynomial of the graph $\overline G(q\,r\,s)$ (Fig. 32a)
is

$$T(\overline G(q\,r\,s))=\overline G(r)\overline G(s)+(y{{y^{q-1}-1}\over {y-1}}+1)\overline G(r+s),$$

\noindent and the Tutte polynomial of the graph $G(p\,\overline
q\,\overline r\,\overline s)$ (Fig. 32b) is

$$T(G(p\,\overline
q\,\overline r\,\overline s))={{x^p-1}\over {x-1}}({{y^q-1}\over
{y-1}}+x-1)({{y^r-1}\over {y-1}}+x-1)({{y^s-1}\over
{y-1}}+x-1)+T(\overline G(q\,r\,s)).$$

\noindent Since graph $G(p\,1,q,r,s)$ (Fig. 32c) resolves into the
graph $G(p\,\overline q\,\overline r\,\overline s)$ and the block
sum of graphs $\overline G(q\,r\,s)$ and $G(p)$, the general formula
for the Tutte polynomial of the graphs $G(p\,1,q,r,s)$ is

$$T(G(p\,1,q,r,s))=T(G(p\,\overline
q\,\overline r\,\overline s))+T(\overline G(q\,r\,s))T(G(p)).$$

\subsection{Family \it{\textbf{p\,q,r\,1,s}}}

\begin{figure}[th]
\centerline{\psfig{file=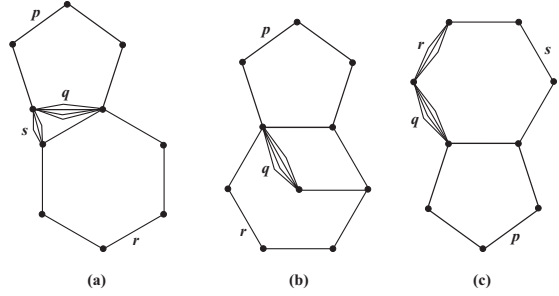,width=3.00in}} \vspace*{8pt}
\caption{Graphs (a) $G(p\,q,r\,1,s)$; (b) $G(p\,1\,1,q\,1,r)$; (c)
$G(p\,1,q,r+s)$.\label{f1.33}}
\end{figure}

Graph $G(p\,q,r\,1,s)$ (Fig. 33a) resolves into the graph
$G(p\,q\,r\,s)$ and the block sum of graphs $G(p\,(q+s))$ and
$G(r)$. The general formula for the Tutte polynomial of the graphs
$G(p\,q,r\,1,s)$ is

$$T(G(p\,q,r\,1,s))=T(G(p\,q\,r\,s))+T(G(p\,(q+s)))T(G(r)).$$

\subsection{Family \it{\textbf{p\,1\,1,q\,1,r}}}

Graph $G(p\,1\,1,q\,1,r)$ (Fig. 33b) resolves into the graph
$G(r\,q\,1\,1\,p)$ and the block sum of graphs $G(p\,1\,(r+1))$ and
$G(\overline q)$. The general formula for the Tutte polynomial of
the graphs $G(p\,1\,1,q\,1,r)$ is

$$T(G(p\,1\,1,q\,1,r))=T(G(r\,q\,1\,1\,p))+T(G(p\,1\,(r+1)))T(G(\overline q)).$$

\subsection{Family \it{\textbf{p\,1,q,r+s}}}

Graph $G(p\,1,q,r+s)$ (Fig. 33c) resolves into the graph
$G((p+s)\,\overline q\,\overline r)$ and the block sum of graphs
$G(s\,\overline q\,\overline r)$ and $G(p)$. The general formula for
the Tutte polynomial of the graphs $G(p\,1,q,r+s)$ is

$$T(G(p\,1,q,r+s))=T(G((p+s)\,\overline
q\,\overline r))+T(G(s\,\overline q\,\overline r))T(G(p)).$$

\subsection{Family \it{\textbf{p,q,r,s+}}}

\begin{figure}[th]
\centerline{\psfig{file=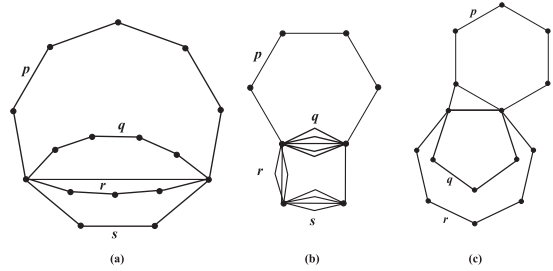,width=3.00in}} \vspace*{8pt}
\caption{Graphs (a) $G(p,q,r,s+)$; (b) $G(p\,q,r,s+)$; (c)
$G(p\,1\,1,q,r+)$.\label{f1.34}}
\end{figure}

Graph $G(p,q,r,s+)$ (Fig. 34a) resolves into the graph $G(p,q,r,s)$
and the block sum of graphs $G(p)$, $G(q)$, $G(r)$ and $G(s)$. The
general formula for the Tutte polynomial of the graphs $G(p,q,r,s+)$
is

$$T(G(p,q,r,s+))=T(G(p,q,r,s))+T(G(p))T(G(q))T(G(r))T(G(s)).$$

\subsection{Family \it{\textbf{p\,q,r,s+}}}

Graph $G(p\,q,r,s+)$ (Fig. 34b) resolves into the graph
$G(p\,q,r,s)$ and the block sum of graphs $G(p\,q)$, $G(\overline
r)$, and $G(\overline s)$. The general formula for the Tutte
polynomial of the graphs $G(p\,q,r,s+)$ is

$$T(G(p\,q,r,s+))=T(G(p\,q,r,s))+T(G(p\,q))T(G(\overline r))T(G(\overline s)).$$

\subsection{Family \it{\textbf{p\,1\,1,q,r+}}}

Graph $G(p\,1\,1,q,r+)$ (Fig. 34c) resolves into the graph
$G(p,q,r++)=G(p,q,r+2)$ and the block sum of graphs $G(q\,1\,r)$ and
$G((p+1))$. The general formula for the Tutte polynomial of the
graphs $G(p\,1\,1,q,r+)$ is

$$T(G(p\,1\,1,q,r+))=T(G(p,q,r+2))+T(G(q\,1\,r))T(G((p+1))).$$

\subsection{Family \it{\textbf{(p,q+)\,(r,s)}}}

\begin{figure}[th]
\centerline{\psfig{file=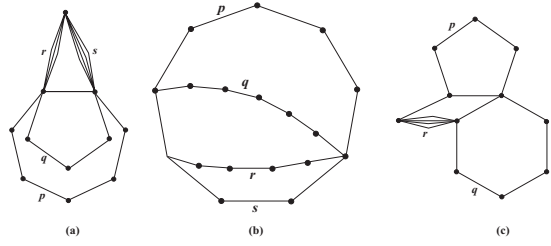,width=3.00in}} \vspace*{8pt}
\caption{Graphs (a) $G((p,q+)\,(r,s))$; (b) $G((p,q)\,1\,(r,s))$;
(c) $G(p\,1,q\,1,r+)$.\label{f1.35}}
\end{figure}

Graph $G((p,q+)\,(r,s))$ (Fig. 35a) resolves into the graph
$G((p,q)\,(r,s))$ and the block sum of graphs $G(p)$, $G(q)$ and
$G((\overline {r+s}))$. The general formula for the Tutte polynomial
of the graphs $G((p,q+)\,(r,s))$ is

$$T(G((p,q+)\,(r,s)))=T(G((p,q)\,(r,s)))+T(G(p))T(G(q))T(G((\overline {r+s}))).$$

\subsection{Family \it{\textbf{(p,q)\,1\,(r,s)}}}

Graph $G((p,q)\,1\,(r,s))$ (Fig. 35b) resolves into the graph
$G(p,q,r,s)$ and the block sum of graphs $G((p+q))$ and $G((r+s))$.
The general formula for the Tutte polynomial of the graphs
$G((p,q)\,1\,(r,s))$ is

$$T(G((p,q)\,1\,(r,s)))=T(G(p,q,r,s))+T(G((p+q)))T(G((r+s))).$$

\subsection{Family \it{\textbf{p\,1,q\,1,r+}}}

Graph $G(p\,1,q\,1,r+)$ (Fig. 35c) resolves into the graph
$G(p\,1,q\,1,r)$ and the block sum of graphs $G((p+1))$, $G((q+1))$,
and $G(\overline r)$. The general formula for the Tutte polynomial
of the graphs $G(p\,1,q\,1,r+)$ is

$$T(G(p\,1,q\,1,r+))=T(G(p\,1,q\,1,r))+T(G((p+1)))T(G((q+1))T(G(\overline r)).$$

\subsection{Family \it{\textbf{(p\,1,q)\,(r,s)}}}

\begin{figure}[th]
\centerline{\psfig{file=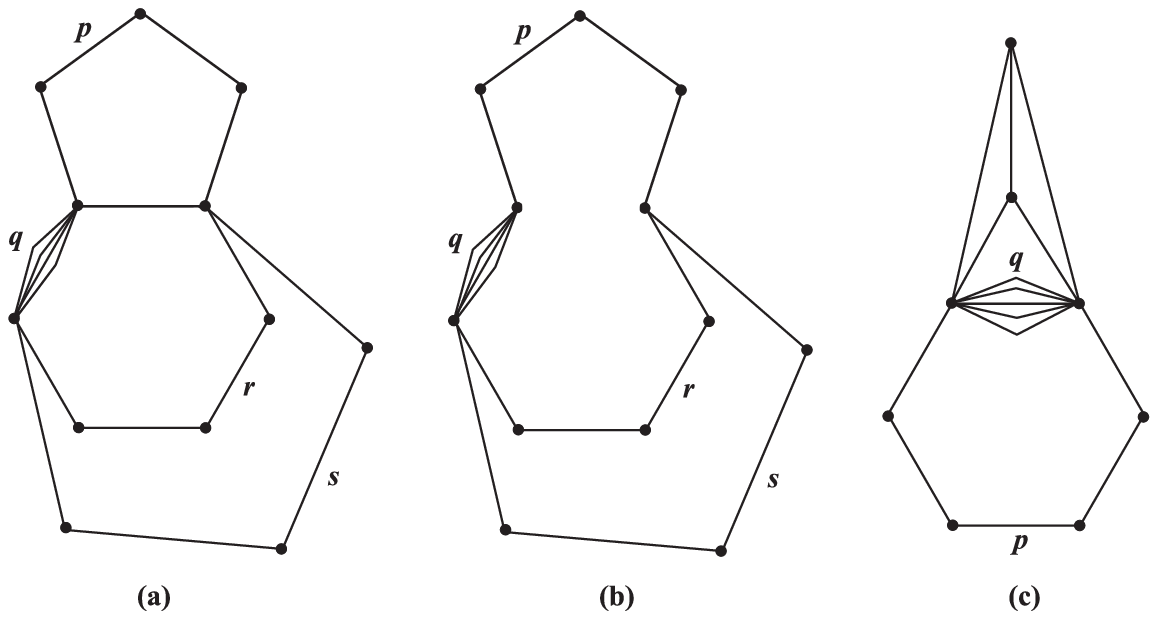,width=3.00in}} \vspace*{8pt}
\caption{Graphs (a) $G((p\,1,q)\,(r,s))$; (b) $G'(p\,\overline
q\,r\,s)$; (c) $G(.p\,q)$.\label{f1.36}}
\end{figure}

Graph $G((p\,1,q)\,(r,s))$ (Fig. 36a) resolves into the graph
$G'(p\,\overline q\,r\,s)$ (Fig. 36b) and the block sum of graphs
$G(s\,q\,r)$, and $G(p)$. The Tutte polynomial of the graphs
$G'(p\,\overline q\,r\,s)$ is

$$T(G'(p\,\overline
q\,r\,s))={{y^q-1}\over {y-1}}T(G(p,r,s)+x^p({{x^{r+s}-1}\over
{x-1}}+y-1),$$

\noindent and the general formula for the Tutte polynomial of the
graphs $G((p\,1,q)\,(r,s))$ is

$$T(G((p\,1,q)\,(r,s)))=T(G'(p\,\overline
q\,r\,s))+T(G(s\,q\,r))T(G(p)).$$

\subsection{Family \it{\textbf{.p\,q}}}

Graph $G(.p\,q)$ (Fig. 36c) resolves into the graphs $G(p\,q,2,2)$
and $G(p,2,2+q)$. The general formula for the Tutte polynomial of
the graphs $G(.p\,q)$ is

$$T(G(.p\,q))=T(G(p\,q,2,2))+T(G(p,2,2+q)).$$

\subsection{Family \it{\textbf{.p\,1\,1}}}

\begin{figure}[th]
\centerline{\psfig{file=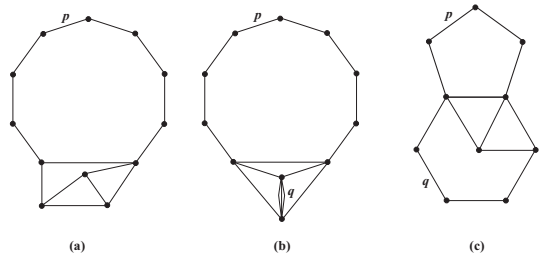,width=3.00in}} \vspace*{8pt}
\caption{Graphs (a) $G(.p\,1\,1)$; (b) $G(.p\,1:q)$; (c)
$G(.p\,1.q)$.\label{f1.37}}
\end{figure}

Graph $G(.p\,1\,1)$ (Fig. 37a) resolves into the graph $G(.p\,1)$
and the block sum of the graph $G((p+1))$ and the graph which
consists from two triangles with the common edge. The general
formula for the Tutte polynomial of the graphs $G(.p\,1\,1)$ is

$$T(G(.p\,1\,1))=T(G(.p\,1))+T(G((p+1)))(x + 2x^2+x^3+y+2xy+y^2).$$

\subsection{Family \it{\textbf{.p\,1:q}}}

Graph $G(.p\,1:q)$ (Fig. 37b) resolves into the graph $G(.p:q\,0)$
and the block sum of the graphs $G(p)$ and $\overline G(q\,2\,2)$.
The general formula for the Tutte polynomial of the graphs
$G(.p\,1:q)$ is

$$T(G(.p\,1:q))=T(G(.p:q\,0))+T(G(p))T(G(q\,2\,2)).$$

\subsection{Family \it{\textbf{.p\,1.q}}}

Graph $G(.p\,1.q)$ (Fig. 37c) resolves into the graphs
$G(p\,1,q\,1,2)$ and $G(p,(q+1),2,1)$. The general formula for the
Tutte polynomial of the graphs $G(.p\,1.q)$ is

$$T(G(.p\,1.q))=T(G(p\,1,q\,1,2))+T(G(p,(q+1),2,1)).$$

\subsection{Family \it{\textbf{.p\,1:q\,0}}}

\begin{figure}[th]
\centerline{\psfig{file=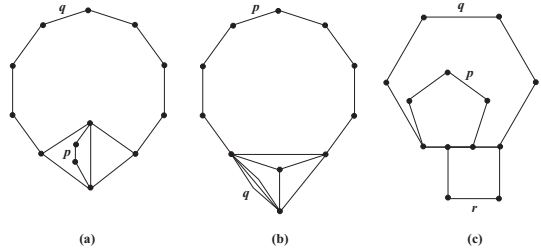,width=3.00in}} \vspace*{8pt}
\caption{Graphs (a) $G(.p\,1:q\,0)$; (b) $G(.p\,1.q\,0)$; (c)
$G(.p.q\,0.r)$.\label{f1.38}}
\end{figure}

Graph $G(.p\,1:q\,0)$ (Fig. 38a) resolves into the graph $G(.p:q)$
and the block sum of graphs $G(q,\overline 2,\overline 2)$ and
$G(p)$. The Tutte polynomial of the graph $G(q,\overline 2,\overline
2)$ is

$$T(G(q,\overline 2,\overline 2))={{x^q-1}\over {x-1}}(x+y)^2+{{y^4-1}\over {y-1}}+x-1.$$

\noindent The general formula for the Tutte polynomial of the graphs
$G(.p\,1:q\,0)$ is

$$T(G(.p\,1:q\,0))=T(G(.p:q))+T(G(q,\overline 2,\overline 2))T(G(p)).$$

\subsection{Family \it{\textbf{.p\,1.q\,0}}}

Graph $G(.p\,1.q\,0)$ (Fig. 38b) resolves into the graph $G(.p.q)$
and the block sum of the graphs $G(1 \overline {(q+1)} \overline 2)$
and $G(p)$. The Tutte polynomial of the graph $G(1\,\overline
{(q+1)}\, \overline 2)$ is

$$T(G(1\,\overline {(q+1)}\, \overline 2))=({{y^{q+1}-1}\over {y-1}}+x-1)(x+y)+{{y^{q+3}-1}\over {y-1}}+x-1.$$

\noindent The general formula for the Tutte polynomial of the graphs
$G(.p\,1.q\,0)$ is

$$T(G(.p\,1.q\,0))=T(G(.p.q))+T(G(1\,\overline {(q+1)}\,
\overline 2))T(G(p)).$$

\subsection{Family \it{\textbf{.p.q\,0.r}}}

Graph $G(.p.q\,0.r)$ (Fig. 38c) resolves into the graph
$G((p+1),q,(r+1))$ and the other graph which resolves into the graph
$G((p+q),r,1)$ and the block sum of graphs $G(q,r,1)$ and $G(p)$.
The general formula for the  Tutte polynomial of the graphs
$G(.p.q\,0.r)$ is

$$T(G(.p.q\,0.r))=T(G((p+1),q,(r+1)))+T(G((p+q),r,1))+T(G(q,r,1))T(G(p)).$$

\subsection{Family \it{\textbf{p:q\,0:r\,0}}}

\begin{figure}[th]
\centerline{\psfig{file=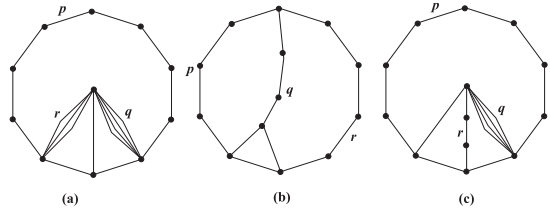,width=3.00in}} \vspace*{8pt}
\caption{Graphs (a) $G(p:q\,0:r\,0)$; (b) $G(p\,0:q\,0:r\,0)$; (c)
$G(.p.q.r)$.\label{f1.39}}
\end{figure}

Graph $G(p:q\,0:r\,0)$ (Fig. 39a) resolves into the graphs
$G(p,\overline {q+1},\overline{r+1})$ and $G((p,2)\,(q,r))$. The
general formula for the  Tutte polynomial of the graphs
$G(p:q\,0:r\,0)$ is

$$T(G(p:q\,0:r\,0))=T(G(p,\overline {q+1},\overline{r+1}))+T(G((p,2)\,(q,r))).$$

\subsection{Family \it{\textbf{p\,0:q\,0:r\,0}}}

Graph $G(p\,0:q\,0:r\,0)$ (Fig. 39b) resolves into the graph
$G((p+1),q,(r+1))$ and the other graph which resolves into the graph
$G(p,(q+1),r)$ and the block sum of the graphs $G(p,q,r)$ and
$G(\overline 1)$. The general formula for the  Tutte polynomial
forthe graphs $G(p\,0:q\,0:r\,0)$ is

$$T(G(p\,0:q\,0:r\,0))=T(G((p+1),q,(r+1)))+T(G(p,(q+1),r))+yT(G(p,q,r)).$$

\subsection{Family \it{\textbf{.p.q.r}}}

Graph $G(.p.q.r)$ (Fig. 39c) resolves into the graphs
$G((p+1)\,1\,r\,q)$ and $G(p\,q\,1\,1\,r)$. The general formula for
the Tutte polynomial of the graphs $G(.p.q.r)$ is

$$T(G(.p.q.r))=T(G((p+1)\,1\,r\,q))+T(G(p\,q\,1\,1\,r)).$$

\subsection{Family \it{\textbf{p:q:r}}}

\begin{figure}[th]
\centerline{\psfig{file=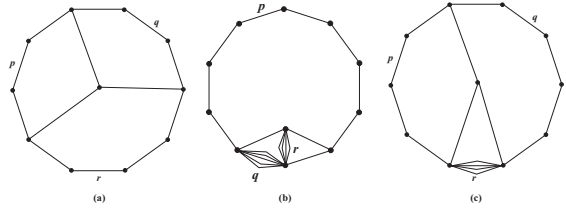,width=3.00in}} \vspace*{8pt}
\caption{Graphs (a) $G(p:q:r)$; (b) $G(.p.q.r\,0)$; (c)
$G(p:q:r\,0)$.\label{f1.40}}
\end{figure}

Graph $G(p:q:r)$ (Fig. 40a) resolves into the graphs $G((p+r),2,q)$
and $G(p\,1\,q\,1\,r)$. The general formula for the Tutte polynomial
of the graphs $G(p:q:r)$ is

$$T(G(p:q:r))=T(G((p+r),2,q))+T(G(p\,1\,q\,1\,r)).$$

\subsection{Family \it{\textbf{.p.q.r\,0}}}

Graph $G(.p.q.r\,0)$ (Fig. 40b) resolves into the graphs
$G((p+1)\,1,q,r)$ and $G(p\,q,(r+1),1)$. The general formula for the
Tutte polynomial of the graphs $G(.p.q.r\,0)$ is

$$T(G(.p.q.r\,0))=T(G((p+1)\,1,q,r))+T(G(p\,q,(r+1),1)).$$

\subsection{Family \it{\textbf{p:q:r\,0}}}

Graph $G(p:q:r\,0)$ (Fig. 40c) resolves into the graphs
$G((p+q)\,r\,2)$ and $G(p\,1,q\,1,r)$. The general formula for the
Tutte polynomial of the graphs $G(p:q:r\,0)$ is

$$T(G(p:q:r\,0))=T(G((p+q)\,r\,2))+T(G(p\,1,q\,1,r)).$$

\subsection{Family \it{\textbf{.(p,q)}}}

\begin{figure}[th]
\centerline{\psfig{file=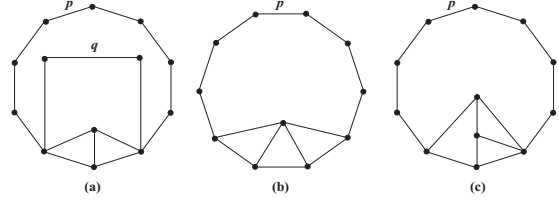,width=3.00in}} \vspace*{8pt}
\caption{Graphs (a) $G(.(p,q))$; (b) $G(8^*p)$; (c)
$G(8^*p\,0)$.\label{f1.41}}
\end{figure}

Graph $G(.(p,q))$ (Fig. 41a) resolves into the graphs $G(p,2,2,q)$
and $G((p,q)\,(2,2))$. The general formula for the Tutte polynomial
of the graphs $G(.(p,q))$ is

$$T(G(.(p,q)))=T(G(p,2,2,q))+T(G((p,q)\,(2,2))).$$

\subsection{Family \it{\textbf{8$^*$p}}}

Graph $G(8^*p)$ (Fig. 41b) resolves into the graphs
$G(2\,1\,p\,1\,2)$ and $G(.p:2\,0)$. The general formula for the
Tutte polynomial of the graphs $G(8^*p)$ is

$$T(G(8^*p))=T(G(2\,1\,p\,1\,2))+T(G(.p:2\,0)).$$

\subsection{Family \it{\textbf{8$^*$p\,0}}}

Graph $G(8^*p\,0)$ (Fig. 41c) resolves into the graphs
$G((p+1)\,1\,1\,1\,2)$ and $G(.p\,1:1)$. The general formula for the
Tutte polynomial of the graphs $G(8^*p\,0)$ is

$$T(G(8^*p\,0))=T(G((p+1)\,1\,1\,1\,2))+T(G(.p\,1)).$$

\subsection{Family \it{\textbf{p\,q\,r\,s\,t}}}

\begin{figure}[th]
\centerline{\psfig{file=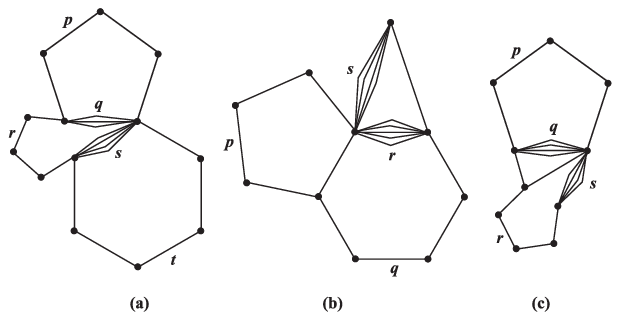,width=3.00in}} \vspace*{8pt}
\caption{Graphs (a) $G(p\,q\,r\,s\,t)$; (b) $G(p\,1\,q\,r\,1\,s)$;
(c) $G(p\,q\,1\,1\,r\,s)$.\label{f1.42}}
\end{figure}

The graphs $G(p\,q\,r\,s\,t)$ illustrated in Fig. 42a  correspond to
the link family \emph{\textbf{p\,q\,r\,s\,t}}. In order to obtain
general formula for the Tutte polynomial we use the relations

$$T(G(p\,q\,r\,s\,t))-T(G(p\,(q-1)\,r\,s\,t))=y^{q-1}T(G(p))T(G(r\,s\,t)).$$

\noindent The general formula for the Tutte polynomial of the graphs
$G(p\,q\,r\,s\,t)$ is

$$T(G(p\,q\,r\,s\,t))={{y^q-1}\over{y-1}}T(G(p))T(G(r\,s\,t))+T(G((p+r)\,s\,t)).$$

\subsection{Family \it{\textbf{p\,1\,q\,r\,1\,s}}}

Graph $G(p\,1\,q\,r\,1\,s)$ (Fig. 42b) resolves into the graph
$G((p+q)\,r\,1\,s)$ and the block sum of the graphs $G(q\,r\,1\,s)$
and $G(p)$.  The general formula for the Tutte polynomial of the
graphs $G(p\,1\,q\,r\,1\,s)$ is

$$T(G(p\,1\,q\,r\,1\,s))=T(G((p+q)\,r\,1\,s))+T(G(q\,r\,1\,s))T(G(p)).$$

\subsection{Family \it{\textbf{p\,q\,1\,1\,r\,s}}}

The graph $G(p\,q\,1\,1\,r\,s)$ illustrated in Fig. 42c corresponds
to the link family \emph{\textbf{p\,q\,1\,1\,r\,s}}. The graph
$G(p\,q\,1\,1\,r\,s)$ resolves into the graph $G'(p\,(\overline
{q+1})\,r\,\overline s)$ and the block sum of the graphs $G(p\,q)$
and $G((r+1)\,s)$. In order to obtain general formula for the Tutte
polynomial of the graph $G'(p\,\overline q\,r\,\overline s)$ (Fig.
42d) we use the relations

$$T(G'(p\,\overline q\,r\,\overline s))-T(G'(p\,(\overline {q-1})\,r\,\overline s))=y^{q-1}T(G(p))T(G(r\,s)).$$

\noindent The Tutte polynomial of the graphs $G'(p\,\overline
q\,r\,\overline s)$ is

$$T(G'(p\,\overline q\,r\,\overline s))={{y^q-1}\over{y-1}}T(G(p))T(G(r\,s))+T(G((p+r)\,s)),$$

\noindent and the general formula for the Tutte polynomial of the
graphs $G(p\,q\,1\,1\,r\,s)$ is

$$T(G(p\,q\,1\,1\,r\,s))=T(G'(p\,(\overline
{q+1})\,r\,\overline s))+T(G(p\,q))T(G((r+1)\,s)).$$

\subsection{Family \it{\textbf{p\,q\,1\,r\,1\,s}}}

\begin{figure}[th]
\centerline{\psfig{file=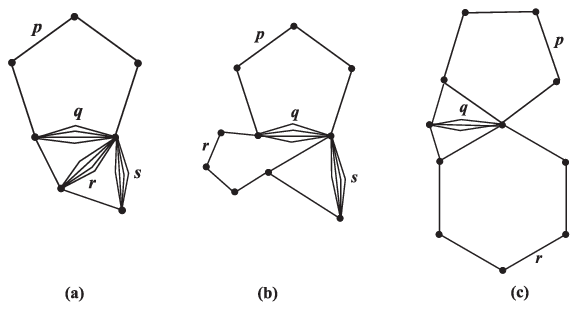,width=3.00in}} \vspace*{8pt}
\caption{Graphs (a) $G(p\,q\,1\,r\,1\,s)$; (b)
$G(p\,q\,r\,1\,1\,s)$; (c) $G(p\,1\,1\,q\,1\,1\,r)$.\label{f1.43}}
\end{figure}

Graph $G(p\,q\,1\,r\,1\,s)$ (Fig. 43a) resolves into the graph
$G(p\,q\,1\,(r+s))$ and the block sum of the graphs $G(p\,q\,1\,r)$
and $G(\overline s)$. The general formula for the Tutte polynomial
of the graphs $G(p\,q\,1\,r\,1\,s)$ is

$$T(G(p\,q\,1\,r\,1\,s))=T(G(p\,q\,1\,(r+s)))+T(G(p\,q\,1\,r))T(G(\overline s)).$$

\subsection{Family \it{\textbf{p\,q\,r\,1\,1\,s}}}

Graph $G(p\,q\,r\,1\,1\,s)$ (Fig. 43b) resolves into the graph
$G'(p\,\overline q\,(r+1)\,\overline s)$ and the block sum of the
graphs $G(p\,q\,r)$ and $G((\overline {s+1}))$. The general formula
for the Tutte polynomial of the graphs $G(p\,q\,r\,1\,1\,s)$ is

$$T(G(p\,q\,r\,1\,1\,s))=T(G'(p\,\overline q\,(r+1)\,\overline s))+T(G(p\,q\,r))T(G((\overline {s+1}))).$$

\subsection{Family \it{\textbf{p\,1\,1\,q\,1\,1\,r}}}

Graph $G(p\,1\,1\,q\,1\,1\,r)$ (Fig. 43c) resolves into the graph
$G((p+1)\,q\,1\,1\,r)$ and the block sum of the graphs
$G(r\,1\,1\,(q+1))$ and $G(p)$. The general formula for the Tutte
polynomial of the graphs  $G(p\,1\,1\,q\,1\,1\,r)$  is

$$T(G(p\,1\,1\,q\,1\,1\,r))=T(G((p+1)\,q\,1\,1\,r))+T(G(r\,1\,1\,(q+1)))T(G(p)).$$

\subsection{Family \it{\textbf{p\,1\,q\,1\,1\,1\,r}}}

\begin{figure}[th]
\centerline{\psfig{file=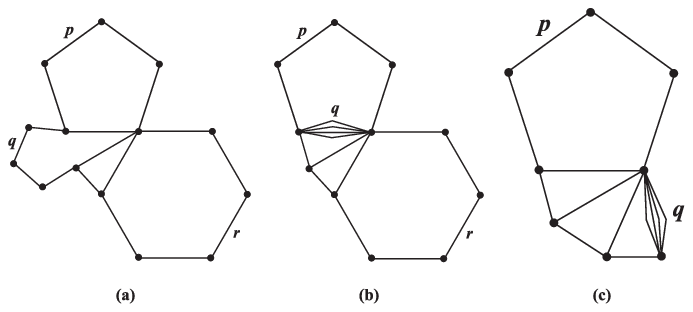,width=3.40in}} \vspace*{8pt}
\caption{Graphs (a) $G(p\,1\,q\,1\,1\,1\,r)$; (b)
$G(p\,q\,1\,1\,1\,1\,r)$; (c)
$G(p\,1\,1\,1\,1\,1\,1\,r)$.\label{f1.44}}
\end{figure}

Graph $G(p\,1\,q\,1\,1\,1\,r)$ (Fig. 44a) resolves into the graph
$G((p+q)\,1\,1\,1\,r)$ and the block sum of the graphs
$G(q\,1\,1\,1\,r)$ and $G(p)$. The general formula for the Tutte
polynomial of the graphs  $G(p\,1\,q\,1\,1\,1\,r)$ is

$$T(G(p\,1\,q\,1\,1\,1\,r))=T(G((p+q)\,1\,1\,1\,r))+T(G(q\,1\,1\,1\,r))T(G(p)).$$

\subsection{Family \it{\textbf{p\,q\,1\,1\,1\,1\,r}}}

Graph $G(p\,q\,1\,1\,1\,1\,r)$ (Fig. 44b) resolves into the graph
$G(p\,q\,1\,1\,(r+1))$ and the block sum of the graphs
$G(p\,q\,1\,1\,1)$ and $G(r)$. The general formula for the Tutte
polynomial of the graphs  $G(p\,q\,1\,1\,1\,1\,r)$  is

$$T(G(p\,q\,1\,1\,1\,1\,r))=T(G(p\,q\,1\,1\,(r+1)))+T(G(p\,q\,1\,1\,1))T(G(r)).$$

\subsection{Family \it{\textbf{p\,1\,1\,1\,1\,1\,1\,q}}}

Graph $G(p\,1\,1\,1\,1\,1\,1\,q)$ (Fig. 44c) resolves into the graph
$G((p+1)\,1\,1\,1\,1\,q)$ and the block sum of the graphs $\overline
G(2\,1\,1\,1\,q)$ and $G(p)$, where $\overline G(2\,1\,1\,1\,q)$ is
the dual of the graph $G(2\,1\,1\,1\,q)$. The general formula for
the Tutte polynomial of the graphs $G(p\,1\,1\,1\,1\,1\,1\,q)$  is

$$T(G(p\,1\,1\,1\,1\,1\,1\,q))=T(G((p+1)\,1\,1\,1\,1\,q))+T(\overline G(2\,1\,1\,1\,q))T(G(p)).$$

\subsection{Family \it{\textbf{p,q,r,s,t}}}

\begin{figure}[th]
\centerline{\psfig{file=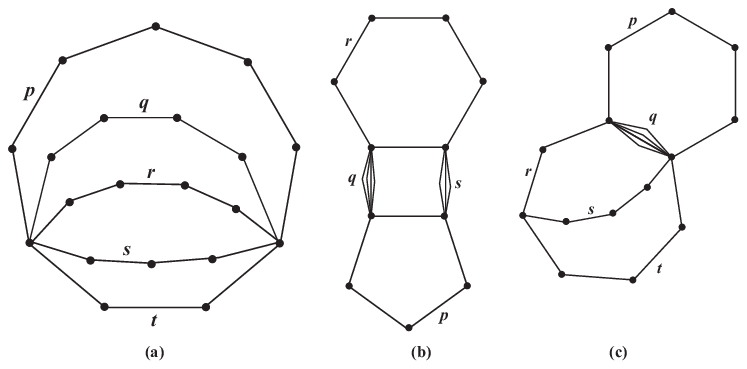,width=3.40in}} \vspace*{8pt}
\caption{Graphs (a) $G(p,q,r,s,t)$; (b) $G(p\,1,q,r\,1,s)$; (c)
$G(p\,q\,r,s,t)$.\label{f1.45}}
\end{figure}

In order to obtain formula for the Tutte polynomial of the graph
$G(p,q,r,s,t)$ (Fig. 45a) we use the relations

$$T(G(p,q,r,s,t))-T(G(p,q,r,s,(t-1)))=x^{t-1}T(G(p,q,r,s)).$$

\noindent Since the Tutte polynomial of the graph $G(p,q,r,s,0)$ is
$$T(G(p,q,r,s,0))=T(G(p))T(G(q))T(G(r))T(G(s)),$$ the general formula
for the Tutte polynomial of the graphs  $G(p,q,r,s,t)$ is

$$T(G(p,q,r,s,t))={{x^t-1}\over
{x-1}}T(G(p,q,r,s))+T(G(p))T(G(q))T(G(r))T(G(s)).$$

\subsection{Family \it{\textbf{p\,1,q,r\,1,s}}}

In order to obtain formula for the Tutte polynomial of the graph
$G(p\,1,q,r\,1,s)$ (Fig. 45b) we use the relations

$$T(G(p\,1,q,r\,1,s))-T(G(p\,1,q,r\,1,(s-1)))=y^{s-1}T(G(p\,1,r\,1,q)).$$

\noindent Since the Tutte polynomial of the graph $G(p\,1,q,r\,1,0)$
is $T(G(p\,1,q,r\,1,0))=T(G(p+1))T(G(\overline q))T(G(r+1))$, the
general formula for the Tutte polynomial of the graphs
$G(p\,1,q,r\,1,s)$ is

$$T(p\,1,q,r\,1,s)={{y^s-1}\over
{y-1}}T(G(p\,1,r\,1,q))+T(G(p+1))T(G(\overline q))T(G(r+1)).$$

\subsection{Family \it{\textbf{p\,q\,r,s,t}}}

In order to obtain formula for the Tutte polynomial of the graph
$G(p\,q\,r,s,t)$ (Fig. 45c) we use the relations

$$T(G(p\,q\,r,s,t))-T(G(p\,q\,r,s,(t-1)))=x^{t-1}T(G(p\,q\,(r+s))).$$

\noindent Since the Tutte polynomial of the graph $G(p\,q\,r,s,0)$
is $T(G(p\,q\,r,s,0))=T(G(p\,q\,r))T(G(s))$, the general formula for
the Tutte polynomial of the graphs $G(p\,q\,r,s,t)$ is

$$T(G(p\,q\,r,s,t))={{x^t-1}\over
{x-1}}T(G(p\,q\,(r+s)))+T(G(p\,q\,r))T(G(s)).$$

\subsection{Family \it{\textbf{p\,1\,1\,q,r,s}}}

\begin{figure}[th]
\centerline{\psfig{file=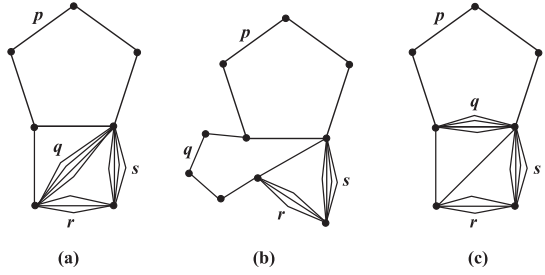,width=3.00in}} \vspace*{8pt}
\caption{Graphs (a) $G(p\,1\,1\,q,r,s)$; (b) $G(p\,1\,q\,1,r,s)$;
(c) $G(p\,q\,1\,1,r,s)$.\label{f1.46}}
\end{figure}

Graph $G(p\,1\,1\,q,r,s)$ (Fig. 46a) resolves into the graph
$G((p+1)\,q,r,s)$ and the block sum of the graphs $\overline
G((q+1)\,r\,s)$ and $G(p)$, where $\overline G((q+1)\,r\,s)$ is the
dual of the graph $G((q+1)\,r\,s)$. The general formula for the
Tutte polynomial of the graphs $G(p\,1\,1\,q,r,s)$  is

$$T(G(p\,1\,1\,q,r,s))=T(G((p+1)\,q,r,s))+T(\overline
G((q+1)\,r\,s))T(G(p)).$$

\subsection{Family \it{\textbf{p\,1\,q\,1,r,s}}}

Graph $G(p\,1\,q\,1,r,s)$ (Fig. 46b) resolves into the graph
$G((p+q)\,1,r,s)$ and the block sum of the graphs $G(q\,1,r,s)$ and
$G(p)$. The general formula for the Tutte polynomial of the graphs
$G(p\,1\,q\,1,r,s)$ is

$$T(G(p\,1\,q\,1,r,s))=T(G((p+q)\,1,r,s))+T(G(q\,1,r,s))T(G(p)).$$

\subsection{Family \it{\textbf{p\,q\,1\,1,r,s}}}

Graph $G(p\,q\,1\,1,r,s)$ (Fig. 46c) resolves into the graph
$G(p\,(q+1),s,r)$ and the block sum of the graphs $\overline
G(s\,1\,r)$ and $G(p\,q)$. The general formula for the Tutte
polynomial of the graphs $G(p\,q\,1\,1,r,s)$ is

$$T(G(p\,q\,1\,1,r,s))=T(G(p\,(q+1),s,r))+T(\overline G(s\,1\,r))T(G(p\,q)).$$

\subsection{Family \it{\textbf{p\,1\,1\,1\,1,q,r}}}

\begin{figure}[th]
\centerline{\psfig{file=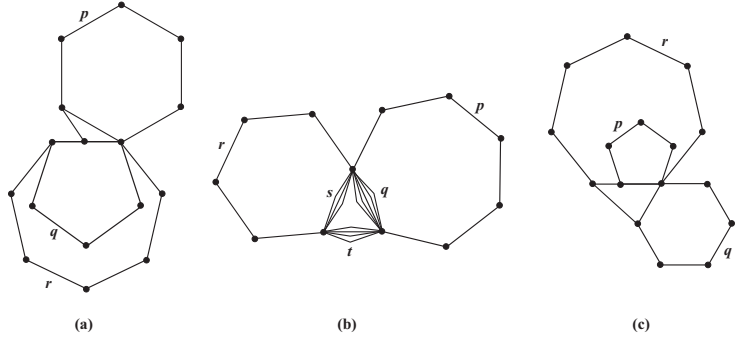,width=4.00in}} \vspace*{8pt}
\caption{Graphs (a) $G(p\,1\,1\,1\,1,q,r)$; (b) $G(p\,q,r\,s,t)$;
(c) $G(p\,1\,1,q\,1\,1,r)$.\label{f1.47}}
\end{figure}

Graph $G(p\,1\,1\,1\,1,q,r)$ (Fig. 47a) resolves into the graph
$G(p\,1\,1,q,r+)$ and the block sum of the graphs $G(p\,1\,2)$ and
$G(q+r)$. The general formula for the Tutte polynomial of the graphs
$G(p\,1\,1\,1\,1,q,r)$ is

$$T(p\,1\,1\,1\,1,q,r)=T(G(p\,1\,1,q,r+))+T(G(p\,1\,2))T(G(q+r)).$$

\subsection{Family \it{\textbf{p\,q,r\,s,t}}}

In order to obtain formula for the Tutte polynomial of the graph
$G(p\,q,r\,s,t)$ (Fig. 47b) we use the relations

$$T(G(p\,q,r\,s,t))-T(G(p\,q,r\,s,(t-1)))=y^{t-1}T(G(p\,(q+s)\,r)).$$

\noindent Since the Tutte polynomial of the graph $G(p\,q,r\,s,0)$
is $T(G(p\,q,r\,s,0))=T(G(p\,q))T(G(r\,s))$, the general formula for
the Tutte polynomial of the graphs $G(p\,q,r\,s,t)$ is

$$T(G(p\,q,r\,s,t))={{y^t-1}\over
{y-1}}T(G(p\,(q+s)\,r))+T(G(p\,q))T(G(r\,s)).$$

\subsection{Family \it{\textbf{p\,1\,1,q\,1\,1,r}}}

Graph $G(p\,1\,1,q\,1\,1,r)$ (Fig. 47c) resolves into the graph
$G(q\,1\,1,p,r+)$ and the block sum of the graphs $G(q\,1\,(r+1))$
and $G(p+1)$. The general formula for the Tutte polynomial of the
graphs $G(p\,1\,1,q\,1\,1,r)$ is

$$T(G(p\,1\,1,q\,1\,1,r))=T(G(q\,1\,1,p,r+))+T(G(q\,1\,(r+1)))T(G(p+1)).$$

\subsection{Family \it{\textbf{p\,q,r\,1,s\,1}}}

\begin{figure}[th]
\centerline{\psfig{file=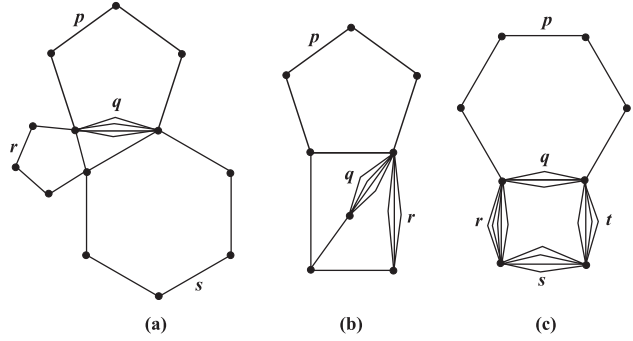,width=3.40in}} \vspace*{8pt}
\caption{Graphs (a) $G(p\,q,r\,1,s\,1)$; (b) $G(p\,1\,1,q\,1,r\,1)$;
(c) $G(p\,q,r,s,t)$.\label{f1.48}}
\end{figure}

Graph $G(p\,q,r\,1,s\,1)$ (Fig. 48a) resolves into the graph
$G(p\,q\,r\,1\,s)$ and the block sum of the graphs $G(p\,(q+1)\,s)$
and $G(r)$. The general formula for the Tutte polynomial of the
graphs $G(p\,q,r\,1,s\,1)$ is

$$T(G(p\,q,r\,1,s\,1))=T(G(p\,q\,r\,1\,s))+T(G(p\,(q+1)\,s))T(G(r)).$$

\subsection{Family \it{\textbf{p\,1\,1,q\,1,r\,1}}}

Graph $G(p\,1\,1,q\,1,r\,1)$ (Fig. 48b) resolves into the graph
$G(p\,1\,1\,(q+r))$, the block sum of graphs $G(p\,1\,1\,q)$ and
$G(\overline r)$, and the block sum of graphs $G(p\,1\,2\,r)$ and
$G(\overline q)$. The general formula for the Tutte polynomial of
the graphs $G(p\,1\,1,q\,1,r\,1)$ is

$$T(G(p\,1\,1,q\,1,r\,1))=T(G(p\,1\,1\,(q+r))+T(G(p\,1\,1\,q))T(G(\overline r))+T(G(p\,1\,2\,r))T(G(\overline q)).$$

\subsection{Family \it{\textbf{p\,q,r,s,t}}}

In order to obtain formula for the Tutte polynomial of the graph
$G(p\,q,r,s,t)$ (Fig. 48c) we use the relations

$$T(G(p\,q,r,s,t))-T(G(p\,q,r,s,(t-1)))=y^{t-1}T(G(p\,q,r,s)).$$

\noindent Since the Tutte polynomial of the graph $G(p\,q,r,s,0)$ is

$$T(G(p\,q,r,s,0))=T(G(\overline r))T(G(\overline s))T(G(p\,q)),$$

\noindent the general formula for the Tutte polynomial of the graphs
$G(p\,q,r,s,t)$ is

$$T(G(p\,q,r,s,t))={{y^t-1}\over
{y-1}}T(G(p\,q,r,s))+T(G(\overline r))T(G(\overline s))T(G(p\,q)).$$

\subsection{Family \it{\textbf{p\,1\,1,q,r,s}}}

\begin{figure}[th]
\centerline{\psfig{file=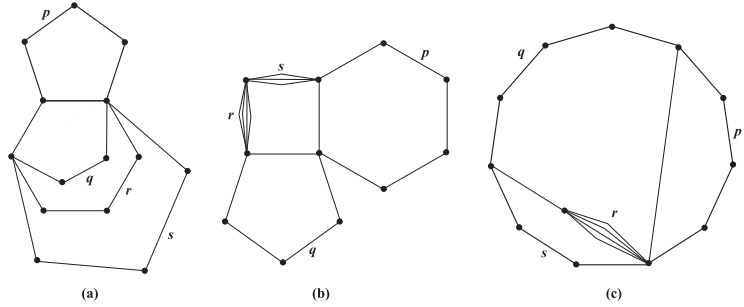,width=4.00in}} \vspace*{8pt}
\caption{Graphs (a) $G(p\,1\,1,q,r,s)$; (b) $G(p\,1,q\,1,r,s)$; (c)
$G(p\,1\,q,r\,1,s)$.\label{f1.49}}
\end{figure}

Graph $G(p\,1\,1,q,r,s)$ (Fig. 49a) resolves into the graph
$G((p+1),q,r,s)$, the block sum of graphs $G(q,r,s)$ and $G(p)$, and
the block sum of graphs $G(p)$, $G(q)$, $G(r)$, and $G(s)$. The
general formula for the Tutte polynomial of the graphs
$G(p\,1\,1,q,r,s)$ is

$$T(G(p\,1\,1,q,r,s))=T(G((p+1),q,r,s))+T(G(q,r,s))T(G(p))+$$

$$T(G(p))T(G(q))T(G(r))T(G(s)).$$

\subsection{Family \it{\textbf{p\,1,q\,1,r,s}}}

Graph $G(p\,1,q\,1,r,s)$ (Fig. 49b) resolves into the graph
$G((p+q)\,\overline r\,\overline s)$, the block sum of graphs
$G(p\,1,r,s)$ and $G(q)$, and the block sum of graphs
$G(q\,\overline r\,\overline s)$ and $G(p)$. The general formula for
the Tutte polynomial of the graphs $G(p\,1\,1,q,r,s)$ is

$$T(G(p\,1,q\,1,r,s))=T(G((p+q)\,\overline r\,\overline s))+T(G(p\,1,r,s))T(G(q))+T(G(q\,\overline r\,\overline s))T(G(p)).$$

\subsection{Family \it{\textbf{p\,1\,q,r\,1,s}}}

Graph $G(p\,1\,q,r\,1,s)$ (Fig. 49c) resolves into the graph
$G(s\,r\,q\,1\,p)$ and the block sum of the graphs $G(s+q,1,p)$ and
$\overline G(r)$. The general formula for the Tutte polynomial of
the graphs $G(p\,1\,q,r\,1,s)$ is

$$T(G(p\,1\,q,r\,1,s))=T(G(s\,r\,q\,1\,p))+T(G(s+q,1,p))T(\overline G(r)).$$

\subsection{Family \it{\textbf{p\,q\,1,r\,1,s}}}

\begin{figure}[th]
\centerline{\psfig{file=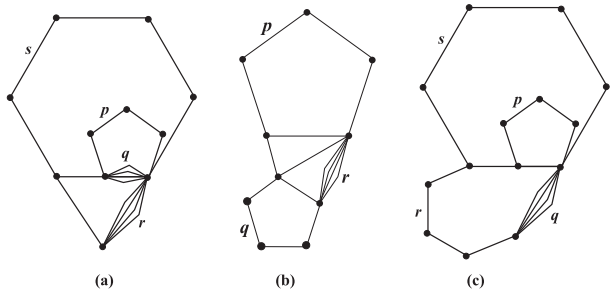,width=3.40in}} \vspace*{8pt}
\caption{Graphs (a) $G(p\,q\,1,r\,1,s)$; (b) $G(p\,1\,1\,1,q\,1,r)$;
(c) $G(p\,1\,1,q\,r,s)$.\label{f1.50}}
\end{figure}

Graph $G(p\,q\,1,r\,1,s)$ (Fig. 50a) resolves into the graph
$G(p\,(q+r)\,s)$, the block sum of graphs $G(p\,q)$ and
$G((s+1)\,r)$, and the block sum of graphs $G(p\,q\,s)$ and
$G(\overline r)$. The general formula for the Tutte polynomial of
the graphs $G(p\,q\,1,r\,1,s)$ is

$$T(G(p\,q\,1,r\,1,s))=T(G(p\,(q+r)\,s))+T(G(p\,q))T(G((s+1)\,r))+T(G(p\,q\,s))T(G(\overline r)).$$

\subsection{Family \it{\textbf{p\,1\,1\,1,q\,1,r}}}

Graph $G(p\,1\,1\,1,q\,1,r)$ (Fig. 50b) resolves into the graph
$G((p+1)\,1,q\,1,r)$ and the block sum of graphs $G(q\,1,2,r)$ and
$G(p)$. The general formula for the Tutte polynomial of the graphs
$G(p\,1\,1\,1,q\,1,r)$ is

$$T(G(p\,1\,1\,1,q\,1,r))=T(G((p+1)\,1,q\,1,r)+T(G(q\,1,2,r))T(G(p)).$$

\subsection{Family \it{\textbf{p\,1\,1,q\,r,s}}}

Graph $G(p\,1\,1,q\,r,s)$ (Fig. 50c) resolves into graph $\overline
G(q\,r,s,(p+1))$, dual of the graph $G(q\,r,s,(p+1))$, the block sum
of graphs $G((r+s)\,q)$ and $G(p)$, and the block sum of graphs
$G(r\,q)$, $G(p)$ and $G(s)$. The general formula for the Tutte
polynomial of the graphs $G(p\,1\,1,q\,r,s)$ is

$$T(G(p\,1\,1,q\,r,s))=T(\overline G(q\,r,s,(p+1))+T(G((r+s)\,q))T(G(p))+$$

$$T(G(r\,q))T(G(s))T(G(p)).$$

\subsection{Family \it{\textbf{p,q,r,s+t}}}

\begin{figure}[th]
\centerline{\psfig{file=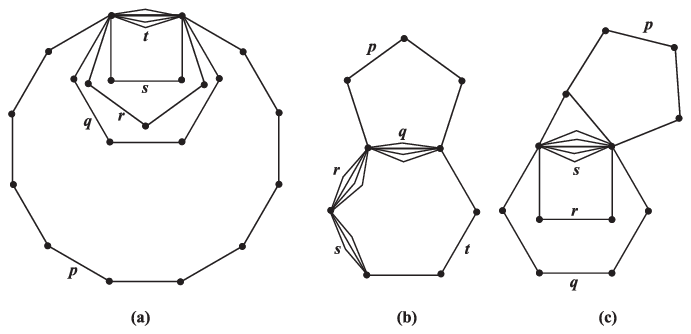,width=3.40in}} \vspace*{8pt}
\caption{Graphs (a) $G(p,q,r,s+t)$; (b) $G(p\,q,r,s+t)$; (c)
$G(p\,1\,1,q,r+s)$.\label{f1.51}}
\end{figure}

In order to obtain formula for the Tutte polynomial of the graph
$G(p,q,r,s+t)$ (Fig. 51a) we use the relations

$$T(G(p,q,r,s+t))-T(G(p,q,r,s+(t-1)))=y^{t-1}T(G(p))T(G(q))T(G(r))T(G(s)).$$

\noindent Since the Tute polynomial of the graph $G(p,q,r,s,0)$ is
$T(G(p,q,r,s,0))=T(G(p,q,r,s))$, the general formula for the Tutte
polynomial of the graphs $G(p,q,r,s+t)$ is

$$T(G(p,q,r,s+t))={{y^t-1}\over
{y-1}}T(G(p))T(G(q))T(G(r))T(G(s))+T(G(p,q,r,s)).$$

\subsection{Family \it{\textbf{p\,q,r,s+t}}}

In order to obtain formula for the Tutte polynomial of the graph
$G(p\,q,r,s+t)$ (Fig. 51b) we use the relations

$$T(G(p\,q,r,s+t))-T(G(p\,q,r,s+(t-1)))=x^{t-1}T(G(p\,q))T(G(\overline r))T(G(\overline s)).$$

\noindent Since the Tutte polynomial of the graph $G(p\,q,r,s+0)$ is
$T(G(p\,q,r,s+0))=T(G(p\,q,r,s))$, the general formula for the Tutte
polynomial of the graphs $G(p\,q,r,s+t)$ is

$$T(G(p\,q,r,s+t))={{x^t-1}\over
{x-1}}T(G(p\,q))T(G(\overline r))T(G(\overline s))+T(G(p\,q,r,s)).$$

\subsection{Family \it{\textbf{p\,1\,1,q,r+s}}}

Graph $G(p\,1\,1,q,r+s)$ (Fig. 51c) resolves into the graph
$G(\overline s\,q\,r\,(p+1))$, dual of the graph $G(s\,\overline
q\,\overline r\,(\overline{p+1}))$ and the block sum of graphs
$G(q\,(s+1)\,r)$ and $G(p)$. The general formula for the Tutte
polynomial of the graphs $G(p\,1\,1,q,r+s)$ is

$$T(G(p\,1\,1,q,r+s))=T(G(\overline
s\,q\,r\,(p+1)))+T(G(q\,(s+1)\,r))T(G(p)).$$

\subsection{Family \it{\textbf{p\,1,q\,1,r+s}}}

\begin{figure}[th]
\centerline{\psfig{file=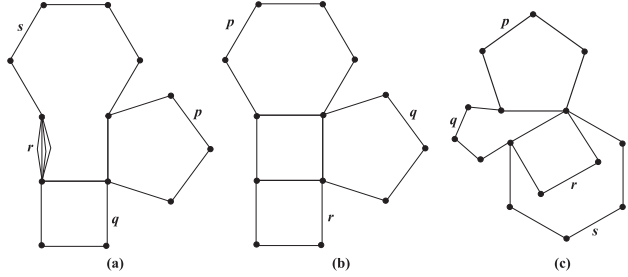,width=3.40in}} \vspace*{8pt}
\caption{Graphs (a) $G(p\,1,q\,1,r+s)$; (b) $G(p\,1,q\,1,r\,1+)$;
(c) $G(p\,1\,q,r,s+)$.\label{f1.52}}
\end{figure}

Graph $G(p\,1,q\,1,r+s)$ (Fig. 52a) resolves into the graph
$G(q\,1,r,1+(p+s-1))$ and the block sum of the graphs
$G(q\,1,r,1+(s-1))$ and $G(p)$. The general formula for the Tutte
polynomial of the graphs $G(p\,1,q\,1,r+s)$ is

$$T(G(p\,1,q\,1,r+s))=T(G(q\,1,r,1+(p+s-1)))+T(G(q\,1,r,1+(s-1)))T(G(p)).$$

\subsection{Family \it{\textbf{p\,1,q\,1,r\,1+}}}

Graph $G(p\,1,q\,1,r\,1+)$ (Fig. 52b) resolves into the graph
$G(p\,1\,(r+1)\,1\,q)$ and the block sum of the graphs
$G(p\,1\,1\,1\,q)$ and $G(r)$. The general formula for the Tutte
polynomial of the graphs $G(p\,1,q\,1,r\,1+)$ is

$$T(G(p\,1,q\,1,r\,1+))=T(G(p\,1\,(r+1)\,1\,q))+T(G(p\,1\,1\,1\,q))T(G(r)).$$

\subsection{Family \it{\textbf{p\,1\,q,r,s+}}}

Graph $G(p\,1\,q,r,s+)$ (Fig. 52c) resolves into the graph
$G((p+q),1,r,s)$ and the block sum of the graphs $G(q,1,r,s)$ and
$G(p)$. The general formula for the Tutte polynomial of the graphs
$G(p\,1\,q,r,s+)$ is

$$T(G(p\,1\,q,r,s+))=T(G((p+q),1,r,s))+T(G(q,1,r,s))T(G(p)).$$

\subsection{Family \it{\textbf{p\,q\,1,r,s+}}}

\begin{figure}[th]
\centerline{\psfig{file=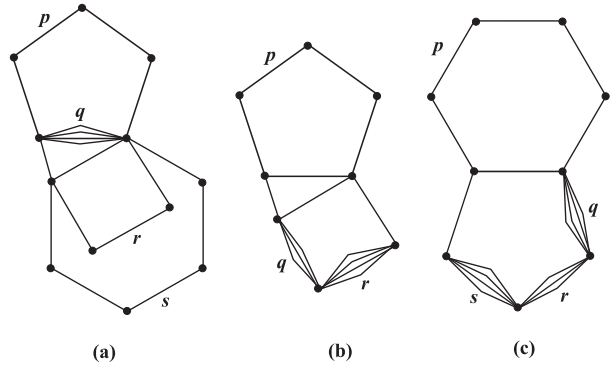,width=3.40in}} \vspace*{8pt}
\caption{Graphs (a) $G(p\,q\,1,r,s+)$; (b) $G(p\,1\,1\,1,q,r+)$; (c)
$G(p\,1,q,r,s+)$.\label{f1.53}}
\end{figure}

Graph $G(p\,q\,1,r,s+)$ (Fig. 53a) resolves into the graph
$G((\overline {q+1})\,r\,s\,p)$, dual of the graph
$G((q+1)\,\overline r\,\overline s\,\overline p)$ and the block sum
of the graphs $G(p\,q)$ and $G(1,r,s)$. The general formula for the
Tutte polynomial of the graphs $G(p\,q\,1,r,s+)$ is

$$T(G(p\,q\,1,r,s+))=T(G((\overline
{q+1})\,r\,s\,p))+T(G(p\,q))T(G(1,r,s)).$$

\subsection{Family \it{\textbf{p\,1\,1\,1,q,r+}}}

Graph $G(p\,1\,1\,1,q,r+)$ (Fig. 53b) resolves into the graph
$G(p\,1\,1\,1,q,r)$ and the block sum of the graphs $G(p\,1\,2)$,
$G(\overline q)$, and $G(\overline r)$. The general formula for the
Tutte polynomial of the graphs $G(p\,1\,1\,1,q,r+)$ is

$$T(G(p\,1\,1\,1,q,r+))=T(G((p\,1\,1\,1,q,r))+T(G(p\,1\,2))T(G(\overline q))T(G(\overline r)).$$

\subsection{Family \it{\textbf{p\,1,q,r,s+}}}

Graph $G(p\,1,q,r,s+)$ (Fig. 53c) resolves into the graph
$G(p\,1,q,r,s)$ and the block sum of the graphs $G(p+1)$,
$G(\overline q)$, $G(\overline r)$, and $G(\overline s)$. The
general formula for the Tutte polynomial of the graphs
$G(p\,1,q,r,s+)$ is

$$T(G(p\,1,q,r,s+))=T(G((p\,1,q,r,s))+T(G(p+1))T(G(\overline q))T(G(\overline r))T(G(\overline s)).$$

\subsection{Family \it{\textbf{p\,q,r\,1,s+}}}

\begin{figure}[th]
\centerline{\psfig{file=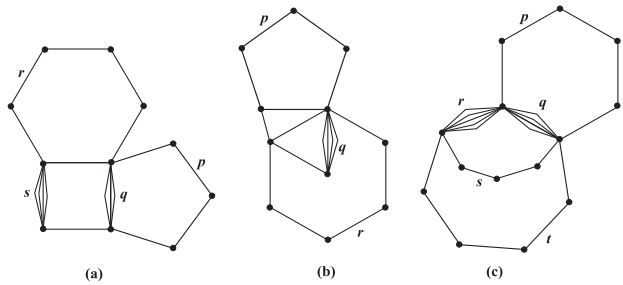,width=3.40in}} \vspace*{8pt}
\caption{Graphs (a) $G(p\,q,r\,1,s+)$; (b) $G(p\,1\,1,q\,1,r+)$; (c)
$G((p\,q,r)\,(s,t))$.\label{f1.54}}
\end{figure}

Graph $G(p\,q,r\,1,s+)$ (Fig. 54a) resolves into the graph
$G(p\,q,r\,1,s)$ and the block sum of the graphs $G(p\,q)$,
$G(r+1)$, and $G(\overline s)$. The general formula for the Tutte
polynomial of the graphs $G(p\,q,r\,1,s+)$ is

$$T(G(p\,q,r\,1,s+))=T(G((p\,q,r\,1,s))+T(G(p\,q))T(G(r+1))T(G(\overline s)).$$

\subsection{Family \it{\textbf{p\,1\,1,q\,1,r+}}}

Graph $G(p\,1\,1,q\,1,r+)$ (Fig. 54b) resolves into the graph
$G(r\,(q+1)\,1\,1\,p)$ and the block sum of the graphs
$G(p\,1\,1\,1\,r)$ and $G(\overline q)$. The general formula for the
Tutte polynomial of the graphs $G(p\,1\,1,q\,1,r+)$ is

$$T(G(p\,1\,1,q\,1,r+))=T(G((r\,(q+1)\,1\,1\,p))+T(G(p\,1\,1\,1\,r))T(G(\overline q)).$$

\subsection{Family \it{\textbf{(p\,q,r)\,(s,t)}}}

In order to obtain formula for the Tutte polynomial of the graph
$G((p\,q,r)\,(s,t))$ (Fig. 54c) we use the relations

$$T(G((p\,q,r)\,(s,t)))-T(G((p\,q,r)\,(s,(t-1))))=x^{t-1}T(G(p\,q\,s\,r)).$$

\noindent Since the Tutte polynomial of the graph
$G((p\,q,r)\,(s,0))$ is
$T(G((p\,q,r)\,(s,0)))=T(G(p\,(q+r)))T(G(s))$, the general formula
for the Tutte polynomial of the graphs $G((p\,q,r)\,(s,t))$ is

$$T(G((p\,q,r)\,(s,t)))={{x^t-1}\over
{x-1}}T(G(p\,q\,s\,r))+T(G(p\,(q+r)))T(G(s)).$$

\subsection{Family \it{\textbf{(p\,1\,1,q)\,(r,s)}}}

\begin{figure}[th]
\centerline{\psfig{file=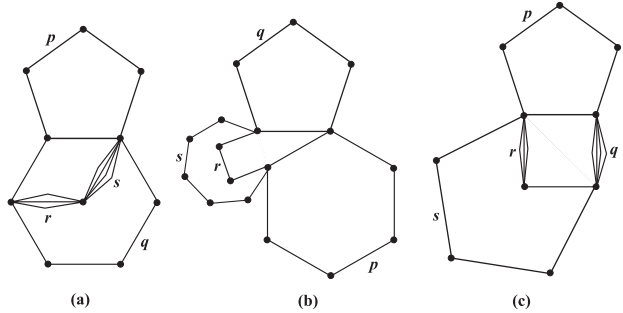,width=3.40in}} \vspace*{8pt}
\caption{Graphs (a) $G((p\,1\,1,q)\,(r,s))$; (b)
$G((p\,1,q\,1)\,(r,s))$; (c) $G((p\,1,q)\,(r\,1,s))$.\label{f1.55}}
\end{figure}

Graph $G((p\,1\,1,q)\,(r,s))$ (Fig. 55a) resolves into the graph
$G((p+1,q)\,(r,s))$ and the block sum of the graphs
$G((1,q)\,(r,s))$ and $G(p)$. The general formula for the Tutte
polynomial of the graphs $G((p\,1\,1,q)\,(r,s))$ is

$$T(G((p\,1\,1,q)\,(r,s)))=T(G(((p+1,q)\,(r,s)))+T(G((1,q)\,(r,s)))T(G(p)).$$

\subsection{Family \it{\textbf{(p\,1,q\,1)\,(r,s)}}}

In order to obtain formula for the Tutte polynomial of the graph
$G((p\,1,q\,1)\,(r,s))$ (Fig. 55b) we use the relations

$$T(G((p\,1,q\,1)\,(r,s)))-T(G((p\,1,q\,1)\,(r,(s-1))))=x^{s-1}T(G(p\,1\,r\,1\,q)).$$

\noindent Since the Tutte polynomial of the graph
$G((p\,1,q\,1)\,(r,0))$ is
$T(G((p\,1,q\,1)\,(r,0)))=T(G(p\,2\,q))T(G(r))$, the general formula
for the Tutte polynomial of the graphs $G((p\,1,q\,1)\,(r,s))$ is

$$T(G((p\,1,q\,1)\,(r,s))={{x^s-1}\over
{x-1}}T(G(p\,1\,r\,1\,q))+T(G(p\,2\,q))T(G(r)).$$

\subsection{Family \it{\textbf{(p\,1,q)\,(r\,1,s)}}}

Graph $G((p\,1,q)\,(r\,1,s))$ (Fig. 55c) resolves into the graph
$G(s\,r,p\,1,q)$ and the block sum of the graphs $G(p\,1\,s\,q)$ and
$G(\overline r)$. The general formula for the Tutte polynomial of
the graphs $G((p\,1,q)\,(r\,1,s))$ is

$$T(G((p\,1,q)\,(r\,1,s)))=T(G(s\,r,p\,1,q))+T(G(p\,1\,s\,q)T(G(\overline r)).$$

\subsection{Family \it{\textbf{(p,q,r)\,(s,t)}}}

\begin{figure}[th]
\centerline{\psfig{file=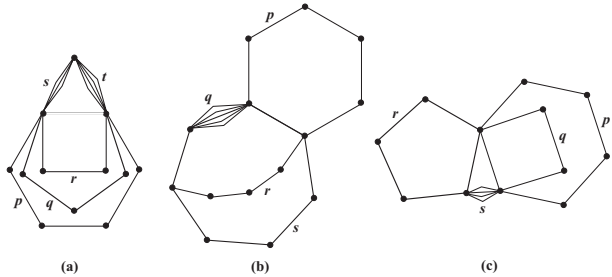,width=3.40in}} \vspace*{8pt}
\caption{Graphs (a) $G((p,q,r)\,(s,t))$; (b) $G((p\,1,q+)\,(r,s))$;
(c) $G((p,q+)\,(r\,1,s))$.\label{f1.56}}
\end{figure}

In order to obtain formula for the Tutte polynomial of the graph
$G((p,q,r)\,(s,t))$ (Fig. 56a) we use the relations

$$T(G((p,q,r)\,(s,t))-T(G(((p-1),q,r)\,(s,t)))=x^{p-1}T(G((q,r)\,(s,t))).$$

\noindent Since the Tutte polynomial of the graph
$G((0,q,r)\,(s,t))$ is
$T(G((0,q,r)\,(s,t)))=T(G(\overline{s+t}))T(G(r))T(G(q))$, the
general formula for the Tutte polynomial of the graphs
$G((p,q,r)\,(s,t))$ is

$$T(G((p,q,r)\,(s,t))={{x^p-1}\over
{x-1}}T(G((q,r)\,(s,t)))+T(G(\overline{s+t}))T(G(q))T(G(r)).$$

\subsection{Family \it{\textbf{(p\,1,q+)\,(r,s)}}}

Graph $G((p\,1,q+)\,(r,s))$ (Fig. 56b) resolves into the graph
$G((p\,1,q)\,(r,s))$ and the block sum of the graphs $G(r+s)$,
$G(p+1)$, and $G(\overline q)$. The general formula for the Tutte
polynomial of the graphs $G((p\,1,q+)\,(r,s))$ is

$$T(G((p\,1,q+)\,(r,s)))=T(G((p\,1,q)\,(r,s)))+T(G(r+s))T(G(p+1))T(G(\overline q)).$$

\subsection{Family \it{\textbf{(p,q+)\,(r\,1,s)}}}

Graph $G((p,q+)\,(r\,1,s))$ (Fig. 56c) resolves into the graph
$G((r\,1,s)\,(p,q))$ and the block sum of the graphs $G(r\,((s+1))$,
$G(p)$, and $G(q)$. The general formula for the Tutte polynomial of
the graphs $G((p,q+)\,(r\,1,s))$ is

$$T(G((p,q+)\,(r\,1,s)))=T(G((r\,1,s)\,(p,q)))+T(G(r\,((s+1)))T(G(p))T(G(q)).$$

\subsection{Family \it{\textbf{(p,q+)\,(r,s+)}}}

\begin{figure}[th]
\centerline{\psfig{file=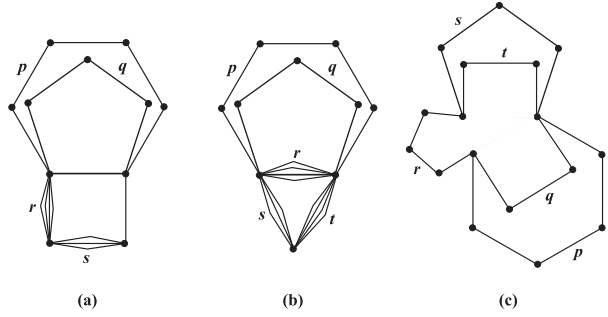,width=3.40in}} \vspace*{8pt}
\caption{Graphs (a) $G((p,q+)\,(r,s+))$; (b) $G((p,q+r)\,(s,t))$;
(c) $G((p,q)\,r\,(s,t))$.\label{f1.57}}
\end{figure}

Graph $G((p,q+)\,(r,s+))$ (Fig. 57a) resolves into the graph
$G((p,q+)\,(r,s))$ and the block sum of the graphs $G(p,q,1)$,
$G(\overline r)$, and $G(\overline s)$. The general formula for the
Tutte polynomial of the graphs $G((p,q+)\,(r,s+))$ is

$$T(G((p,q+)\,(r,s+))=T(G((p,q+)\,(r,s))+T(G(p,q,1)T(G(\overline r))T(G(\overline s)).$$

\subsection{Family \it{\textbf{(p,q+r)\,(s,t)}}}

In order to obtain formula for the Tutte polynomial of the graph
$G((p,q+r)\,(s,t))$ (Fig. 57b) we use the relations

$$T(G((p,q+r)\,(s,t))-T(G((p,q+r)\,(s,(t-1))))=y^{t-1}T(G(p\,(r+s)\,q)).$$

\noindent Since the Tutte polynomial of the graph
$G((p,q+r)\,(s,0))$ is
$T(G((p,q+r)\,(s,0)))=T(G(p\,r\,q))T(G(\overline s))$, the general
formula for the Tutte polynomial of the graphs $G((p,q+r)\,(s,t))$
is

$$T(G((p,q+r)\,(s,t))={{y^t-1}\over
{y-1}}T(G(p\,(r+s)\,q))+T(G(p\,r\,q))T(G(\overline s)).$$

\subsection{Family \it{\textbf{(p,q)\,r\,(s,t)}}}

In order to obtain formula for the Tutte polynomial of the graph
$G((p,q)\,r\,(s,t))$ (Fig. 57c) we use the relations

$$T(G((p,q)\,r\,(s,t))-T(G((p,q)\,r\,(s,(t-1))))=x^{t-1}T(G(p,q,(r+s))).$$

\noindent Since the Tutte polynomial of the graph
$G((p,q)\,r\,(s,0))$ is $T(G((p,q)\,r\,(s,0)))=T(G(p,q,r))T(G(s))$,
the general formula for the Tutte polynomial of the graphs
$G((p,q)\,r\,(s,t))$ is

$$T(G((p,q)\,r\,(s,t)))={{x^t-1}\over
{x-1}}T(G(p,q,(r+s)))+T(G(p,q,r)))T(G(s)).$$

\subsection{Family \it{\textbf{(p\,1,q)\,1\,(r,s)}}}

\begin{figure}[th]
\centerline{\psfig{file=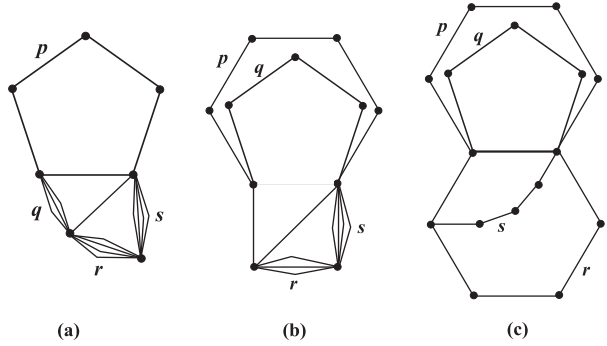,width=3.40in}} \vspace*{8pt}
\caption{Graphs (a) $G((p\,1,q)\,1\,(r,s))$; (b)
$G((p,q)\,1\,1\,(r,s))$; (c) $G((p,q+)\,1\,(r,s))$.\label{f1.58}}
\end{figure}

Graph $G((p\,1,q)\,1\,(r,s))$ (Fig. 58a) resolves into the graph
$G(p\,1,q,r,s)$ and the block sum of the graphs $G(p\,(q+1))$ and
$G(\overline {r+s})$. The general formula for the Tutte polynomial
of the graphs $G((p\,1,q)\,1\,(r,s))$ is

$$T(G((p\,1,q)\,1\,(r,s))=T(G(p\,1,q,r,s)+T(G(p\,(q+1)))T(G(\overline {r+s})).$$

\subsection{Family \it{\textbf{(p,q)\,1\,1\,(r,s)}}}

Graph $G((p,q)\,1\,1\,(r,s))$ (Fig. 58b) resolves into the graph
$G((p,q+)\,(r,s))$ and the block sum of the graphs $\overline
G(r\,1\,s)$ and $G(p+q)$. The general formula for the Tutte
polynomial of the graphs is $G((p,q)\,1\,1\,(r,s))$ is

$$T(G((p,q)\,1\,1\,(r,s))=T(G((p,q+)\,(r,s))+T(\overline
G(r\,1\,s))T(G(p+q)).$$

\subsection{Family \it{\textbf{(p,q+)\,1\,(r,s)}}}

Graph $G((p,q+)\,1\,(r,s))$ (Fig. 58c) resolves into the graph
$G(p,q,r,s,1)$ and the block sum of the graphs $G(p,q,1)$ and
$G(r+s)$. The general formula for the Tutte polynomial of the graphs
is $G((p,q+)\,1\,(r,s))$ is

$$T(G((p,q+)\,1\,(r,s)))=T(G(p,q,r,s,1))+T(G(p,q,1))T(G(r+s)).$$

\subsection{Family \it{\textbf{(p,q),r,(s,t)}}}

\begin{figure}[th]
\centerline{\psfig{file=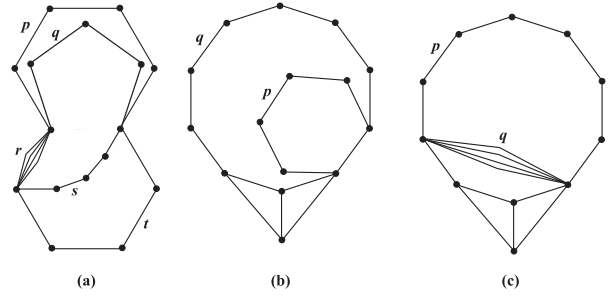,width=3.40in}} \vspace*{8pt}
\caption{Graphs (a) $G((p,q),r,(s,t))$; (b) $G(.p\,1\,q)$; (c)
$G(.p\,q\,1)$.\label{f1.59}}
\end{figure}

In order to obtain formula for the Tutte polynomial of the graph
$G((p,q),r,(s,t))$ (Fig. 59a) we use the relations

$$T(G((p,q),r,(s,t)))-T(G((p,q),(r-1),(s,t)))=y^{r-1}T(G(p,q,s,t)).$$

\noindent Since the Tutte polynomial of the graph $G((p,q),0,(s,t))$
is $T(G((p,q),0,(s,t)))=T(G(p+q))T(G(s+t))$, the general formula for
the Tutte polynomial of the graphs is $G((p,q),r,(s,t))$ is

$$T(G((p,q),r,(s,t)))={{y^r-1}\over
{y-1}}T(G(p,q,s,t))+T(G(p+q))T(G(s+t)).$$

\subsection{Family \it{\textbf{.p\,1\,q}}}

Graph $G(.p\,1\,q)$ (Fig. 59b) resolves into the graph $G(.(p+q))$
and the block sum of the graphs $G(.q)$ and $G(p)$. The general
formula for the Tutte polynomial of the graphs is $G(.p\,1\,q)$ is

$$T(G(.p\,1\,q))=T(G(.(p+q)))+T(G(.q))T(G(p)).$$

\subsection{Family \it{\textbf{.p\,q\,1}}}

In order to obtain formula for the Tutte polynomial of the graph
$G(.p\,q\,1)$ (Fig. 59c) we use the relations

$$T(G(.p\,q\,1))-T(G(.p\,(q-1)\,1))=y^{q-1}T(G(p))T(G(.1)),$$

\noindent where $T(G(.1))=2x+3x^2+x^3+2y+4xy+3y^2+y^3.$

\noindent Since the Tutte polynomial of the graph $G(.p\,0\,1)$ is
$T(G(.p\,0\,1))=T(G(.(p+1)))$, the general formula for the Tutte
polynomial of the graphs is $G(.p\,q\,1)$ is

$$T(G(.p\,q\,1))={{y^q-1}\over
{y-1}}T(G(p))T(G(.1))+T(G(.(p+1))).$$

\subsection{Family \it{\textbf{.p\,1\,1\,1}}}

\begin{figure}[th]
\centerline{\psfig{file=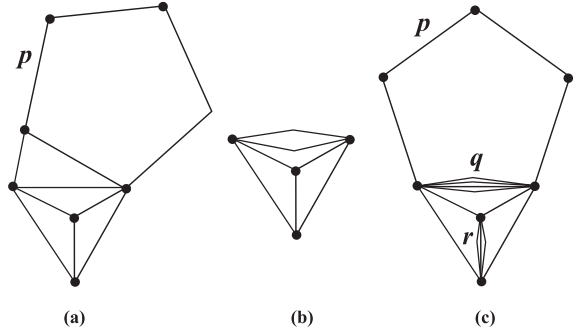,width=3.40in}} \vspace*{8pt}
\caption{Graphs (a) $G(.p\,1\,1\,1)$; (b) $G'$; (c)
$G(.p\,q:r)$.\label{f1.60}}
\end{figure}

Graph $G(.p\,1\,1\,1)$ (Fig. 60a) resolves into the graph
$G(.(p+1)\,1)$ and the block sum of the graphs $G(p)$ and $G'$ (Fig.
60b) with the Tutte polynomial
$$T(G')=2x+3x^2+x^3+2y+5xy+x^2y+4y^2+2xy^2+3y^3+y^4.$$

\noindent The general formula for the Tutte polynomial of the graphs
is $G(.p\,1\,1\,1)$ is

$$T(G(.p\,1\,1\,1))=T(G(.(p+1)\,1))+T(G')T(G(p)).$$

\subsection{Family \it{\textbf{.p\,q:r}}}

In order to obtain formula for the Tutte polynomial of the graph
$G(.p\,q:r)$ (Fig. 60c) we use the relations

$$T(G(.p\,q:r))-T(G(.p\,(q-1):r))=y^{q-1}T(\overline G(1\,(r-1),2,2))T(G(p)).$$

\noindent Since the Tutte polynomial of the graph $G(.p\,(0):r)$ is
$T(G(.p\,(0):r))=T(G(.p:r\,0))$, the general formula for the Tutte
polynomial of the graphs is $G(.p\,q:r)$ is

$$T(G(.p\,q:r))={{y^q-1}\over
{y-1}}T(G(p))T(\overline G(1\,(r-1),2,2))+T(G(.p:r\,0)).$$

\subsection{Family \it{\textbf{.p\,q.r}}}

\begin{figure}[th]
\centerline{\psfig{file=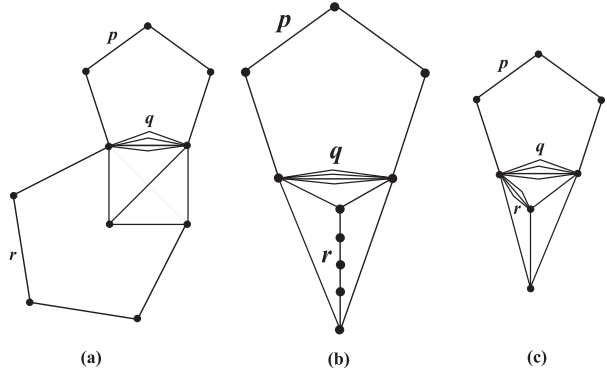,width=3.40in}} \vspace*{8pt}
\caption{Graphs (a) $G(.p\,q.r)$; (b) $G(.p\,q:r\,0)$; (c)
$G(.p\,q.r\,0)$.\label{f1.61}}
\end{figure}

Graph $G(.p\,q.r)$ (Fig. 61a) resolves into the graphs
$G(p\,q\,1\,1\,(r+1))$ and $G(p\,(q+2)\,r)$, and the block sum of
the graphs $G(p\,q\,r)$ and $G(2)$ with the Tutte polynomial
$T(G(2))=x+y.$

\noindent The general formula for the Tutte polynomial of the graphs
is $G(.p\,q.r)$ is

$$T(G(.p\,q.r))=T(G(p\,q\,1\,1\,(r+1)))+T(G(p\,(q+2)\,r))+T(G(p\,q\,r))T(G(2)).$$

\subsection{Family \it{\textbf{.p\,q:r\,0}}}

In order to obtain formula for the Tutte polynomial of the graph
$G(.p\,q:r\,0)$ (Fig. 61b) we use the relations

$$T(G(.p\,q:r\,0))-T(G(.(p-1)\,q:r\,0))=x^{p-1}T(G(.r:q\,0)).$$

\noindent Since the Tutte polynomial of the graph $G(.(0)\,q:r\,0)$
is $T(G(.(0)\,q:r\,0))=y^qT(G(r\,\overline 2\,\overline 2))$, the
general formula for the Tutte polynomial of the graphs
$G(.p\,q:r\,0)$ is

$$T(G(.p\,q:r\,0))={{x^p-1}\over
{x-1}}T(G(.r:q\,0))+y^qT(G(r\,\overline 2\,\overline 2)).$$

\subsection{Family \it{\textbf{.p\,q.r\,0}}}

Graph $G(.p\,q.r\,0)$ (Fig. 61c) resolves into the graphs
$G(2\,r\,1\,q\,p)$ and $G(((p,1)+q)\,(r,2))$. The general formula
for the Tutte polynomial of the graphs $G(.p\,q.r\,0)$ is

$$T(G(.p\,q.r\,0))=T(G(2\,r\,1\,q\,p))+T(G(((p,1)+q)\,(r,2))).$$

\subsection{Family \it{\textbf{.p\,1\,1:q}}}

\begin{figure}[th]
\centerline{\psfig{file=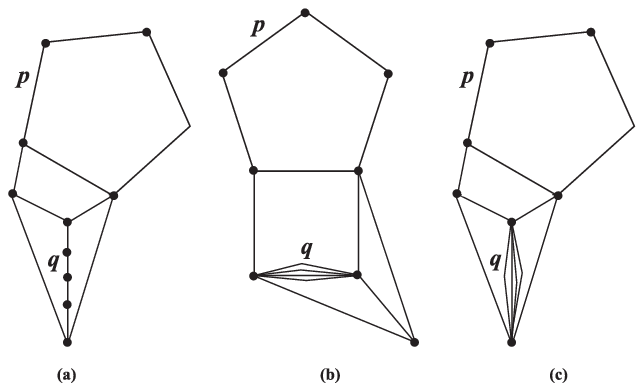,width=3.40in}} \vspace*{8pt}
\caption{Graphs (a) $G(.p\,1\,1:q)$; (b) $G(.p\,1\,1.q)$; (c)
$G(.p\,1\,1:q\,0)$.\label{f1.62}}
\end{figure}

Graph $G(.p\,1\,1:q)$ (Fig. 62a) resolves into the graph
$G(.p\,1:q\,0)$ and the block sum of the graphs $G(2,q,2)$ and
$G(p+1)$. The general formula for the Tutte polynomial of the graphs
$G(.p\,1\,1:q)$ is

$$T(G(.p\,1\,1:q))=T(G(.p\,1:q\,0))+T(2,q,2)T(G(p+1)).$$

\subsection{Family \it{\textbf{.p\,1\,1.q}}}

Graph $G(.p\,1\,1.q)$ (Fig. 62b) resolves into the graph
$G((p+1):1:q\,0)$ and the block sum of the graphs $G(.1:q\,0)$ and
$G(p)$. The general formula for the Tutte polynomial of the graphs
$G(.p\,1\,1.q)$ is

$$T(G(.p\,1\,1.q))=T(G((p+1):1:q\,0))+T(.1:q\,0)T(G(p)).$$

\subsection{Family \it{\textbf{.p\,1\,1:q\,0}}}

Graph $G(.p\,1\,1:q\,0)$ (Fig. 62c) resolves into graph $G(.p\,1:q)$
and the block sum of graphs $G(2\,q\,2)$ and $G(p+1)$. The general
formula for the Tutte polynomial of the graphs $G(.p\,1\,1:q\,0)$ is

$$T(G(.p\,1\,1:q\,0))=T(G(.p\,1:q)+T(2\,q\,2)T(G(p+1)).$$

\subsection{Family \it{\textbf{.p\,1\,1.q\,0}}}

\begin{figure}[th]
\centerline{\psfig{file=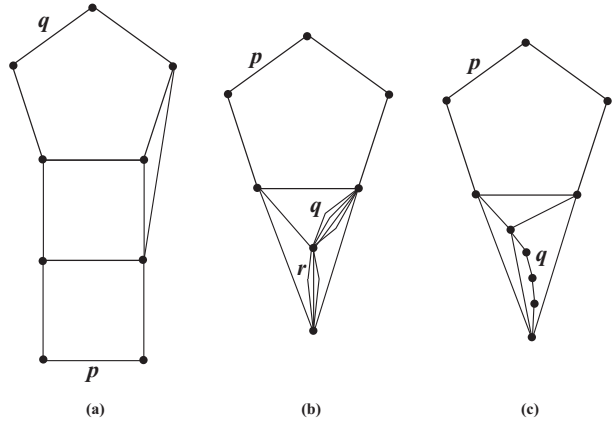,width=3.40in}} \vspace*{8pt}
\caption{Graphs (a) $G(.p\,1\,1.q\,0)$; (b) $G(.p\,1.q\,0.r)$; (c)
$G(.p\,1:q\,1)$.\label{f1.63}}
\end{figure}

Graph $G(.p\,1\,1.q\,0)$ (Fig. 63a) resolves into the graphs
$G(p\,1,q\,1,2)$, $G(p,(q+1),2,1)$, and the block sum of the graphs
$G((q+1)\,1\,2)$ and $G(p+1)$. The general formula for the Tutte
polynomial of the graphs $G(.p\,1\,1.q\,0)$ is

$$T(G(.p\,1\,1.q\,0))=T(G(p\,1,q\,1,2))+T(G(p,(q+1),2,1))+T((q+1)\,1\,2)T(G(p+1)).$$

\subsection{Family \it{\textbf{.p\,1.q\,0.r}}}

Graph $G(.p\,1.q\,0.r)$ (Fig. 63b) resolves into the graph
$G(.p.q.r\,0)$ and the block sum of the graphs $\overline
G((q+1),r,2)$ and $G(p)$. The general formula for the Tutte
polynomial of the graphs $G(.p\,1.q\,0.r)$ is

$$T(G(.p\,1.q\,0.r))=T(G(.p.q.r\,0))+T(\overline G((q+1),r,2))T(G(p)).$$

\subsection{Family \it{\textbf{.p\,1:q\,1}}}

Graph $G(.p\,1:q\,1)$ (Fig. 63c) resolves into the graph
$G(.q\,1:p\,0)$ and the block sum of the graphs $G(q\,1,2,2)$ and
$G(p)$. The general formula for the Tutte polynomial of the graphs
$G(.p\,1:q\,1)$ is

$$T(G(.p\,1:q\,1))=T(G(.q\,1:p\,0))+T(q\,1,2,2)T(G(p)).$$

\subsection{Family \it{\textbf{.p\,1.q\,1}}}

\begin{figure}[th]
\centerline{\psfig{file=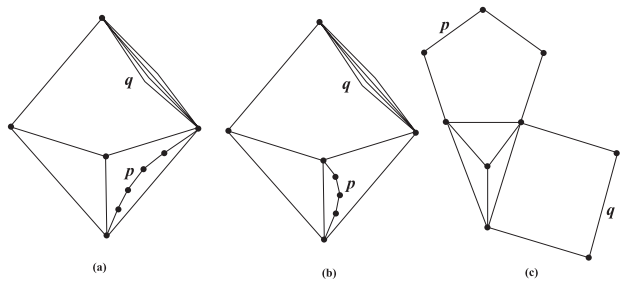,width=3.40in}} \vspace*{8pt}
\caption{Graphs (a) $G(.p\,1.q\,1)$; (b) $G(.p\,1:q\,1\,0)$; (c)
$G(.p\,1.q\,1\,0)$.\label{f1.64}}
\end{figure}

Graph $G(.p\,1.q\,1)$ (Fig. 64a) resolves into the graph
$G(.p\,1.q\,0)$ and the block sum of the graphs $G(p\,1\,1\,1\,2)$
and $G(\overline q)$. The general formula for the Tutte polynomial
of the graphs $G(.p\,1.q\,1)$ is

$$T(G(.p\,1.q\,1))=T(G(.p\,1.q\,0))+T(p\,1\,1\,1\,2)T(G(\overline q)).$$

\subsection{Family \it{\textbf{.p\,1:q\,1\,0}}}

Graph $G(.p\,1:q\,1\,0)$ (Fig. 64b) resolves into the graph
$G(.p\,1:q)$ and the block sum of the graphs $G(p,2,2,1)$ and
$G(\overline q)$. The general formula for the Tutte polynomial of
the graphs $G(.p\,1:q\,1\,0)$ is

$$T(G(.p\,1:q\,1\,0))=T(G(.p\,1:q))+T(p,2,2,1)T(G(\overline q)).$$

\subsection{Family \it{\textbf{.p\,1.q\,1\,0}}}

Graph $G(.p\,1.q\,1\,0)$ (Fig. 64c) resolves into the graph
$G(.q\,1.p)$ and the block sum of the graphs $G(q\,2\,1\,2)$ and
$G(p)$. The general formula for the Tutte polynomial of the graphs
$G(.p\,1.q\,1\,0)$ is

$$T(G(.p\,1.q\,1\,0))=T(G(.q\,1.p))+T(q\,2\,1\,2)T(G(p)).$$

\subsection{Family \it{\textbf{.p\,1\,0.q.r}}}

\begin{figure}[th]
\centerline{\psfig{file=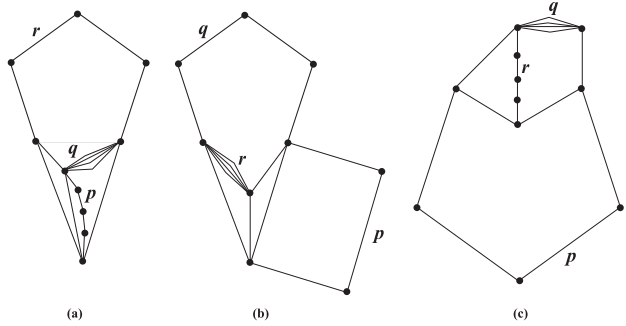,width=3.40in}} \vspace*{8pt}
\caption{Graphs (a) $G(.p\,1\,0.q.r)$; (b) $G(p\,1.q.r)$; (c)
$G(.p.q\,1\,0.r)$.\label{f1.65}}
\end{figure}

Graph $G(.p\,1\,0.q.r)$ (Fig. 65a) resolves into the graph
$G(.r.q.p)$ and the block sum of the graphs
$G(r\,(\overline{q+1})\,\overline 2)$ and $G(p)$. The general
formula for the Tutte polynomial of the graphs $G(.p\,1\,0.q.r)$ is

$$T(G(.p\,1\,0.q.r))=T(G(.r.q.p))+T(G(r\,(\overline{q+1})\,\overline 2))T(G(p)).$$

\subsection{Family \it{\textbf{.p\,1.q.r}}}

Graph $G(.p\,1.q.r)$ (Fig. 65b) resolves into the graph
$G(.p:r\,0.q\,0)$ and the block sum of the graphs $G(q\,1,r,2)$ and
$G(p)$. The general formula for the Tutte polynomial of the graphs
$G(.p\,1.q.r)$ is

$$T(G(.p\,1.q.r))=T(G(.p:r\,0.q\,0))+T(G(q\,1,r,2))T(G(p)).$$

\subsection{Family \it{\textbf{.p.q\,1\,0.r}}}

Graph $G(.p.q\,1\,0.r)$ (Fig. 65c) resolves into the graph
$G(.p.q.r)$ and the block sum of the graphs $G((p+1)\,1\,(r+1))$ and
$G(\overline q)$. The general formula for the Tutte polynomial of
the graphs $G(.p.q\,1\,0.r)$ is

$$T(G(.p.q\,1\,0.r))=T(G(.p.q.r))+T(G((p+1)\,1\,(r+1)))T(G(\overline q)).$$

\subsection{Family \it{\textbf{.p.q\,1.r}}}

\begin{figure}[th]
\centerline{\psfig{file=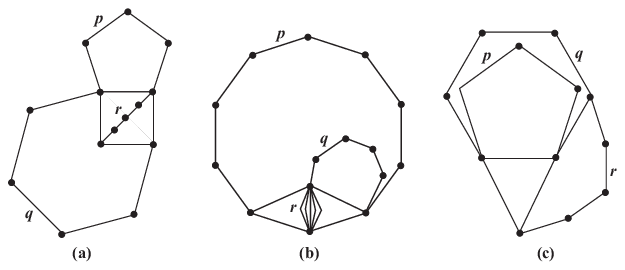,width=3.40in}} \vspace*{8pt}
\caption{Graphs (a) $G(.p.q\,1.r)$; (b) $G(.p.q\,1.r\,0)$; (c)
$G(.p\,1.q.r\,0)$.\label{f1.66}}
\end{figure}

Graph $G(.p.q\,1.r)$ (Fig. 66a) resolves into the graph
$G(.p.q\,0.r)$ and the block sum of the graphs $G(p\,1\,1\,1\,r)$
and $G(q)$. The general formula for the Tutte polynomial of the
graphs $G(.p.q\,1.r)$ is

$$T(G(.p.q\,1.r))=T(G(.p.q\,0.r))+T(G(p\,1\,1\,1\,r)T(G(q)).$$

\subsection{Family \it{\textbf{.p.q\,1.r\,0}}}

Graph $G(.p.q\,1.r\,0)$ (Fig. 66b) resolves into the graphs
$G((p+1)\,1\,q\,r)$ and $G((p,q)\,(1,(r+1)))$, and the block sum of
the graphs $G(p\,1\,1\,(r+1))$ and $G(q)$. The general formula for
the Tutte polynomial of the graphs $G(.p.q\,1.r\,0)$ is

$$T(G(.p.q\,1.r\,0))=T(G((p+1)\,1\,q\,r))+T(G((p,q)\,(1,(r+1))))+T(G(p\,1\,1\,(r+1)))T(G(q)).$$

\subsection{Family \it{\textbf{.p\,1.q.r\,0}}}

Graph $G(.p\,1.q.r\,0)$ (Fig. 66c) resolves into the graphs
$G(p\,1\,1\,1\,(q+r))$ and $G(r\,1\,1,q,p+)$. The general formula
for the Tutte polynomial of the graphs $G(.p\,1.q.r\,0)$ is

$$T(G(.p\,1.q.r\,0))=T(G(p\,1\,1\,1\,(q+r)))+T(G(r\,1\,1,q,p+)).$$

\subsection{Family \it{\textbf{p\,1\,0:q\,0:r\,0}}}

\begin{figure}[th]
\centerline{\psfig{file=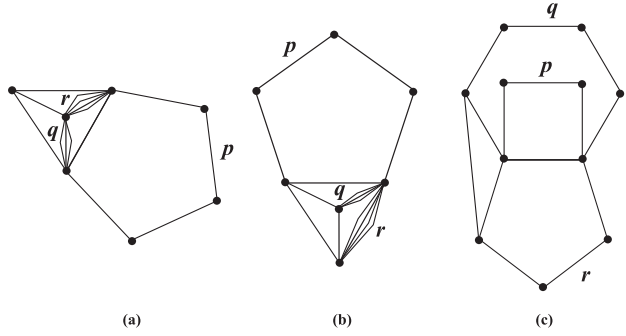,width=3.40in}} \vspace*{8pt}
\caption{Graphs (a) $G(p\,1\,0:q\,0:r\,0)$; (b) $G(p\,1:q:r)$; (c)
$G(p\,1:q\,0:r\,0)$.\label{f1.67}}
\end{figure}

Graph $G(p\,1\,0:q\,0:r\,0)$ (Fig. 67a) resolves into the graphs
$G(p\,1,(q+1),(r+1))$ and $G((p,2+)\,(q,r))$. The general formula
for the Tutte polynomial of the graphs $G(p\,1\,0:q\,0:r\,0)$ is

$$T(G(p\,1\,0:q\,0:r\,0))=T(G(p\,1,(q+1),(r+1)))+T(((p,2+)\,(q,r))).$$

\subsection{Family \it{\textbf{p\,1:q:r}}}

Graph $G(p\,1:q:r)$ (Fig. 67b) resolves into the graphs
$G(p\,q\,1\,r)$ and $G(p\,(\overline {q+r})\,\overline 2)$, the
block sum of the graphs $G((p+1)\,r)$ and $G(\overline q)$, and the
block sum of the graphs $\overline G((q+1)\,1\,(r+1))$ and $G(p)$.
The general formula for the Tutte polynomial of the graphs
$G(p\,1:q:r)$ is

$$T(G(p\,1:q:r))=T(G(p\,q\,1\,r))+T(G(p\,(\overline {q+r})\,\overline 2))+T(G((p+1)\,r))T(G(\overline q))+$$

$$T(\overline G((q+1)\,1\,(r+1)))T(G(p)).$$

\subsection{Family \it{\textbf{p\,1:q\,0:r\,0}}}

Graph $G(p\,1:q\,0:r\,0)$ (Fig. 67c) resolves into the graphs
$G(p\,0:q\,0:r\,0)$ and the block sum of the graphs
$G(q\,1\,1\,1\,r)$ and $G(p)$. The general formula for the Tutte
polynomial of the graphs $G(p\,1:q\,0:r\,0)$ is

$$T(G(p\,1:q\,0:r\,0))=T(G(p\,0:q\,0:r\,0))+T(G(q\,1\,1\,1\,r))T(G(p)).$$

\subsection{Family \it{\textbf{p:q:r\,1\,0}}}

\begin{figure}[th]
\centerline{\psfig{file=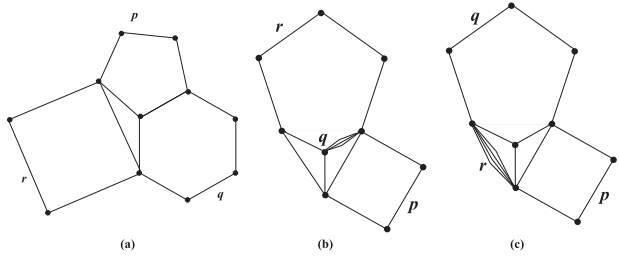,width=3.40in}} \vspace*{8pt}
\caption{Graphs (a) $G(p:q:r\,1\,0)$; (b) $G(p\,1:q:r\,0)$; (c)
$G(p\,1\,0:q:r\,0)$.\label{f1.68}}
\end{figure}

Graph $G(p:q:r\,1\,0)$ (Fig. 68a) resolves into the graphs
$G(p\,1,q\,1,r\,1)$ and $G((p+q),r,2,1)$. The general formula for
the Tutte polynomial of the graphs $G(p:q:r\,1\,0)$ is

$$T(G(p:q:r\,1\,0))=T(G(p\,1,q\,1,r\,1))+T(G((p+q),r,2,1)).$$

\subsection{Family \it{\textbf{p\,1:q:r\,0}}}

Graph $G(p\,1:q:r\,0)$ (Fig. 68b) resolves into the graphs
$G((r+1)\,q\,1\,1\,p)$ and $G((p,r+)\,(2,q))$. The general formula
for the Tutte polynomial of the graphs $G(p\,1:q:r\,0)$ is

$$T(G(p\,1:q:r\,0))=T(G((r+1)\,q\,1\,1\,p))+T(G((p,r+)\,(2,q))).$$

\subsection{Family \it{\textbf{p\,1\,0:q:r\,0}}}

Graph $G(p\,1\,0:q:r\,0)$ (Fig. 68c) resolves into the graph
$G(p:q:r\,0)$ and the block sum of the graphs $G(q\,r\,1\,2)$ and
$G(p)$. The general formula for the Tutte polynomial of the graphs
$G(p\,1\,0:q:r\,0)$ is

$$T(G(p\,1\,0:q:r\,0))=T(G(p:q:r\,0))+T(G(q\,r\,1\,2))T(G(p)).$$

\subsection{Family \it{\textbf{.p.q.r.s}}}

\begin{figure}[th]
\centerline{\psfig{file=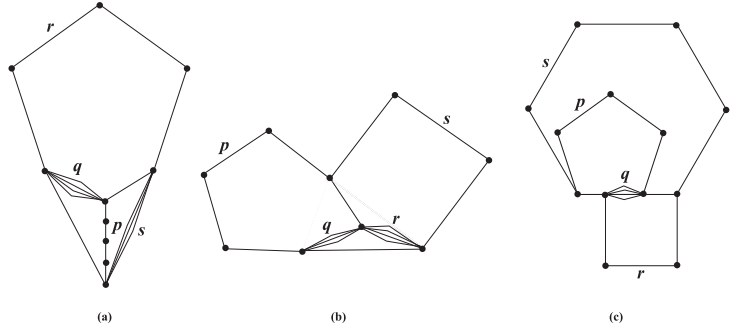,width=4.00in}} \vspace*{8pt}
\caption{Graphs (a) $G(.p.q.r.s)$; (b) $G(.p.q.r\,0.s\,0)$; (c)
$G(.p.q.r.s\,0)$.\label{f1.69}}
\end{figure}

In order to obtain formula for the Tutte polynomial of the graph
$G(.p.q.r.s)$ (Fig. 69a) we use the relations

$$T(G(.p.q.r.s))-T(G(.p.q.(r-1).s))=x^{r-1}(T(G(p\,q\,1\,s))+T(G((p+1)\,s))T(G(\overline q))).$$

\noindent Since the Tutte polynomial of the graph $G(.p.q.(0).s)$ is
$T(G(.p.q.(0).s))=T(G(p\,(\overline {q+1})\,(\overline {s+1})))$,
the general formula for the Tutte polynomial of the graphs
$G(.p.q.r.s)$ is

$$T(G(.p.q.r.s))={{x^r-1}\over
{x-1}}(T(G(p\,q\,1\,s))+T(G((p+1)\,s))T(G(\overline
q)))+T(G(p\,(\overline {q+1})\,(\overline {s+1}))).$$

\subsection{Family \it{\textbf{.p.q.r\,0.s\,0}}}

Graph $G(.p.q.r\,0.s\,0)$ (Fig. 69b) resolves into the graphs
$G(p\,q\,1\,r\,s)$ and $G((p+s)\,1,q,r)$. The general formula for
the Tutte polynomial of the graphs $G(.p.q.r\,0.s\,0)$ is

$$T(G(.p.q.r\,0.s\,0))=T(G(p\,q\,1\,r\,s))+T(G((p+s)\,1,q,r)).$$

\subsection{Family \it{\textbf{.p.q.r.s\,0}}}

Graph $G(.p.q.r.s\,0)$ (Fig. 69c) resolves into the graphs
$G((r+s)\,1\,p\,q)$ and $G(p,r,s+q)$, and the block sum of the
graphs $G(r\,q)$ and $G(p+s)$. The general formula for the Tutte
polynomial of the graphs $G(.p.q.r.s\,0)$ is

$$T(G(.p.q.r.s\,0))=T(G((r+s)\,1\,p\,q))+T(G(p,r,s+q))+T(G(r\,q))T(G(p+s)).$$

\subsection{Family \it{\textbf{.p.q\,0.r.s\,0}}}

\begin{figure}[th]
\centerline{\psfig{file=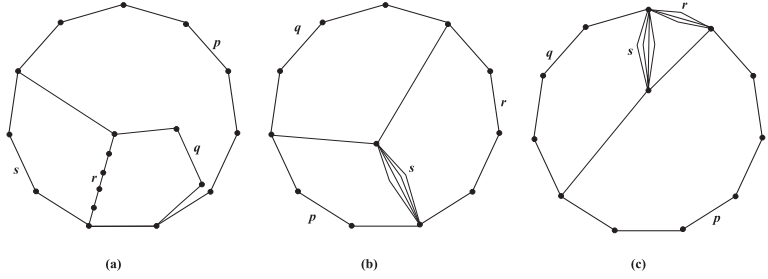,width=4.20in}} \vspace*{8pt}
\caption{Graphs (a) $G(.p.q\,0.r.s\,0)$; (b) $G(p\,0.q.r\,0.s\,0)$;
(c) $G(p\,0.q.r.s\,0)$.\label{f1.70}}
\end{figure}

Graph $G(.p.q\,0.r.s\,0)$ (Fig. 70a) resolves into the graphs
$G(p,q,r,s)$ and $G((p+q)\,1\,(r+s))$, and the block sum of the
graphs $G(p+s)$ and $G(q+r)$. The general formula for the Tutte
polynomial of the graphs $G(.p.q\,0.r.s\,0)$ is

$$T(G(.p.q\,0.r.s\,0))=T(G(p,q,r,s))+T(G((p+q)\,1\,(r+s)))+T(G(p+s))T(G(q+r)).$$

\subsection{Family \it{\textbf{p\,0.q.r\,0.s\,0}}}

Graph $G(p\,0.q.r\,0.s\,0)$ (Fig. 70b) resolves into the graphs
$G(r\,s\,q\,1\,p)$ and $G((p+r)\,s\,q)$, and the block sum of the
graphs $G(p+q+r)$ and $G(\overline s)$. The general formula for the
Tutte polynomial of the graphs $G(p\,0.q.r\,0.s\,0)$ is

$$T(G(p\,0.q.r\,0.s\,0))=T(G(r\,s\,q\,1\,p))+T(G((p+r)\,s\,q))+T(G(p+q+r))T(G(\overline s)).$$

\subsection{Family \it{\textbf{p\,0.q.r.s\,0}}}

Graph $G(p\,0.q.r.s\,0)$ (Fig. 70c) resolves into the graphs
$G((p+q)\,r\,1\,s)$ and $G(q\,s,p\,1,r)$. The general formula for
the Tutte polynomial of the graphs $G(p\,0.q.r.s\,0)$ is

$$T(G(p\,0.q.r.s\,0))=T(G((p+q)\,r\,1\,s))+T(G(q\,s,p\,1,r)).$$

\subsection{Family \it{\textbf{p.q\,0.r.s\,0}}}

\begin{figure}[th]
\centerline{\psfig{file=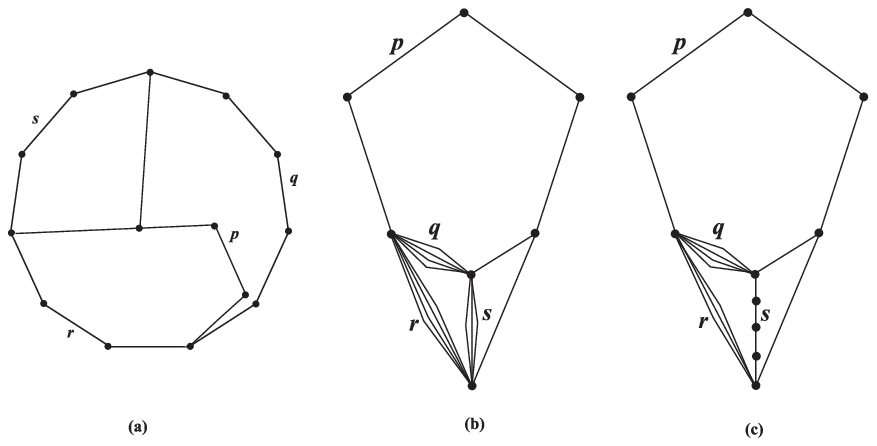,width=4.20in}} \vspace*{8pt}
\caption{Graphs (a) $G(p.q\,0.r.s\,0)$; (b) $G(p.q.r\,0.s)$; (c)
$G(p.q.r\,0.s\,0)$.\label{f1.71}}
\end{figure}

Graph $G(p.q\,0.r.s\,0)$ (Fig. 71a) resolves into the graphs
$G(p,q,(r+s))$ and $G(p,r,(q+s))$, the block sum of the graphs
$G(p,q,r)$ and $G(s)$, and the block sum of the graph $G(q+r+s)$ and
chain graph of the length $p$. The general formula for the Tutte
polynomial of the graphs $G(p.q\,0.r.s\,0)$ is

$$T(G(p.q\,0.r.s\,0))=T(G(p,q,(r+s)))+T(G(p,r,(q+s))) +T(G(p,q,r))T(G(s))+$$
$$x^pT(G(q+r+s)).$$

\subsection{Family \it{\textbf{p.q.r\,0.s}}}

Graph $G(p.q.r\,0.s)$ (Fig. 71b) resolves into the graphs
$G((p+1)\,q,r,s)$ and $G(p\,r,q,(s+1))$. The general formula for the
Tutte polynomial of the graphs $G(p.q.r\,0.s)$ is

$$T(G(p.q.r\,0.s))=T(G((p+1)\,q,r,s))+T(G(p\,r,q,(s+1))).$$

\subsection{Family \it{\textbf{p.q.r\,0.s\,0}}}

Graph $G(p.q.r\,0.s\,0)$ (Fig. 71c) resolves into the graphs
$G((p+1)\,q\,s\,r)$ and $G(p\,r,s\,1,q)$. The general formula for
the Tutte polynomial of the graphs $G(p.q.r\,0.s\,0)$ is

$$T(G(p.q.r\,0.s\,0))=T(G((p+1)\,q\,s\,r))+T(G(p\,r,s\,1,q)).$$

\subsection{Family \it{\textbf{p.q.r.s}}}

\begin{figure}[th]
\centerline{\psfig{file=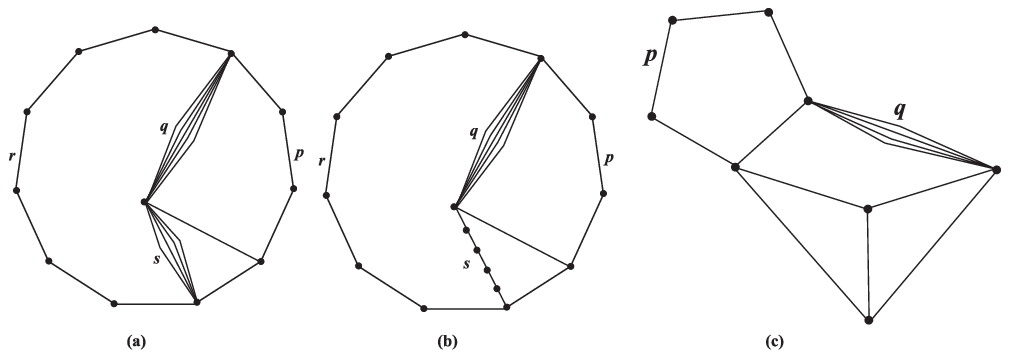,width=4.60in}} \vspace*{8pt}
\caption{Graphs (a) $G(p.q.r.s)$; (b) $G(p.q.r.s\,0)$; (c)
$G(.(p\,1,q))$.\label{f1.72}}
\end{figure}

Graph $G(p.q.r.s)$ (Fig. 72a) resolves into the graphs
$G((p+1)\,q\,r\,s)$ and $G((p,r)\,(q,(s+1)))$. The general formula
for the Tutte polynomial of the graphs $G(p.q.r.s)$ is

$$T(G(p.q.r.s))=T(G((p+1)\,q\,r\,s))+T(G((p,r)\,(q,(s+1)))).$$

\subsection{Family \it{\textbf{p.q.r.s\,0}}}

Graph $G(p.q.r.s\,0)$ (Fig. 72b) resolves into the graphs
$G'(s\,\overline q\,(p+1)\,r)$ and $G(p\,q\,r\,1\,s)$. The general
formula for the Tutte polynomial of the graphs $G(p.q.r.s\,0)$ is

$$T(G(p.q.r.s\,0))=T(G'(s\,\overline q\,(p+1)\,r))+T(G(p\,q\,r\,1\,s)).$$

\subsection{Family \it{\textbf{.(p\,1,q)}}}

Graph $G(.(p\,1,q))$ (Fig. 72c) resolves into the graphs
$G((p\,1,q)\,(2,2))$ and $G(p\,1,q,2,2)$. The general formula for
the Tutte polynomial of the graphs $G(.(p\,1,q))$ is

$$T(G(.(p\,1,q)))=T(G((p\,1,q)\,(2,2)))+T(G(p\,1,q,2,2)).$$

\subsection{Family \it{\textbf{.(p,q)\,1}}}

\begin{figure}[th]
\centerline{\psfig{file=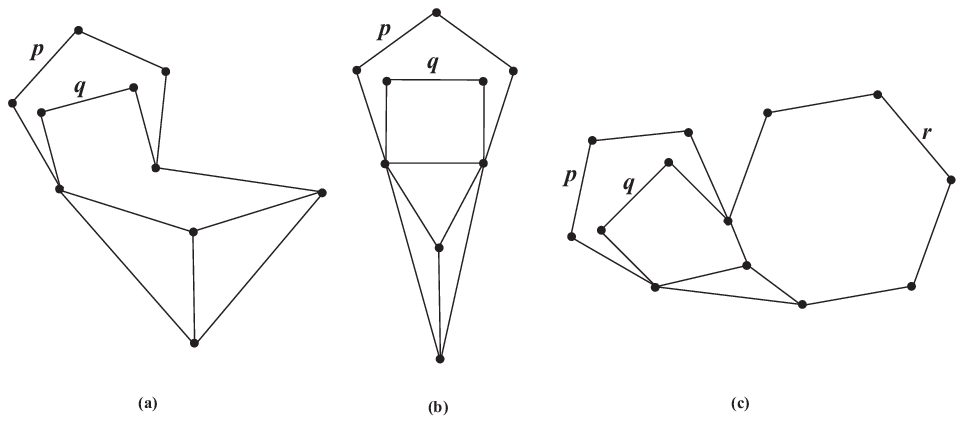,width=4.60in}} \vspace*{8pt}
\caption{Graphs (a) $G(.(p,q)\,1)$; (b) $G(.(p,q+))$; (c)
$G(.(p,q).r)$.\label{f1.73}}
\end{figure}

Graph $G(.(p,q)\,1)$ (Fig. 73a) resolves into the graph $G(.(p,q))$
and the block sum of the graphs $G(p+q)$ and $G'$, where $G'$ is the
graph which consists of two triangles with a common edge, with the
Tutte polynomial $T(G')=x+2x^2+x^3+y+2xy+y^2.$ The general formula
for the Tutte polynomial of the graphs $G(.(p,q)\,1)$ is

$$T(G(.(p,q)\,1))=T(G(.(p,q)))+T(G(p+q))(x+2x^2+x^3+y+2xy+y^2).$$

\subsection{Family \it{\textbf{.(p,q+)}}}

Graph $G(.(p,q+))$ (Fig. 73b) resolves into graph $G(.(p,q))$ and
the block sum of the graphs $G(p)$, $G(q)$, and $\overline
G(2\,1\,2)$ with the Tutte polynomial $T(\overline
G(2\,1\,2))=x+x^2+y+2xy+2y^2+y^3.$ The general formula for the Tutte
polynomial of the graphs $G(.(p,q+))$ is

$$T(G(.(p,q+)))=T(G(.(p,q)))+T(G(p))T(G(q))(x+x^2+y+2xy+2y^2+y^3).$$

\subsection{Family \it{\textbf{.(p,q).r}}}

Graph $G(.(p,q).r)$ (Fig. 73c) resolves into the graphs $G(p,q,r,2)$
and $G'(r\,\overline 2\,p\,q)$, the block sum of the graphs $G(p+q)$
and $G(r+2)$, and the block sum of the graphs $G(r\,2)$, $G(p)$, and
$G(q)$. The general formula for the Tutte polynomial of the graphs
$G(.(p,q).r)$ is

$$T(G(.(p,q).r))=T(G(p,q,r,2))+T(G'(r\,\overline 2\,p\,q))+T(G(p+q))T(G(r+2))+$$
$$T(G(r\,2))T(G(p))T(G(q)).$$

\subsection{Family \it{\textbf{.(p,q).r\,0}}}

\begin{figure}[th]
\centerline{\psfig{file=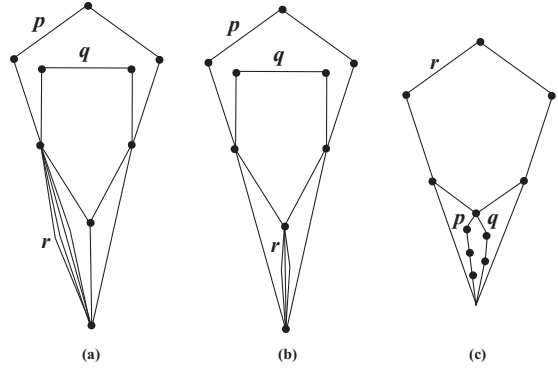,width=3.00in}} \vspace*{8pt}
\caption{Graphs (a) $G(.(p,q).r\,0)$; (b) $G(.(p,q):r)$; (c)
$G(.(p,q):r\,0)$.\label{f1.74}}
\end{figure}

Graph $G(.(p,q).r\,0)$ (Fig. 74a) resolves into the graphs
$G((p,q)\,((r+1),2))$ and $G(p,q,2+r)$, and the block sum of the
graphs $G(p,q,2)$ and $G(\overline r)$. The general formula for the
Tutte polynomial of the graphs $G(.(p,q).r\,0)$ is

$$T(G(.(p,q).r\,0))=T(G((p,q)\,((r+1),2)))+T(G(p,q,2+r))+T(G(p,q,2))T(G(\overline r)).$$

\subsection{Family \it{\textbf{.(p,q):r}}}

In order to obtain formula for the Tutte polynomial of the graph
$G(.(p,q):r)$ (Fig. 74b) we use the relations

$$T(G(.(p,q):r))-T(G(.(p,(q-1)):r))=x^{q-1}T(G(.p:r\,0)).$$

\noindent Since the Tutte polynomial of the graph $G(.(p,0):r)$ is
$T(G(.(p,0):r))=T(\overline G(r,2,2))$, the general formula for the
Tutte polynomial of the graphs $G(.(p,q):r)$ is

$$T(G(.(p,q):r))={{x^q-1}\over {x-1}}T(G(.p:r\,0))+T(\overline G(r,2,2))T(G(p)).$$

\subsection{Family \it{\textbf{.(p,q):r\,0}}}

In order to obtain formula for the Tutte polynomial of the graph
$G(.(p,q):r\,0)$ (Fig. 74c) we use the relations

$$T(G(.(p,q):r\,0))-T(G(.(p,q):(r-1)\,0))=x^{r-1}T(G(p,q,2,2)).$$

\noindent Since the Tutte polynomial of the graph $G(.(p,q):(0)\,0)$
is $T(G(.(p,q):(0)\,0))=T(G((p,q)\,(2,2)))$, the general formula for
the Tutte polynomial of the graphs $G(.(p,q):r\,0)$ is

$$T(G(.(p,q):r\,0))={{x^r-1}\over {x-1}}T(G(p,q,2,2))+T(G((p,q)\,(2,2))).$$

\subsection{Family \it{\textbf{8$^*$p\,0::q\,0}}}

\begin{figure}[th]
\centerline{\psfig{file=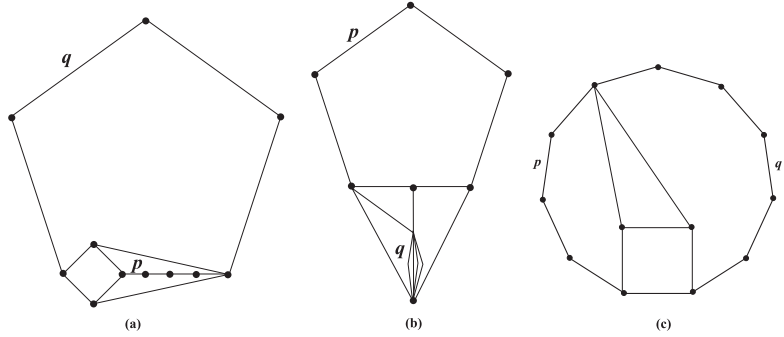,width=4.20in}} \vspace*{8pt}
\caption{Graphs (a) $G(8^*p\,0::q\,0)$; (b) $G(8^*p\,0:.q\,0)$; (c)
$G(8^*p\,0:q\,0)$.\label{f1.75}}
\end{figure}

In order to obtain formula for the Tutte polynomial of the graph
$G(8^*p\,0::q\,0)$ (Fig. 75a) we use the relations

$$T(G(8^*p\,0::q\,0))-T(G(8^*p\,0::(q-1)\,0))=x^{q-1}T(G(.p:2)).$$

\noindent Since the graph $G(8^*p\,0::(0)\,0)$ resolves into the
graph $G(p\,2\,1\,2)$ and the block sum of the graphs $G((p+1)\,2)$
and $G(2)$, the general formula for the Tutte polynomial of the
graphs $G(8^*p\,0::q\,0)$ is

$$T(G(8^*p\,0::q\,0))={{x^q-1}\over {x-1}}T(G(.p:2))+T(G(p\,2\,1\,2))+T(G((p+1)\,2))(x+y).$$

\subsection{Family \it{\textbf{8$^*$p\,0:.q\,0}}}

Graph $G(8^*p\,0:.q\,0)$ (Fig. 75b) resolves into the graphs
$G(.p.2.q\,0)$ and $G((p,2+)\,(q,1))$, and the block sum of the
graphs $G(2\,q)$ and $G(p+2)$. The general formula for the Tutte
polynomial of the graphs $G(8^*p\,0:.q\,0)$ is

$$T(G(8^*p\,0:.q\,0))=T(G(.p.2.q\,0))+T(G((p,2+)\,(q,1)))+T(G(2\,q))T(G(p+2)).$$

\subsection{Family \it{\textbf{8$^*$p\,0:q\,0}}}

Graph $G(8^*p\,0:q\,0)$ (Fig. 75c) resolves into the graphs
$G(.(p,q))$ and $G((p+1)\,1\,1\,1\,(q+1))$. The general formula for
the Tutte polynomial of the graphs $G(8^*p\,0:q\,0)$ is

$$T(G(8^*p\,0:q\,0))=T(G(.(p,q)))+T(G((p+1)\,1\,1\,1\,(q+1))).$$

\subsection{Family \it{\textbf{8$^*$p\,0.q\,0}}}

\begin{figure}[th]
\centerline{\psfig{file=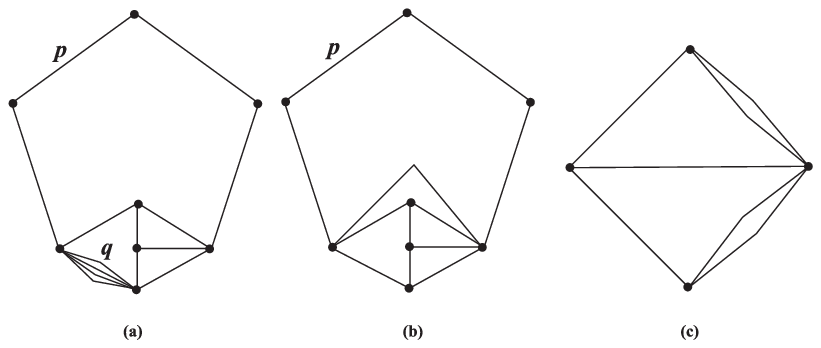,width=4.20in}} \vspace*{8pt}
\caption{Graphs (a) $G(8^*p\,0.q\,0)$; (b) $G(8^*p\,1)$; (c)
$G'$.\label{f1.76}}
\end{figure}

Graph $G(8^*p\,0.q\,0)$ (Fig. 76a) resolves into the graphs
$G(.1.q.p.2)$ and $G((2\,1,q)\,(p,2))$. The general formula for the
Tutte polynomial of the graphs $G(8^*p\,0.q\,0)$ is

$$T(G(8^*p\,0.q\,0))=T(G(.1.q.p.2))+T(G((2\,1,q)\,(p,2))).$$

\subsection{Family \it{\textbf{8$^*$p\,1}}}

Graph $G(8^*p\,1)$ (Fig. 76b) resolves into the graph $G(8^*p\,0)$
and the block sum of the graphs $G(p)$ and $G'$ (Fig. 76c) with the
Tutte polynomial $T(G')=x+2x^2+x^3+y+4xy+2x^2y+3y^2+3xy^2+3y^3+y^4.$
The general formula for the Tutte polynomial of the graphs
$G(8^*p\,1)$ is

$$T(G(8^*p\,1))=T(G(8^*p\,0))+T(G(p))(x+2x^2+x^3+y+4xy+2x^2y+3y^2+3xy^2+3y^3+y^4).$$

\subsection{Family \it{\textbf{8$^*$p\,1\,0}}}

\begin{figure}[th]
\centerline{\psfig{file=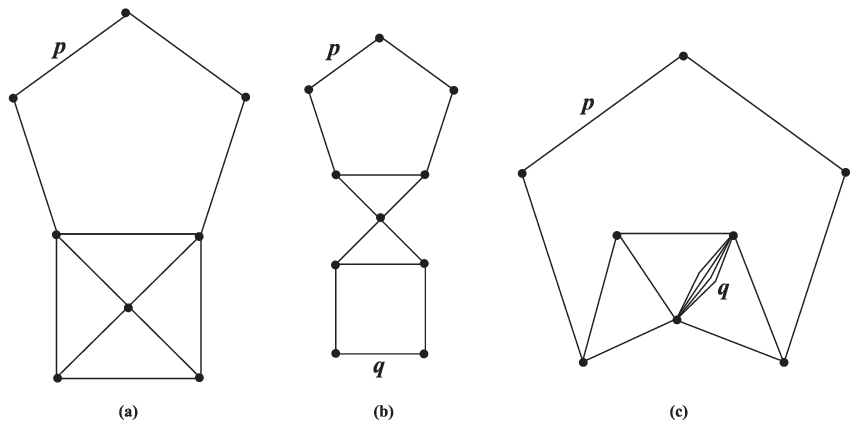,width=4.00in}} \vspace*{8pt}
\caption{Graphs (a) $G(8^*p\,1\,0)$; (b) $G(8^*p::q)$; (c)
$G(8^*p:.q)$.\label{f1.77}}
\end{figure}

Graph $G(8^*p\,1\,0)$ (Fig. 77a) resolves into the graph $G(8^*p)$
and the block sum of the graphs $G(p)$ and $G(.1:2\,0)$ with the
Tutte polynomial
$T(G(.1:2\,0))=2x+3x^2+x^3+2y+5xy+x^2y+4y^2+2xy^2+3y^3+y^4).$ The
general formula for the Tutte polynomial of the graphs
$G(8^*p\,1\,0)$ is

$$T(G(8^*p\,1\,0))=T(G(8^*p))+T(G(p))(2x+3x^2+x^3+2y+5xy+x^2y+4y^2+2xy^2+3y^3+y^4).$$

\subsection{Family \it{\textbf{8$^*$p::q}}}

Graph $G(8^*p::q)$ (Fig. 77b) resolves into the graphs
$G((p+1)\,1\,1\,1\,(q+1))$ and $G(1\,0.p.q\,0.2\,0)$. The general
formula for the Tutte polynomial of the graphs $G(8^*p::q)$ is

$$T(G(8^*p::q))=T(G((p+1)\,1\,1\,1\,(q+1)))+T(G(1\,0.p.q\,0.2\,0)).$$

\subsection{Family \it{\textbf{8$^*$p:.q}}}

Graph $G(8^*p:.q)$ (Fig. 77c) resolves into the graphs
$G((p+1)\,1\,1\,1\,1\,q)$ and $G(1\,0.p.1\,0.(q+1)\,0)$. The general
formula for the Tutte polynomial of the graphs $G(8^*p:.q)$ is

$$T(G(8^*p:.q))=T(G((p+1)\,1\,1\,1\,1\,q))+T(G(1\,0.p.1\,0.(q+1)\,0)).$$

\subsection{Family \it{\textbf{8$^*$p:q}}}

\begin{figure}[th]
\centerline{\psfig{file=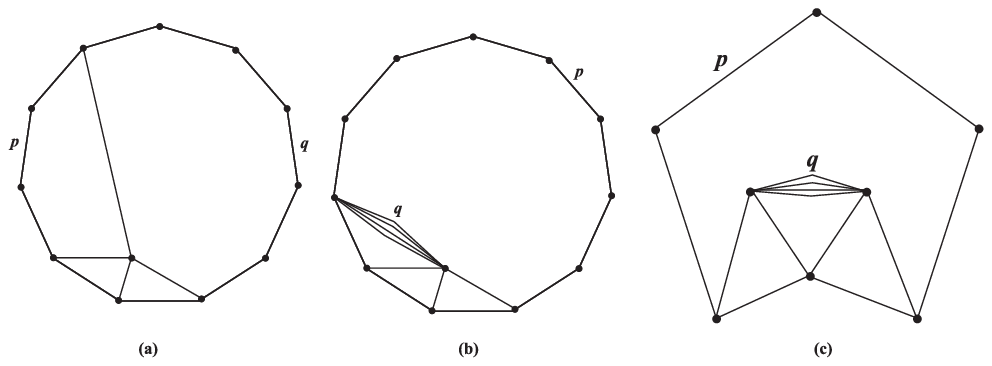,width=4.20in}} \vspace*{8pt}
\caption{Graphs (a) $G(8^*p:q)$; (b) $G(8^*p.q)$; (c)
$G(8^*p::q\,0)$.\label{f1.78}}
\end{figure}

Graph $G(8^*p:q)$ (Fig. 78a) resolves into the graphs $G(.(p+q))$
and $G(p\,1\,2\,1\,q)$, and the block sum of the graphs $G(p\,2)$
and $G(q\,2)$. The general formula for the Tutte polynomial of the
graphs $G(8^*p:q)$ is

$$T(G(8^*p:q))=T(G(.(p+q)))+T(G(p\,1\,2\,1\,q))+T(G(p\,2))T(G(q\,2)).$$

\subsection{Family \it{\textbf{8$^*$p.q}}}

Graph $G(8^*p.q)$ (Fig. 78b) resolves into the graphs $G(.p.q.1\,0)$
and $G(p\,(\overline {q+1})\,\overline 3)$, the block sum of the
graphs $G((p+1)\,(q+1))$ and $G(2)$, and the block sum of the graphs
$G(2\,1,q,1+(p-1))$ and $G(1)$. The general formula for the Tutte
polynomial of the graphs $G(8^*p.q)$ is

$$T(G(8^*p.q))=T(G(.p.q.1\,0))+T(G(p\,(\overline {q+1})\,\overline 3))+T(G((p+1)\,(q+1)))(x+y)+$$
$$xT(G(2\,1,q,1+(p-1))).$$

\subsection{Family \it{\textbf{8$^*$p::q\,0}}}

Graph $G(8^*p::q\,0)$ (Fig. 78c) resolves into the graphs
$G(.(p+1).q)$ and $G(p\,1\,1\,1,q,2)$. The general formula for the
Tutte polynomial of the graphs $G(8^*p::q\,0)$ is

$$T(G(8^*p::q\,0))=T(G(.(p+1).q))+T(G(p\,1\,1\,1,q,2)).$$

\subsection{Family \it{\textbf{8$^*$p:.q\,0}}}

\begin{figure}[th]
\centerline{\psfig{file=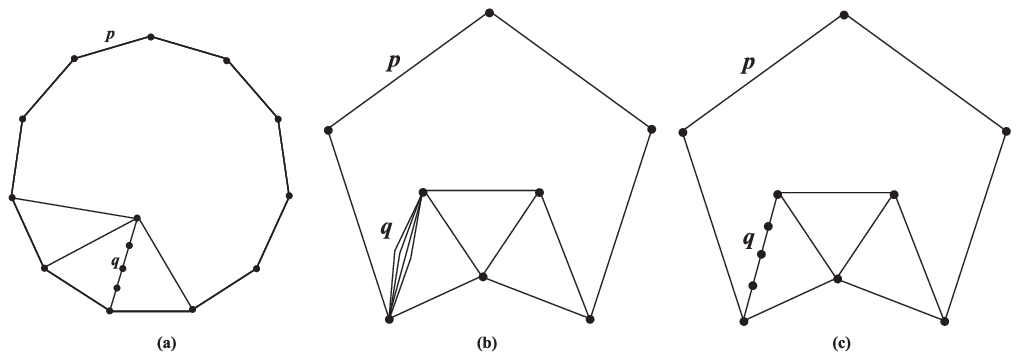,width=4.20in}} \vspace*{8pt}
\caption{Graphs (a) $G(8^*p:.q\,0)$; (b) $G(8^*p:q\,0)$; (c)
$G(8^*p.q\,0)$.\label{f1.79}}
\end{figure}

Graph $G(8^*p:.q\,0)$ (Fig. 79a) resolves into the graphs
$G(.q\,1:p\,0)$ and $G((q+1)\,1\,p\,1\,2)$. The general formula for
the Tutte polynomial of the graphs $G(8^*p:.q\,0)$ is

$$T(G(8^*p:.q\,0))=T(G(.q\,1:p\,0))+T(G((q+1)\,1\,p\,1\,2)).$$

\subsection{Family \it{\textbf{8$^*$p:q\,0}}}

Graph $G(8^*p:q\,0)$ (Fig. 79b) resolves into the graphs
$G(p\,1\,1\,1\,1\,(q+1))$, $G(p\,\overline q\,\overline 2\,\overline
2)$ and $G'(p\,\overline q\,2\,2)$. The general formula for the
Tutte polynomial of the graphs $G(8^*p:q\,0)$ is

$$T(G(8^*p:q\,0))=T(G(p\,1\,1\,1\,1\,(q+1)))+T(G(p\,\overline q\,\overline 2\,\overline
2))+T(G'(p\,\overline q\,2\,2)).$$

\subsection{Family \it{\textbf{8$^*$p.q\,0}}}

In order to obtain formula for the Tutte polynomial of the graph
$G(8^*p.q\,0)$ (Fig. 79c) we use the relations

$$T(G(8^*p.q\,0))-T(G(8^*p.(q-1)\,0))=x^{q-1}T(G(.(p+1))).$$

\noindent Since the Tutte polynomial of the graph $G(8^*p.(0)\,0)$
is $T(G(8^*p.(0)\,0))=T(G(p\,1\,1\,1\,1\,2))$, the general formula
for the Tutte polynomial of the graphs $G(8^*p.q\,0)$ is

$$T(G(8^*p.q\,0))={{x^q-1}\over {x-1}}T(G(.(p+1)))+T(G(p\,1\,1\,1\,1\,2)).$$

\subsection{Family \it{\textbf{9$^*$.p}}}

\begin{figure}[th]
\centerline{\psfig{file=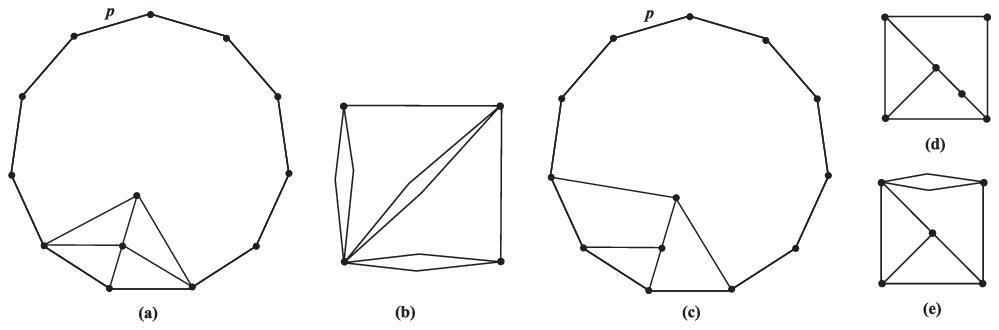,width=4.20in}} \vspace*{8pt}
\caption{Graphs (a) $G(9^*.p)$; (b) $G'$; (c) $G(9^*p)$; (d) $G''$;
(e) $G'''$.\label{f1.80}}
\end{figure}

In order to obtain formula for the Tutte polynomial of the graph
$G(9^*.p)$ (Fig.80a) we use the relations

$$T(G(9^*.p))-T(G(9^*.(p-1)))=x^{p-1}T(Wh(5)).$$

\noindent where $Wh(5)$ is the wheel graph with the Tutte polynomial
$$T(Wh(5))=3x+6x^2+4x^3+x^4+3y+9xy+4x^2y+6y^2+4xy^2+4y^3+y^4.$$

\noindent Since $G'=G(9^*.(0))$ is the graph $G'$ (Fig. 80b) with
the Tutte polynomial
$$T(G')=x+2x^2+x^3+y+4xy+3x^2y+3y^2+5xy^2+4y^3+2xy^3+3y^4+y^5,$$ the general formula
for the Tutte polynomial of the graphs $G(9^*.p)$ is

$$T(G(9^*.p))={{x^p-1}\over {x-1}}T(Wh(5))+T(G').$$

\subsection{Family \it{\textbf{9$^*$p}}}

In order to obtain formula for the Tutte polynomial of the graph
$G(9^*p)$ (Fig. 80c) we use the relations

$$T(G(9^*p))-T(G(9^*(p-1)))=x^{p-1}T(G'),$$

\noindent where $G''$ is the graph (Fig. 80d)  with the Tutte
polynomial
$$T(G')=2x+5x^2+5x^3+3x^4+x^5+2y+6xy+4x^2y+x^3y+3y^2+2xy^2+y^3.$$

\noindent Since $G'''=G(9^*(0))$ is the graph (Fig. 80e) with the
Tutte polynomial
$$T(G''')=2x+4x^2+3x^3+x^4+2y+7xy+5x^2y+x^3y+5y^2+5xy^2+4y^3+y^4,$$ the general formula
for the Tutte polynomial of the graphs $G(9^*p)$ is

$$T(G(9^*p))={{x^p-1}\over {x-1}}T(G'')+T(G''').$$

\subsection{Family \it{\textbf{9$^*$.p\,0}}}

\begin{figure}[th]
\centerline{\psfig{file=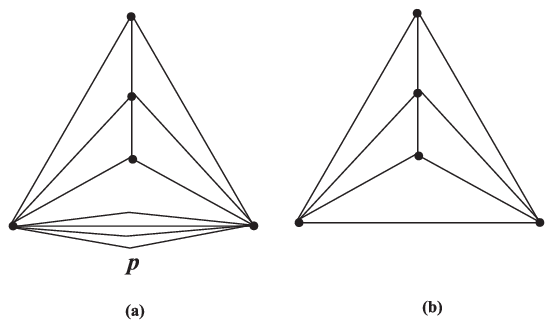,width=3.00in}} \vspace*{8pt}
\caption{Graphs (a) $G(9^*.p\,0)$; (b) $G''$.\label{f1.81}}
\end{figure}

In order to obtain formula for the Tutte polynomial of the graph
$G(9^*.p\,0)$ (Fig. 81a) we use the relations

$$T(G(9^*.p\,0))-T(G(9^*(p-1)\,0))=y^{p-1}T(G'),$$

\noindent where $G'$ is the graph with the Tutte polynomial

$$T(G')=x+2x^2+x^3+y+4xy+3x^2y+3y^2+5xy^2+4y^3+2xy^3+3y^4+y^5.$$

\noindent Since $G''=9^*.(1)$ is the graph $G''$ (Fig. 81b) with the
Tutte polynomial
$$T(G'')=4x+8x^2+5x^3+x^4+4y+13xy+7x^2y+9y^2+9xy^2+8y^3+2xy^3+4y^4+y^5,$$ the general formula
for the Tutte polynomial of the graphs $G(9^*.p\,0)$ is

$$T(G(9^*.p\,0))={{y^{p-1}-1}\over {y-1}}yT(G')+T(G'').$$

\section{Conclusion}

In conclusion, we provide explicit formulae for the Tutte
polynomials of all families of $KL$s derived from source links with
at most $10$ crossings. The Jones polynomials of all alternating and
non-alternating knots and links which belong to the families derived
from source links up to $10$ crossings can be obtained by
substituting  $x\rightarrow -x$ and $y\rightarrow -{1\over x}$ [Bo, 
ChaShro].

\begin{figure}[th]
\centerline{\psfig{file=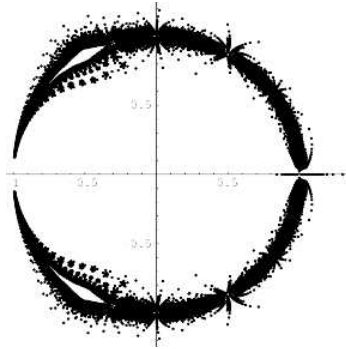,width=1.80in}} \vspace*{8pt}
\caption{Plot of zeros of Jones polynomials for the family $p\,q$
($2\le p\le 50$, $2\le q\le 50$).\label{f1.82}}
\end{figure}

\begin{figure}[th]
\centerline{\psfig{file=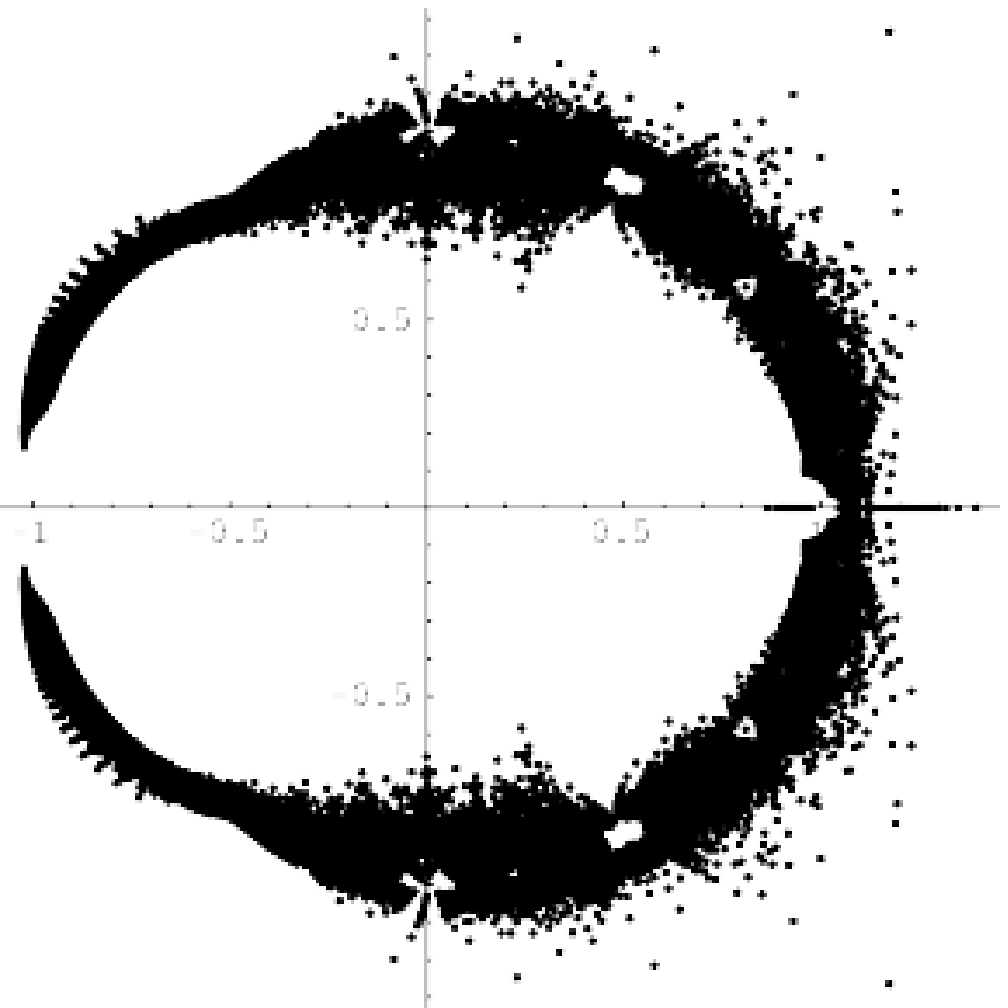,width=1.80in}} \vspace*{8pt}
\caption{Plot of zeros of Jones polynomials for the family of
alternating pretzel $KL$s $p,q,r$ ($2\le p\le 20$, $2\le q\le 20$,
$2\le r\le 20$).\label{f1.83}}
\end{figure}

\begin{figure}[th]
\centerline{\psfig{file=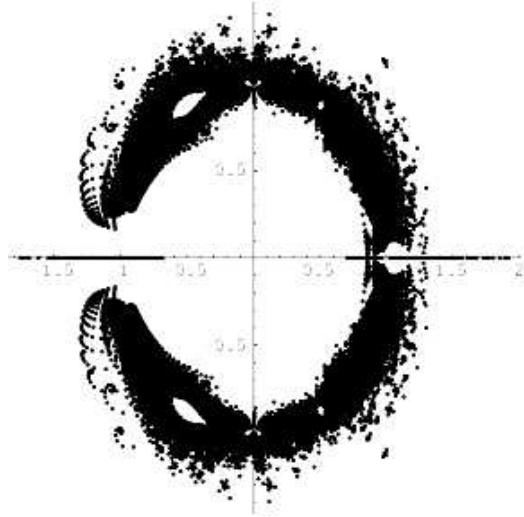,width=2.70in}} \vspace*{8pt}
\caption{Plot of zeros of Jones polynomials for the family of
non-alternating pretzel $KL$s $p,q,-r$ ($2\le p\le 20$, $2\le q\le
20$, $2\le r\le 20$).\label{f1.83}}
\end{figure}

Obtained results can be used for studying the Tutte and Jones
polynomials of $KL$ families. For example, the leading coefficient
of the Tutte polynomial is equal to 1 for all alternating algebraic
$KL$s, and greater then 1 for non-algebraic ones. In particular, it
is invariant for the Tutte polynomials of graphs corresponding to
knots and links inside a family derived from a basic polyhedron.
Additionally they can be used for studying zeros of the Jones
polynomials for the whole family. For example, we plot all zeroes of
the Jones polynomial [ChaShro, JiZha, WuWa] for $KL$ family,
referred to as the characteristic "portrait of family".

The plot of zeroes of Jones polynomials for the alternating link
family $p\,q$ ($2\le p\le 50$, $2\le q\le 50$) is shown in Fig. 82,
the plot of zeroes of Jones polynomials for the alternating link
family $p,q,r$ ($2\le p\le 20$, $2\le q\le 20$, $2\le r\le 20$) is
shown in Fig. 83, and the plot of zeroes of Jones polynomials for
the non-alternating link family $p,q,-r$ ($2\le p\le 20$, $2\le q\le
20$, $2\le r\le 20$) is shown in Fig. 84. More detailed results of
this kind will be given in the forthcoming paper.

The paper is accompanied with the {\it Mathematica} package
downloadable from the address

\medskip

{\tt http://www.mi.sanu.ac.rs/vismath/Tutte.htm}

\medskip

\noindent which provides readers with the possibility to compute all
the results for the given $KL$ families. In particular, the function
{\bf GraphFam} computes the graph corresponding to a $KL$. Functions
{\bf TutteFam} and {\bf JonesFam} compute Tutte and Jones
polynomials of the graphs or links, respectively. Moreover, there
are three functions which provide additional information about zeros
of the Jones polynomial: {\bf Zeros} computes zeros of a particular
Jones polynomial, {\bf ZeroSum} outputs the sum of their absolute
values, and {\bf Portrait} produces the plot of zeros of the Jones
polynomials of a given family. All functions are based on general
formulas, so there are no restrictions on the number of crossings of
$KL$s except for hardware ones.

\bigskip
\bigskip

{\bf References}

\medskip

\medskip

\noindent [Bo] Bollob\' as, B. (1998) {\it Modern Graph Theory}.
Springer, Berlin.

\medskip

\noindent [Cau] Caudron, A. (1982) Classification des n\oe uds et
des enlancements, Public. Math. d'Orsay 82. Univ. Paris Sud, Dept.
Math., Orsay.

\medskip

\noindent [Char] Chartrand, G. (1985) {\it Introductory Graph
Theory}. Dover, New York.

\medskip

\noindent [ChaKo] Champanerkar, A. and Kofman, I. (2005) On the
Mahler measure of Jones polynomials under twisting, Algebraic \@
Geometric Topology {\bf 5}, 1--22.

\medskip

\noindent [ChaShro] Chang, S.~-C. and Shrock, R. (2001) Zeroes of
Jones Polynomials for Families of Knots and Links, Physica A {\bf
296} (2001) 483--494 (arXiv:math-ph/0103043v2).

\medskip

\noindent [Con] Conway, J. (1970) An enumeration of knots and links
and some of their related properties, in {\it Computational Problems
in Abstract Algebra}, Proc. Conf. Oxford 1967 (Ed. J. Leech),
329--358, Pergamon Press, New York.

\medskip

\noindent [ElMe] Ellis-Monaghan, J. and Merino, C. (2008) Graph
Polynomials and their Applications I: The Tutte Polynomial,
arXiv:0803.3079v1 [math].

\medskip

\noindent [Emm] Emmert-Streib, F. (2006) Algorithmic Computation of
Knot Polynomials of Secondary Structure Elements of Proteins,
Journal of Computational Biology, {\bf 13}, 8, 1503--1512.

\medskip

\noindent [Har] Harary, F. (1994) {\it Graph Theory}. Reading, MA:
Addison-Wesley.

\medskip

\noindent [JiZha] Jin, X. and Zhang, F. (2003) Zeroes of the Jones
polynomials for families of pretzel links, Physica A {\bf 328}, 3,
391--408.

\medskip

\noindent [JaSa] Jablan, S.~V. and Sazdanovi\' c, R. (2007) {\it
LinKnot- Knot Theory by Computer}. World Scientific, New Jersey,
London, Singapore (http://math.ict.edu.rs/).

\medskip

\noindent [Kau] Kauffman, L.~H. (1989) A Tutte polynomial for signed
graphs, Discrete Applied Mathematics, {\bf 25}, 105--127.

\medskip

\noindent [KuMu] Kurpita, B. and Murasugi, K. (1998) Knots and
Graphs, Chaos, Solitons \& Fractals, {\bf 9}, 4/5, 623--643.

\medskip

\noindent [Li] Lin, X-S., Zeroes of Jones polynomial,
http://math.ucr.edu/$\sim $xl/abs-jk.pdf

\medskip

\noindent [Ro] Rolfsen, D. (1976) Knots and Links, Publish \& Perish
Inc., Berkeley (American Mathematical Society, AMS Chelsea
Publishing, 2003).

\medskip

\noindent [Thi] Thistlethwaite, M. (1987) A spanning tree expansion
of the Jones polynomial, Topology {\bf 26}, 297--309.

\medskip

\noindent [ThoPe] Thompson, B., Pearce, D.J. (2008) Visualising the
Tutte Polynomial Computation, http://homepages.mcs.vuw.ac.nz/$\sim
$djp/files/TPH-SIENZ07.pdf

\medskip

\noindent [WuWa] Wu, F.Y. and Wang, J. (2001) Zeroes of Jones
polynomial, Physica A {\bf 269}, 483--494 (arXiv:cond-mat/0105013v2
[cond-mat.stat-mech]).

\end{document}